\documentclass[12pt,reqno]{amsart}
\usepackage{cases}
\usepackage{amscd}
\usepackage{amsfonts}
\usepackage{amssymb}
\usepackage{amsmath}
\usepackage{graphicx}
\usepackage{epstopdf}
\usepackage{subfigure}
\usepackage{amsthm}
\usepackage{algorithm} 
\usepackage{algorithmic} 
\usepackage{multirow} 
\usepackage{stmaryrd}
\usepackage{caption} 
\theoremstyle{remark}
\newtheorem{example}{\textbf{Example}}[section]
\numberwithin{equation}{section}

\usepackage{color}
\usepackage{datetime}
\usepackage{tikz}

\makeatletter
\newcommand\figcaption{\def\@captype{figure}\caption}
\newcommand\tabcaption{\def\@captype{table}\caption}
\makeatother

\makeatletter
\newenvironment{breakablealgorithm}
  {
   \begin{center}
     \refstepcounter{algorithm}
     \hrule height.8pt depth0pt \kern2pt
     \renewcommand{\caption}[2][\relax]{
       {\raggedright\textbf{\ALG@name~\thealgorithm} ##2\par}%
       \ifx\relax##1\relax 
         \addcontentsline{loa}{algorithm}{\protect\numberline{\thealgorithm}##2}%
       \else 
         \addcontentsline{loa}{algorithm}{\protect\numberline{\thealgorithm}##1}%
       \fi
       \kern2pt\hrule\kern2pt
     }
  }{
     \kern2pt\hrule\relax
   \end{center}
  }
\makeatother
\usepackage{geometry}
\usepackage{algorithm}
\usepackage{multirow}

\def\bq{\begin{equation}}
\def\eq{\end{equation}}
\def\bqs{\begin{equation*}}
\def\eqs{\end{equation*}}
\def\bsqs{\begin{subequations}}
\def\esqs{\end{subequations}}
\def\ba{\begin{aligned}}
\def\ea{\end{aligned}}
\def\br{\begin{eqnarray}}
\def\er{\end{eqnarray}}
\def\brr{\bq\begin{array}{rlll}}
\def\err{\end{array}\eq}

\def\text#1{\hbox{#1}}
\newtheorem{thm}{Theorem}[section]
\newtheorem{lem}[thm]{Lemma}
\newtheorem{coro}{Corollary}[section]
\newtheorem{rem}[thm]{Remark}

\newcommand{\bsub}{\begin{subequations}}
\newcommand{\esub}{\end{subequations}$\!$}

\makeatletter
\newcommand{\definetitlefootnote}[1]{%
  \newcommand\addtitlefootnote{%
    \makebox[0pt][l]{$^{\bigstar}$}%
    \footnote{\protect\@titlefootnotetext}
  }%
  \newcommand\@titlefootnotetext{\spaceskip=\z@skip $^{\bigstar}$#1}%
}
\makeatother

\title[An AFEM for elliptic equations with line Dirac sources]{An adaptive finite element method for two-dimensional elliptic equations with line Dirac sources \addtitlefootnote}

\author[H.~Cao, H.~Li, N.~Yi, P. Yin]{Huihui Cao$^{\dag}$,  Hengguang Li$^\ddag$, Nianyu Yi$^\dag$, Peimeng Yin$^{*}$}

\address{$\dag$ Hunan Key Laboratory for Computation and Simulation in Science and Engineering, School of Mathematics and Computational Science, Xiangtan University, Xiangtan 411105, Hunan, P.R.China} \email{201721511145@smail.xtu.edu.cn (H. Cao);\  yinianyu@xtu.edu.cn (N. Yi).}
\address{$^\ddag$ Wayne State University, Department of Mathematics, Detroit, Michigan 48202, USA.} \email{li@wayne.edu.}
\address{$^\S$ Multiscale Methods and Dynamics Group, Computer Science and Mathematics Division, Oak Ridge National Laboratory, Oak Ridge, Tennessee 37831, USA.
}\email{yinp@ornl.gov}


\thanks{$^*$ Corresponding author.}

\definetitlefootnote{This manuscript has been authored in part by UT-Battelle, LLC, under contract DE-AC05-00OR22725 with the US Department of Energy (DOE). The US government retains and the publisher, by accepting the article for publication, acknowledges that the US government retains a nonexclusive, paid-up, irrevocable, worldwide license to publish or reproduce the published form of this manuscript, or allow others to do so, for US government purposes. DOE will provide public access to these results of federally sponsored research in accordance with the DOE Public Access Plan (http://energy.gov/downloads/doe-public-access-plan).}

\begin{document}

\begin{abstract} 
In this paper, we propose a novel adaptive finite element method for an elliptic equation with line Dirac delta functions as a source term. We first study the well-posedness and global regularity of the solution in the whole domain.
Instead of regularizing the singular source term and using the classical residual-based a posteriori error estimator, we propose a novel a posteriori estimator based on an equivalent transmission problem with zero source term and nonzero flux jumps on line fractures. The transmission problem is defined in the same domain as the original problem excluding on line fractures, and the solution is therefore shown to be more regular. The estimator relies on meshes conforming to the line fractures and its edge jump residual essentially uses the flux jumps of the transmission problem on line fractures. The error estimator is proven to be both reliable and efficient, an adaptive finite element algorithm is proposed based on the error estimator and the bisection refinement method. Numerical tests show that quasi-optimal convergence rates are achieved even for high order approximations and the adaptive meshes are only locally refined at singular points.
\end{abstract}

\maketitle

\bigskip



\section{Introduction}

We are interested in the adaptive finite element method for the elliptic boundary value problem
\begin{equation}
\label{eq:Possion}
-\Delta u   = \sum_{l=1}^N g_l\delta_{\gamma_l}       \quad \text{in }  \Omega,   
\qquad u =0        \quad \text{on }  \partial \Omega,
\end{equation}
where $\Omega \subset \mathbb{R}^2$ is a polygonal domain, $\gamma_l$, $l=1,\cdots, N$ are disjoint or intersecting line fractures strictly contained in $\Omega$, $g_l\in H^{\beta_l}(\gamma_l)$ with $\beta_l\geq  0$, and $g_l\delta_{\gamma_l}$ in source term $\sum_{l=1}^N g_l\delta_{\gamma_l}$ is a line Dirac measure on line fracture $\gamma_l$ satisfying
\bq\label{deltadef}
\langle g_l\delta_{\gamma_l}, v \rangle = \int_{\gamma_l} g_l(s) v(s)ds, \qquad \forall\ v|_{\gamma_l} \in L^2(\gamma_l).
\eq
Although $g_l \in H^{\beta_l}(\gamma_l) \subset L^2(\gamma_l)$, the line Dirac measure $\sum_{l=1}^N g_l\delta_\gamma \not \in L^2(\Omega)$.

The model (\ref{eq:Possion}) has been widely used to describe monophasic flows in porous media, tissue perfusion or drug delivery by a network of blood vessels \cite{DAngelo12}, and it also has applications in elliptic optimal control problems \cite{Gong14}.
The solution of the elliptic problem (\ref{eq:Possion}) is smooth in a large part of the domain, but it becomes singular in the region close to line fractures $\gamma_l$ and the region close to the vertices of the domain \cite{LWYZ21}. 
The corner singularity has been well understood in the literature \cite{Apel99, Dauge88, Grisvard85, Kondratiev67, Li10} and references therein, we shall focus on the regularity of the solution near line fractures $\gamma_l$. The smoothness of the source term can be obtained by the duality argument \cite{LionsMaganesVolI},
thus the regularity of solution for problem (\ref{eq:Possion}) follows from the standard elliptic regularity theory \cite{Grisvard92, Alinhac07}.

Finite element methods for the second-order elliptic equations with singular source terms date back to the 1970s, but the main focus was on point Dirac delta sources (see e.g., \cite{Babuska71, Scott73, Scott76, Casas85, ABR06, HW12}). More recently, singular sources on complex geometry \cite{Gong14, HR19, LWYZ21, HL21, DAngelo08, DAngelo12, Ariche16}, including one-dimensional (1D) fracture sources, have attracted more attention. The finite element method was studied in \cite{HR19} for problems involving a $\mathcal{C}^2$ closed fracture strictly contained in the domain, and later an adaptive finite element method was proposed to improve the convergence rate \cite{HL21}. As a controlled equation in an optimal control problem, the boundary value problem (\ref{eq:Possion}) with a single $\mathcal{C}^2$ curve fracture was solved in \cite{Gong14} by the linear finite element method. 

Due to the lack of regularity, the finite element method for problem (\ref{eq:Possion}) has only a convergence rate $h^{\frac{1}{2}}$ on general quasi-uniform meshes. 
Later on, in order to improve the convergence rate for problem (\ref{eq:Possion}) with one line segment fracture and the coefficient function $g_1=\text{const}$, Li et al. \cite{LWYZ21} studied the regularities in both Sobolev space and weighted Sobolev space, and a finite element algorithm was proposed to approximate the singular solution at the optimal convergence rate on graded meshes, which were densely refined only at the endpoints of the line fractures.
The graded finite element algorithm in \cite{LWYZ21} can be applied to problem (\ref{eq:Possion}), but the grading parameter (used to generate graded meshes) depends on the smoothness of functions $g_l$, and it could be complicated to calculate and may vary case by case for different functions $g_l$ in order to generate graded meshes on which the finite element solutions are optimal.

An alternative way to obtain optimal finite element solutions for problem (\ref{eq:Possion}) is by the adaptive finite element methods (AFEMs), which are effective numerical methods for problems with singularities.  AFEMs usually consist of four steps (see e.g., \cite{D96, MNS02}),
$$
\text{SOLVE $\rightarrow$ ESTIMATE $\rightarrow$ MARK $\rightarrow$ REFINE,}
$$
which generates a sequence of meshes, on which the finite element approximations converge to the solution of the target problem.
An essential ingredient of the AFEMs is a posteriori error estimator, which is a computable quantity that depends on the finite element approximation and known data, and provides information about the size and the distribution of the error of the numerical approximation. Therefore, it can be used to guide mesh adaption and as an error estimation.
For results on the a posteriori error estimations of finite element analysis for the second order elliptic problems with an $L^2$ source term can be found in \cite{AO00, V96} and references therein.

Elliptic problems with point Dirac delta source term were sufficiently studied by the AFEMs, for which the residual-based a posteriori error estimators were widely employed to guide the mesh adaptions and as the finite element solution error estimations \cite{BDD2004,D96,MNS00,S07,V96}.
Due to the singularity of the point Dirac delta source term,  it was generally regularized to an $L^2(\Omega)$ or $L^p(\Omega)$ function with $1<p<\infty$ by projecting the source term to a polynomial space.
Therefore, the residual-based a posteriori error estimator for the Poisson problem with an $L^2$ source term \cite{V1994,AO00} can be applied.

Recently, the regularization techniques of projecting the source term to an $L^2(\Omega)$ or $L^p(\Omega)$ function were also applied to the elliptic problem with line Dirac delta source term \cite{HL21, MMR22}.
The resulted residual-based a posteriori error estimators were also effective in proposing adaptive finite element algorithms, which don't rely on specific meshes.
However, the associated adaptive finite element solutions involve not only the discretization error but also the regularization error \cite{HL21, MMR22}, and the error estimators might lead to over-refinement on adaptive meshes or low convergence rates for high order approximations. 

Motivated by the performance of the finite element solutions on graded meshes for which the grading parameters are involved, and of AFEMs based on regularized source terms for which the meshes are generally over-refined even for low order approximations, in this work we propose a novel residual-based a posteriori error estimator, which is of high order convergence rates and the adaptive meshes are only locally refined near the singularities of the solution.

Instead of regularizing the singular line Dirac source term in problem (\ref{eq:Possion}), we transfer the problem (\ref{eq:Possion}) to an equivalent interface problem with zero source term and nonzero flux jumps on line fractures $\gamma_l$. More specifically, the coefficients $g_l$ in the line Dirac source term are transferred to the flux jumps on line fractures. The new transferred problem is known as the transmission problem \cite{Li10}, which is defined in the same domain as the original problem excluding on line fractures.
The solution of problem (\ref{eq:Possion}) excluding on the line fractures solves the transmission problem, and it is shown that the solution becomes more regular after the transmission, which implies the finite element solutions for problem (\ref{eq:Possion}) would have a higher convergence rate if the meshes conform to the line fractures. Compared with the convergence rate on general quasi-uniform meshes, the finite element method for problem (\ref{eq:Possion}) has a better convergence rate $h^{\min\{\alpha,\beta+\frac{1}{2},1\}}$ on conforming quasi-uniform meshes, where $\beta=\min_l\{\beta_l\}$ and $\alpha<\frac{\pi}{\omega}$ with $\omega$ being the largest interior angle of the polygonal domain $\Omega$.

Our residual-based a posteriori error estimator is proposed based on the transmission problem.
First, we triangulate the mesh conforming to line fractures $\gamma_{l}$, namely, $\gamma_{l}$ is the union of some edges in the triangulation. Second, the error estimator consists of element residual with zero source and edge residuals involving the difference with the flux jumps $g_l$ on line fractures.
We derive the reliability and efficiency of the proposed a posteriori error estimator with novel skills in handling the edge residual. 
Based on the derived error estimator and bisection mesh refinement method, we propose an adaptive finite element algorithm. The quasi-optimal convergence rates can be numerically achieved for finite element approximations with the adaptive meshes only locally refined at the singular points.

As far as we have known, this is the first work using the transmission problem to construct a posteriori error estimator for problems with Dirac source terms.
It would be interesting to apply the proposed AFEM to problem (\ref{eq:Possion}) with curved line segments, and to explore the applications in three-dimension, we will leave these topics to our future work.
The rest of the paper is organized as follows. In Section \ref{sec-2}, we discuss the well-posedness and global regularity of equation \eqref{eq:Possion} in Sobolev spaces. In Section \ref{31}, we introduced the transmission problem associated with problem (\ref{eq:Possion}) and investigate its well-posedness and regularity, and also showed its relationship with problem (\ref{eq:Possion}).  In Section \ref{sec-3}, we propose a novel residual-based a posteriori error estimator, show its reliability and efficiency, and propose an adaptive finite element algorithm. 
In Section \ref{sec-4}, we present various numerical test results to validate the theoretical findings.

Throughout this paper, $C > 0$ denotes a generic constant that may be different at different occurrences. It will depend on the computational domain, but not on the functions involved and mesh parameters.

\section{Well-posedness and regularity in Sobolev spaces} \label{sec-2} 

Denote by $H^m(\Omega)$, $m\geq 0$,  the Sobolev space that consists of functions whose $i$th ($0\leq i\leq m$) derivatives are square integrable. 
Denote by $H^1_0(\Omega)\subset H^1(\Omega)$  the subspace consisting of functions with  zero trace on the boundary  $\partial\Omega$.
For $s>0$, let $s=m+t$, where $m\in \mathbb Z_{\geq 0}$ and $0<t<1$.
Recall that for   $D\subseteq \mathbb{R}^d$,
the fractional order Sobolev space $H^s(D)$ consists of distributions $v$ in $D$ satisfying
$$
\|v\|^2_{H^s(D)}:=\|v\|^2_{H^m(D)} + \sum_{|\alpha|= m}\int_{D}\int_{D} \frac{|\partial^\alpha v(x) - \partial^\alpha v(y)|^2 }{|x-y|^{d+2t}} dxdy <\infty,
$$
where $\alpha=(\alpha_1, \ldots, \alpha_d)\in\mathbb Z^d_{\geq0}$ is a multi-index such that $\partial^\alpha=\partial_{x_1}^{\alpha_1}\cdots\partial^{\alpha_d}_{x_d}$ and $|\alpha|=\sum_{i=1}^d\alpha_i$.
We denote by $H_0^s(D)$ the closure of $C_0^\infty(D)$ in $H^s(D)$, and $H^{-s}(D)$ the dual space of $H_0^s(D)$.
Let $\widetilde{H}^{s}(D)$ be the space of all $v$ defined in $D$ such that $\tilde{v}\in {H}^{s}(\mathbb{R}^d)$, where $\tilde{v}$ is the extension of $v$ by zero outside $D$. 


\subsection{Trace estimates}
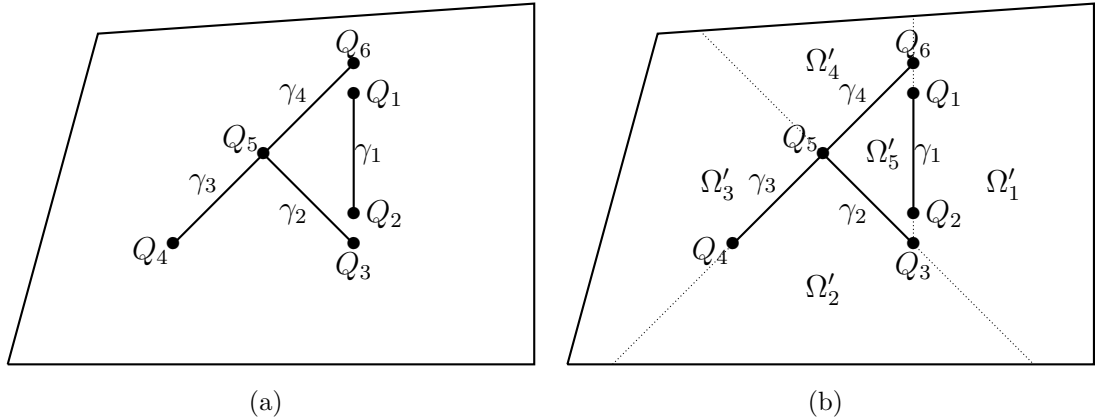
\begin{figure}
\centering
\subfigure[]{
\begin{tikzpicture}[scale=0.2]
\draw[thick]
(-19,-12) -- (-13,10) -- (16,12) -- (16,-12) -- (-19,-12);

\draw[thick] (-8,-4) -- (4,8);
\draw[thick] (-2,2) -- (4,-4);
\draw[thick] (4,6) -- (4,-2);

\draw[ultra thick]  (-8,-4) node {$\bullet$};
\draw[ultra thick]  (4,-4) node {$\bullet$};
\draw[ultra thick]  (4,8) node {$\bullet$};
\draw[ultra thick]  (-2,2) node {$\bullet$};
\draw[ultra thick]  (4,6) node {$\bullet$};
\draw[ultra thick]  (4,-2) node {$\bullet$};

\draw (-6,0) node {$\gamma_3$};
\draw (0,6) node {$\gamma_4$};
\draw (5,2) node {$\gamma_1$};
\draw (0,-2) node {$\gamma_2$};

\draw (-9.5,-4.5) node {$Q_4$};
\draw (4,-5.5) node {$Q_3$};
\draw (4,9.3) node {$Q_6$};
\draw (-3.5,3) node {$Q_5$};
\draw (8,6) node[anchor = east] {$Q_1$};
\draw (8,-2) node[anchor = east] {$Q_2$};
\end{tikzpicture}
}
\subfigure[ ]{
\begin{tikzpicture}[scale=0.2]
\draw[thick]
(-19,-12) -- (-13,10) -- (16,12) -- (16,-12) -- (-19,-12);

\draw[thick] (-8,-4) -- (4,8);
\draw[thick] (-2,2) -- (4,-4);
\draw[thick] (4,6) -- (4,-2);

\draw[densely dotted] (-16,-12) -- (-8,-4);
\draw[densely dotted] (-10,10) -- (-2,2);
\draw[densely dotted] (4,-4) -- (12,-12);
\draw[densely dotted] (4,6) -- (4,144/13);
\draw[densely dotted] (4,-2) -- (4,-4);

\draw[ultra thick]  (-8,-4) node {$\bullet$};
\draw[ultra thick]  (4,-4) node {$\bullet$};
\draw[ultra thick]  (4,8) node {$\bullet$};
\draw[ultra thick]  (-2,2) node {$\bullet$};
\draw[ultra thick]  (4,6) node {$\bullet$};
\draw[ultra thick]  (4,-2) node {$\bullet$};

\draw (-6,0) node {$\gamma_3$};
\draw (0,6) node {$\gamma_4$};
\draw (5,2) node {$\gamma_1$};
\draw (0,-2) node {$\gamma_2$};

\draw (-9.5,-4.5) node {$Q_4$};
\draw (4,-5.5) node {$Q_3$};
\draw (4,9.3) node {$Q_6$};
\draw (-3.5,3) node {$Q_5$};
\draw (8,6) node[anchor = east] {$Q_1$};
\draw (8,-2) node[anchor = east] {$Q_2$};

\draw (-2,-7) node {$\Omega_2^\prime$};
\draw (-9,0) node {$\Omega_3^\prime$};
\draw (-2,8) node {$\Omega_4^\prime$};
\draw (10,0) node {$\Omega_1^\prime$};
\draw (2,2) node {$\Omega_5^\prime$};
\end{tikzpicture}
}
\vspace*{-15pt}
    \caption{(a) Domain $\Omega$ containing four line fractures $\gamma_1,\,\gamma_2,\,\gamma_3$ and $\gamma_4$. (b) $\Omega$ is decomposed into five sub-domains $\{\Omega_j^\prime\}_{j = 1}^5$ by $\gamma_i,\,i = 1,\,2,\,3,\,4$ .}
    \label{fig:Omega2}
\end{figure}

A sketch drawing of the domain $\Omega$ with several line fractures is given in Figure \ref{fig:Omega2}(a). To obtain the trace estimates on line fractures, we first introduce the trace estimate on a general polygonal domain with no line fracture.

\begin{lem} \cite{Ding96, M20} \label{trace1}
Let $\Omega'$ be a polygonal domain with no line fracture, then the trace operator
$$
\vartheta \ : \ H^s(\Omega') \ \rightarrow \ H^{s-\frac{1}{2}}(\partial\Omega')
$$
is bounded for $\frac{1}{2}<s<\frac{3}{2}$.
\end{lem}

\begin{lem}\label{trace}
For the domain $\Omega$ with line segment fractures $\gamma_l$, $l=1,\cdots, N$, it follows that the trace operator
$$
\vartheta \ : \ H^s(\Omega) \ \rightarrow \ H^{s-\frac{1}{2}}(\cup_{l=1}^N\gamma_l)
$$
is bounded for $\frac{1}{2}<s<\frac{3}{2}$.
\end{lem}
\begin{proof}
By extending line fractures $\gamma_l$ appropriately to the boundary of the domain $\Omega$ or another line fracture and denoting the extended line fractures by $\gamma'_l$, which partition the domain $\Omega$ into $M$ polygonal subdomains $\Omega'_j$, $1\leq j \leq M$ and $\gamma'_l$ is shared by neighboring subdomains $\Omega'_j$ (see Figure \ref{fig:Omega2} (b)).
For any $v \in H^s(\Omega)$, it follows
$$
v \in H^s(\Omega'_j), \quad j=1,\cdots, M,
$$
satisfying
$$
\|v\|^2_{H^s(\Omega)} = \sum_{j=1}^{M} \|v\|^2_{H^s(\Omega'_j)} = \|v\|^2_{H^s(\cup_{j=1}^M \Omega'_j)}.
$$
By Lemma \ref{trace1}, if $\frac{1}{2}<s<\frac{3}{2}$, it follows for $l=1,\cdots, N$,
$$
\|v\|_{H^{s-\frac{1}{2}}(\gamma_l)} \leq \|v\|_{H^{s-\frac{1}{2}}(\gamma'_l)} \leq C \|v\|_{H^s(\cup_{j=1}^M \Omega'_j)} = C \|v\|_{H^s(\Omega)}.
$$
Therefore, the conclusion holds.
\end{proof}

\subsection{Well-posedness and regularity}

We have the following result regarding the line Dirac measure $\sum_{l=1}^N g_l\delta_{\gamma_l}$.
\begin{lem}
\label{lemma2-1}
For $\epsilon>0$, the line Dirac measure $\sum_{l=1}^N g_l\delta_{\gamma_l} \in H^{-\frac{1}{2}- \epsilon}(\Omega)$ satisfying
$$
\left\|\sum_{l=1}^N g_l\delta_{\gamma_l}\right\|_{H^{-\frac{1}{2}-\epsilon}(\Omega)} \leq C \sum_{l=1}^N \|g_l\|_{L^2(\gamma_l)}.
$$
\end{lem}

\begin{proof}
The proof is based on the duality pairing (e.g., \cite{LionsMaganesVolI}). For $v \in H^{\frac{1}{2}+\epsilon}(\Omega)$,  by H\"older's inequality and Lemma \ref{trace}, we have for $l=1, \cdots, N$,
\bqs
\ba
\langle g_l\delta_{\gamma_l}, v \rangle = &\int_{\gamma_l} g_l(s)v(s)ds
\leq C \|g_l\|_{L^2(\gamma_l)}\|v\|_{L^2(\gamma_l)}
\leq C \|g_l\|_{L^2(\gamma_l)}\|v\|_{H^{\epsilon}(\gamma_l)}
\leq C \|g_l\|_{L^2(\gamma_l)}\|v\|_{H^{\frac{1}{2}+\epsilon}(\Omega)}.
\ea
\eqs
Therefore, by definition, we have
\[
\left\|\sum_{l=1}^N g_l\delta_{\gamma_l}\right\|_{H^{-\frac{1}{2}-\epsilon}(\Omega)} := \sup \left\{\left\langle \sum_{l=1}^N g_l\delta_{\gamma_l}, v \right\rangle \ : \  \|v\|_{H^{\frac{1}{2}+\epsilon}} = 1 \right\} \leq C \sum_{l=1}^N \|g_l\|_{L^2(\gamma_l)}.
\]
\end{proof}

The variational formulation for problem (\ref{eq:Possion}) is to find $u\in H_0^1(\Omega)$, such that
\begin{eqnarray}\label{eqn.weak}
a(u, v):=\int_\Omega\nabla u\cdot \nabla v dx=\left\langle \sum_{l=1}^N g_l\delta_{\gamma_l}, v \right\rangle, \quad \forall\ v\in H^1_0(\Omega).
\end{eqnarray}
By Lemma \ref{lemma2-1}, the variational formulation (\ref{eqn.weak}) is well-posed.


Therefore, we have the following global regularity estimate.
\begin{lem}\label{thm2-2}
For $\epsilon>0$, the elliptic boundary value problem \eqref{eq:Possion} admits a unique solution $u \in H^{\frac{3}{2}- \epsilon}(\Omega) \cap H^1_0(\Omega)$ satisfying
\bq
\|u\|_{H^{\frac{3}{2}- \epsilon}(\Omega)} \leq C \sum_{l=1}^N \|g_l\|_{L^2(\gamma_l)}.
\eq
\end{lem}

\begin{proof}
The gives
$$
\|u\|_{H^{\frac{3}{2}- \epsilon}(\Omega)}  \leq C \left\|\sum_{l=1}^N g_l\delta_{\gamma_l}\right\|_{H^{-\frac{1}{2}-\epsilon}(\Omega)} \leq C \sum_{l=1}^N \|g_l\|_{L^2(\gamma_l)}.
$$
\end{proof}

\begin{rem}\label{superposition}
Since problem (\ref{eq:Possion}) is a linear problem, so that the solution $u$ of problem (\ref{eq:Possion}) 
can be obtained by summing of solutions of the following problems with one line Dirac source term for $l=1,\cdots, N$,
\begin{equation}
\label{eq:Possion2}
-\Delta u_l   = g_l\delta_{\gamma_l}      \quad \text{in }  \Omega,    \qquad\qquad
u_l           =0        \quad \text{on }  \partial \Omega.
\end{equation}
By the superposition principle, one has
$$
u = \sum_{l=1}^N u_l.
$$
The estimate in Lemma \ref{thm2-2} can also be obtained by first obtaining the estimates for problem (\ref{eq:Possion2}), and then taking the summation of all these estimates.
\end{rem}

Based on Lemma \ref{thm2-2}, we find that no matter how smooth the functions $g_l$ are, the solution of problem (\ref{eq:Possion}) is merely in $H^{\frac{3}{2}- \epsilon}(\Omega)$ for $\forall \epsilon>0$ due to the appearance of the singular line Dirac measure $\sum_{l=1}^N g_l\delta_{\gamma_l}$ in the source term.
Then, by Lemma \ref{thm2-2} and the Sobolev imbedding Theorem \cite{M20}, we have the following result.
\begin{coro}\label{co1}
For $\epsilon>0$, the solution $u$ of problem \eqref{eq:Possion} is H\"{o}lder continuous $u\in \mathcal C^{0,1/2-\epsilon}(\Omega)$. In particular, the solution $u\in \mathcal C^0({\Omega})$.
\end{coro}

By Corollary \ref{co1}, we know that the solution of problem (\ref{eq:Possion}) is continuous across line fractures $\gamma_l$, $l=1, \cdots, N$. 
Next, we introduce the transmission problem of problem (\ref{eq:Possion}) to investigate the normal derivatives of $u$ across line fractures.

\section{The transmission problem}\label{31}

Let $\mathbf{n}^{\pm}$ be the outward unit normal of the region on each side of the fracture $\gamma_l$. For a function $v$, we denote $v^\pm$ (resp. $\partial_{\mathbf{n}^{\pm}} v^\pm$) the traces of $v$ (resp. $\nabla v$) evaluated on the fracture $\gamma_l$ from the region on each side.
We define the jump of $v$ across $\gamma_l$ by $[v]=v^+-v^-$ and the jump of its normal derivatives (or flux jumps) on $\gamma_l$ by $[\partial_{\mathbf {n}} v]=\mathbf {n}^+ \cdot \nabla v^++\mathbf {n}^- \cdot \nabla v^-$.

Based on the observation of the solution and weak solution of problem (\ref{eq:Possion}), we introduce the following interface problem,
\begin{subequations}\label{eq:2d}
\begin{align}
-\Delta w=0 &\quad \mbox{ in   }\Omega \setminus \cup_{l=1}^N\gamma_l,\\
[w]=0 &\quad \mbox{ on   } \gamma_l, \ l=1,\cdots, N,\\
[\partial_\mathbf{n}w]= g_l &\quad \mbox{ on   } \gamma_l, \ l=1,\cdots, N,\\
w=0 &\quad \mbox{ on   } \partial \Omega.
\end{align}
\end{subequations}
The interface problem (\ref{eq:2d}) is known as the transmission problem of the elliptic problem (\ref{eq:Possion}) \cite{Li10}.

We define a space
$$
V=\left \{v\in H^1(\Omega\setminus \cup_{l=1}^N\gamma_l): \ v|_{\partial\Omega}=0, \ [v]|_{\gamma_l}=0, \quad l=1,\cdots, N \right \}.
$$
Similar to \cite{Bramble1996}, multiplying a test function $v \in H_0^1(\Omega)$ on both sides of (\ref{eq:2d}a), and applying the Green's formula together with the interface and boundary conditions (\ref{eq:2d}b-d), we have
\bqs
\ba
-\int_{\Omega\setminus\cup_{l=1}^N\gamma_l} \Delta w v dx = & \int_{\Omega\setminus \cup_{l=1}^N\gamma_l} \nabla w \cdot \nabla v dx - \sum_{l=1}^N\int_{\gamma_l} [\partial_{\mathbf{n}}w]  v ds =0,
\ea
\eqs
thus the variational formulation for the transmission problem (\ref{eq:2d}) is to find $w \in V$ such that
\bq\label{weaktrans}
\int_{\Omega\setminus \cup_{l=1}^N\gamma_l} \nabla w \cdot \nabla v dx = \sum_{l=1}^N\int_{\gamma_l} g_l v ds,\quad \forall v \in H_0^1(\Omega).
\eq
\begin{lem}\label{transweaklem}
The weak formulation (\ref{weaktrans}) admits a unique solution
$w\in V$ satisfying
\bq\label{wurelate}
w=u|_{\Omega\setminus \cup_{l=1}^N\gamma_l},
\eq
where $u$ is the solution of the weak formulation (\ref{eqn.weak}).
\end{lem}
\begin{proof}
By Lemma \ref{trace} and (\ref{deltadef}), we have for $\forall v \in H_0^1(\Omega)$,
$$
\left\langle\sum_{l=1}^N g_l\delta_{\gamma_l}, v \right\rangle = \sum_{l=1}^N\int_{\gamma_l} g_l v ds.
$$
Since $u$ is the solution of (\ref{eqn.weak}), so $w$ given in (\ref{wurelate}) satisfying $w \in V$ and solves (\ref{weaktrans}), whose existence and uniqueness therefore follow from the well-posedness of the weak formulation (\ref{eqn.weak}).
\end{proof}

Lemma \ref{transweaklem} indicates that the solution of problem (\ref{eq:Possion})
solves the transmission problem (\ref{eq:2d}) at least in $H^{\frac{3}{2}-\epsilon}(\Omega\setminus \cup_{l=1}^N\gamma_l) \cap V$.

To investigate the regularity of the transmission problem (\ref{eq:2d}), we first consider the following interface problem,
\begin{subequations}\label{eq:2d0}
\begin{align}
-\Delta z=0 &\quad \mbox{ in   }\Omega \setminus \Gamma_0,\\
[z]=0 &\quad \mbox{ on   } \Gamma_0,\\
[\partial_\mathbf{n}z]= g &\quad \mbox{ on   } \Gamma_0,\\
z=0 &\quad \mbox{ on   } \partial \Omega,
\end{align}
\end{subequations}
where $\Gamma_0$ is a closed sufficiently smooth curve strictly contained in $\Omega$, and $g \in H^{\beta}(\Gamma_0)$ with $\beta\geq 0$. 
For problem (\ref{eq:2d0}), we recall the following result from \cite{Grisvard92, CZ98, ABGL20}.
\begin{lem}\label{regzext}
Let $z$ be the solution of the problem (\ref{eq:2d0}), then it follows $z \in H^{\min\{1+\alpha, \beta+\frac{3}{2}\}}(\Omega\setminus\Gamma_0)$ satisfying
\bq
\|z\|_{H^{\min\{1+\alpha, \beta+\frac{3}{2}\}}(\Omega\setminus\Gamma_0)} \leq C \|g\|_{H^{\beta}(\Gamma_0)},
\eq
where $\alpha<\frac{\pi}{\omega}$ with $\omega$  the largest interior angle of the polygonal domain $\Omega$.
\end{lem}

Next, we introduce the following result from \cite[Theorem 1.2.15 and Theorem 1.2.16]{Grisvard92}.
\begin{lem}\label{gweight}
For a point $x\in \gamma_l$, let $\rho_l(x)$ be multiplication of the distances of $x$ to the endpoints of $\gamma_l$. Then one has $\frac{v}{\rho_l^s} \in L^2(\gamma_l)$ for all $v \in H^{s}(\gamma_l)$ when $s \in (0,\frac{1}{2})$, and 
one also has $\frac{D^\nu v}{\rho_l^{s-|\nu|}}\in L^2(\Omega)$ for all $v \in H_0^s(\Omega)$ and $|\nu|\leq s$ provided $s-\frac{1}{2}$ is not an integer.
\end{lem}

\begin{lem}\label{gextension}
For any $\epsilon>0$, we have the following results,\\
(i) if $g_l\in H^{\beta_l}(\gamma_l)$, it follows $g_l \in \widetilde{H}^{\min\{\beta_l, \frac{1}{2}-\epsilon\}}(\gamma_l)$;\\
(ii) if $g_l\in H_0^{\beta_l}(\gamma_l)$ and $\beta_l+\frac{1}{2}>0$ is not an integer, it follows $g_l \in \widetilde{H}^{\beta_l}(\gamma_l)$;\\
(iii) if $g_l\in H_0^{\beta_l}(\gamma_l)$ and $\beta_l+\frac{1}{2}>0$ is an integer, it follows $g_l \in \widetilde{H}^{\beta_l-\epsilon}(\gamma_l)$.
\end{lem}
\begin{proof}
The proofs of (i) and (ii) follow from Lemma \ref{gweight} and the definition of the fractional space $\widetilde{H}^{s}(D)$ defined at the beginning of Section \ref{sec-2}. For (iii), we have $g_l \in H_0^{\beta_l}(\gamma_l) \subset H_0^{\beta_l-\epsilon}(\gamma_l)$, then the conclusion holds by applying Lemma \ref{gweight} and the definition of $\widetilde{H}^{s}(D)$.
\end{proof}


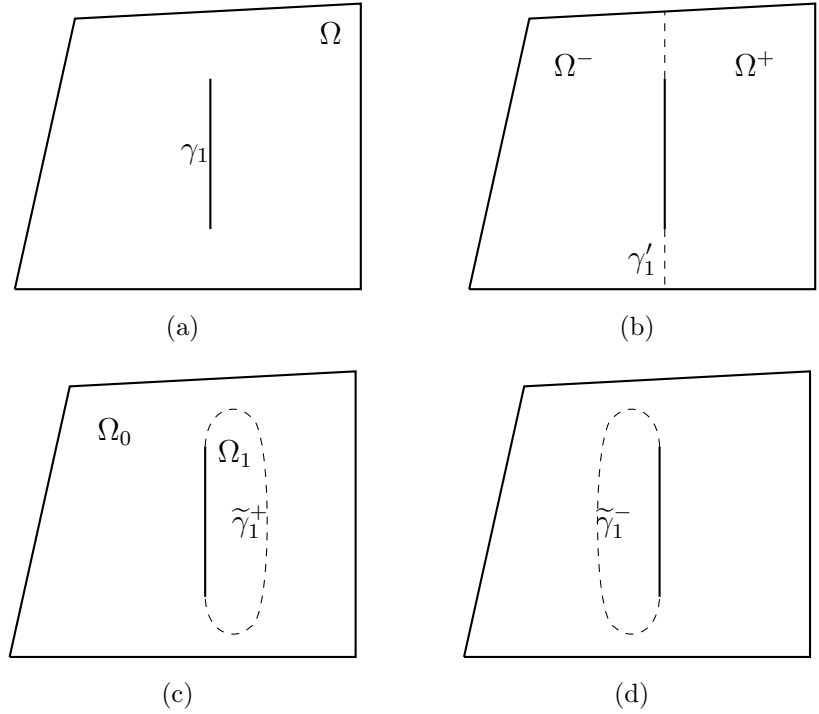
\begin{figure}
\centering
\subfigure[]{
\begin{tikzpicture}[scale=0.2]
\draw[thick]
(-13,-9) -- (10,-9) -- (10,10) -- (-9,9) -- (-13,-9);
\draw[thick] (0,-5) -- (0,5);
\draw (-1,0) node {$\gamma_1$};
\draw (8,8) node {$\Omega$};
\end{tikzpicture}
}\hspace{10mm}
\subfigure[]{
\begin{tikzpicture}[scale=0.2]
\draw[thick]
(-13,-9) -- (10,-9) -- (10,10) -- (-9,9) -- (-13,-9);
\draw[thick] (0,-5) -- (0,5);
\draw[dashed] (0,-5) -- (0,-9);
\draw[dashed] (0,5) -- (0,180/19);
\draw (-1.5,-7) node {$\gamma'_1$};
\draw (6,6) node {$\Omega^+$};
\draw (-6,6) node {$\Omega^-$};
\end{tikzpicture}
}\\
\subfigure[ ]{
\begin{tikzpicture}[scale=0.2]
\draw[thick]
(-13,-9) -- (10,-9) -- (10,10) -- (-9,9) -- (-13,-9);
\draw[thick] (0,-5) -- (0,5);
\draw (3,0) node {$\widetilde{\gamma}_1^+$};
\draw[dashed] (0,5)
.. controls (0.1,7.5) and (2,8) .. (3,7)
.. controls (4.5, 6.5) and (4.5,-6.5) .. (3,-7)
.. controls (2, -8) and (0.1,-7.5) .. (0,-5);
\draw (2,4.5) node {$\Omega_1$};
\draw (-6,6) node {$\Omega_0$};
\end{tikzpicture}
}\hspace{10mm}
\subfigure[ ]{
\begin{tikzpicture}[scale=0.2]
\draw[thick]
(-13,-9) -- (10,-9) -- (10,10) -- (-9,9) -- (-13,-9);
\draw[thick] (0,-5) -- (0,5);
\draw (-3,0) node {$\widetilde{\gamma}_1^-$};
\draw[dashed] (0,5)
.. controls (-0.1,7.5) and (-2,8) .. (-3,7)
.. controls (-4.5, 6.5) and (-4.5,-6.5) .. (-3,-7)
.. controls (-2, -8) and (-0.1,-7.5) .. (0,-5);
\end{tikzpicture}
}
\vspace*{-15pt}
    \caption{The extension of the line fracture $\gamma_1$.}
    \label{fig:Omega3}
\end{figure}

\begin{thm}\label{wreglem}
For $\forall \epsilon>0$, let $w$ be the solution of the transmission problem (\ref{eq:2d}), if $g_l \in H^{\beta_l}(\gamma_l)$, $l=1,\cdots, N$, then it follows
\bq\label{transreg}
\|w\|_{H^{\min\{\alpha+1, \beta+\frac{3}{2}, 2-\epsilon\}}(\Omega\setminus \cup_{l=1}^N\gamma_l)} \leq C \sum_{l=1}^N \|g_l\|_{{H}^{\beta_l}(\gamma_l)}.
\eq
Further, if all $g_l \in H_0^{\beta_l}(\gamma_l)$, $l=1,\cdots, N$, it follows
\bq\label{transreg1}
\ba
\|w\|_{H^{\min\{\alpha+1,\beta+\frac{3}{2}\}}(\Omega\setminus \cup_{l=1}^N\gamma_l)} \leq & C \sum_{l=1}^N \|g_l\|_{H^{\beta_l}(\gamma_l)},\quad \text{if $\beta+\frac{1}{2}>0$ is not an integer};\\
\|w\|_{H^{\min\{\alpha+1,\beta+\frac{3}{2}-\epsilon\}}(\Omega\setminus \cup_{l=1}^N\gamma_l)} \leq & C \sum_{l=1}^N \|g_l\|_{H^{\beta_l}(\gamma_l)},\quad \text{if $\beta+\frac{1}{2}>0$ is an integer};\\
\ea
\eq
Here, $\beta=\min_l\{\beta_l\}$ and $\alpha<\frac{\pi}{\omega}$ with $\omega$ being the largest interior angle of the polygonal domain $\Omega$.
\end{thm}
\begin{proof}
We first prove the case with only one line fracture $\gamma_1$ as shown in Figure \ref{fig:Omega3}(a). We extend $\gamma_1$ to $\gamma'_1$ which has two points of intersection with the boundary $\partial \Omega$, then $\Omega$ is partitioned into two open subdomains $\Omega^-$ and $\Omega^+$ (see Figure \ref{fig:Omega3}(b)). In $\Omega^+$, we extend the line fracture $\gamma_1$ to a closed $\mathcal C^2$ curve $\widetilde{\gamma}_1^+$ partitioning $\Omega$ into two subdomains $\Omega_0$ and $\Omega_1$ as shown in Figure \ref{fig:Omega3}(c), and extend $g_1$ on $\gamma_1$ to $\widetilde g_1$ on $\widetilde{\gamma}_1^+$ satisfying
\begin{equation}\label{g_extend}
\widetilde{g}_1=
\left\{
\begin{aligned}
g_1\quad  & \text{on } \gamma_1, \\
0\quad  & \text{on }  \widetilde{\gamma}_1^+\setminus\gamma_1.
\end{aligned}
\right.
\end{equation}
Then the transmission problem (\ref{eq:2d}) is equivalent to the following problem
\begin{subequations}\label{eq:2d+}
\begin{align}
-\Delta w=0 &\quad \mbox{ in   }\Omega \setminus \widetilde{\gamma}^+_1,\\
[w]=0 &\quad \mbox{ on   } \widetilde{\gamma}^+_1,\\
[\partial_\mathbf{n}w]= \widetilde g_1 &\quad \mbox{ on   } \widetilde{\gamma}^+_1,\\
w=0 &\quad \mbox{ on   } \partial \Omega.
\end{align}
\end{subequations}
Note that $\Omega^-\subset \Omega_0$, so by Lemma \ref{regzext} and Lemma \ref{gextension}, if $g_1 \in H^{\beta_1}(\gamma_1)$, we have
\bq
\ba
\|w\|_{H^{\min\{1+\alpha, \beta_1+\frac{3}{2}, 2-\epsilon\}}(\Omega^-)} \leq & 
\|w\|_{H^{\min\{1+\alpha, \beta_1+\frac{3}{2}, 2-\epsilon\}}(\Omega_0)} 
\leq C \|\widetilde g_1\|_{H^{\min\{\beta_1, \frac{1}{2}-\epsilon\}}(\widetilde \gamma_1^+)}\\
= & C\|g_1\|_{\widetilde{H}^{\min\{\beta_1, \frac{1}{2}-\epsilon\}}(\gamma_1)} 
\leq C \|g_1\|_{H^{\beta_1}(\gamma_1)};
\ea
\eq
if $g_1 \in H_0^{\beta_1}(\gamma_1)$ and $\beta_1+\frac{1}{2}>0$ is not an integer, we have
\bq
\|w\|_{H^{\min\{1+\alpha, \beta_1+\frac{3}{2}\}}(\Omega^-)} \leq \|w\|_{H^{\min\{1+\alpha, \beta_1+\frac{3}{2}\}}(\Omega_0)} \leq C \|\widetilde g_1\|_{H^{\beta_1}(\widetilde \gamma_1^+)}= C\|g_1\|_{H^{\beta_1}(\gamma_1)},
\eq
and if $g_1 \in H_0^{\beta_1}(\gamma_1)$ and $\beta_1+\frac{1}{2}>0$ is an integer, we have
\bq
\|w\|_{H^{\min\{1+\alpha, \beta_1+\frac{3}{2}-\epsilon\}}(\Omega^-)} \leq \|w\|_{H^{\min\{1+\alpha, \beta_1+\frac{3}{2}-\epsilon\}}(\Omega_0)} \leq C \|\widetilde g_1\|_{H^{\beta_1-\epsilon}(\widetilde \gamma_1^+)}\leq C\|g_1\|_{H^{\beta_1}(\gamma_1)},
\eq

Similarly, we can also extend the line fracture $\gamma_1$ to a closed sufficiently smooth curve $\widetilde{\gamma}_1^-$ in $\Omega^-$ as shown in Figure \ref{fig:Omega3}(d), and obtain similar estimates on $\Omega^+$. It can be observed that $w$ is smooth in the neighborhood of $\gamma'_1\setminus \gamma_1$. 
Thus, it follows that if $g_1 \in H^{\beta_1}(\gamma_1)$,
\bq\label{wreg1f}
\|w\|_{H^{\min\{1+\alpha, \beta_1+\frac{3}{2}, 2-\epsilon\}}(\Omega \setminus \gamma_1)} \leq C \|g_1\|_{H^{\beta_1}(\gamma_1)};
\eq
and if $g_1 \in H_0^{\beta_1}(\gamma_1)$, 
\bq\label{wreg1f0}
\ba
\|w\|_{H^{\min\{1+\alpha, \beta_1+\frac{3}{2}\}}(\Omega \setminus \gamma_1)} \leq & C\|g_1\|_{H^{\beta_1}(\gamma_1)}, \quad \text{if $\beta_1+\frac{1}{2}>0$ is not an integer},\\
\|w\|_{H^{\min\{1+\alpha, \beta_1+\frac{3}{2}-\epsilon\}}(\Omega \setminus \gamma_1)} \leq & C\|g_1\|_{H^{\beta_1}(\gamma_1)}, \quad \text{if $\beta_1+\frac{1}{2}>0$ is an integer}.
\ea
\eq
We can apply the regularity estimate (\ref{wreg1f}) or (\ref{wreg1f0}) to multiple line fractures case and obtain the estimate (\ref{transreg}) or (\ref{transreg1}) by using the superposition principle as discussed in Remark \ref{superposition}.
\end{proof}

By Lemma \ref{transweaklem} and (\ref{eq:2d}b), we can extend the solution $w$ of the transmission problem (\ref{eq:2d}) from  $\Omega\setminus \cup_{l=1}^N\gamma_l$ to the whole domain $\Omega$ by taking
\begin{equation}\label{u_alter}
\bar{w}:=
\left\{
\begin{aligned}
w\quad  & \text{in } \Omega\setminus \cup_{l=1}^N\gamma_l, \\
w^+(=w^-)\quad  & \text{on }  \gamma_l, \ l=1, \cdots, N.
\end{aligned}
\right.
\end{equation}
It is obvious that the extended solution
\bq\label{barwcon}
\bar{w} \in \mathcal C^0(\Omega) \cap H_0^1(\Omega),
\eq
and
\bq\label{wrela}
\|\bar{w}\|_{H^1(\Omega)} = \|\bar{w}\|_{V}=\|w\|_{V}.
\eq
Therefore, (\ref{weaktrans}) can be extended to the weak formulation
\bq\label{weaktransnew}
a(\bar{w},v)= \sum_{l=1}^N\int_{\gamma_l} g_l v ds, \quad v\in H_0^1(\Omega).
\eq

\begin{thm}\label{equiv}
Let $u$ be the solution of problem (\ref{eq:Possion}), and $w$ be the solution of the transmission problem (\ref{eq:2d}), then it follows
\bq\label{exeq}
u \equiv \bar{w} \quad \text{in } \Omega,
\eq
where $\bar{w}$ is the extended solution of $w$ in $\Omega$ by (\ref{u_alter}).
\end{thm}
\begin{proof}
We set $\tilde u = u-\bar{w}$ and subtract (\ref{weaktransnew}) from (\ref{eqn.weak}), we have that
$$
a(\tilde u, v) = 0, \quad v \in H_0^1(\Omega).
$$
Set $v=\tilde u \in H_0^1(\Omega)$, we further have
$$
C\|\tilde u\|^2_{H^1(\Omega)}\leq a(\tilde u, \tilde u) = 0,
$$
which gives
$$
\tilde u =0, \quad \text{in } H_0^1(\Omega).
$$
Thus, by Lemma \ref{thm2-2} we have
\bq\label{dwbdd}
\|\nabla w\|_{L^2(\Omega)} \leq C\|\bar{w}\|_{H^1(\Omega)}=C\|u\|_{H^1(\Omega)} \leq C \sum_{l=1}^N \|g_l\|_{L^2(\gamma_l)}.
\eq
Next, we consider closed region $R_\epsilon$ enclosing all line fractures such that $\Omega = R_\epsilon \cup (\Omega \setminus R_\epsilon)$, and denote $\mathbf{n}_\epsilon$ the unit outward norm vector of $\Omega \setminus R_\epsilon$ (inward for $R_\epsilon$) on $\partial R_\epsilon$. 
For $\forall v\in C_0^\infty(\Omega)$,
\bqs
\ba
- \int_\Omega \Delta \bar{w} v dx =  -\lim_{R_\epsilon\rightarrow \cup_{l=1}^N \gamma_l}\left(\int_{\Omega \setminus R_\epsilon} \Delta w v dx+ \int_{R_\epsilon} \Delta w v dx\right) = -\lim_{R_\epsilon\rightarrow \cup_{l=1}^N \gamma_l} \int_{R_\epsilon} \Delta w v dx,
\ea
\eqs
where we have used (\ref{eq:2d}a) in the second equality, namely, $\Delta w = 0$ in $\Omega \setminus R_\epsilon$.

Then for $\forall v\in C_0^\infty(\Omega)$ we have,
\bqs
\ba
-\int_\Omega \Delta \tilde u v dx =& -\int_\Omega \Delta u v dx + \int_\Omega \Delta \bar{w} v dx = -\int_\Omega \Delta u v dx + \lim_{R_\epsilon\rightarrow \cup_{l=1}^N \gamma_l}\int_{R_\epsilon} \Delta w v dx.
\ea
\eqs
Applying (\ref{eq:Possion}) to the first term and Green's formula to the second term on the right hand side of the equation above, we have
\bqs
\ba
-\int_\Omega \Delta \tilde u v dx =& \sum_{l=1}^N\int_{\gamma_l} g_l v ds - \lim_{R_\epsilon\rightarrow \cup_{l=1}^N \gamma_l}\left(\int_{\partial R_\epsilon} \partial_{\mathbf{n}_\epsilon} w v ds + \int_{R_\epsilon} \nabla w \nabla v dx\right).
\ea
\eqs
By (\ref{dwbdd}) and the boundedness of $\nabla v$, we have
\bqs
\ba
\left|\int_{R_\epsilon} \nabla w \nabla v dx\right| \leq \|\nabla w\|_{L^2(R_\epsilon)} \|\nabla v\|_{L^2(R_\epsilon)} \leq \|\nabla w\|_{L^2(\Omega)} \|\nabla v\|_{L^\infty(\Omega)} |R_\epsilon| \rightarrow 0,
\ea
\eqs
as $R_\epsilon\rightarrow \cup_{l=1}^N \gamma_l$.\\
It can be observed
\bqs
\ba
\int_{\partial R_\epsilon} \partial_{\mathbf{n}_\epsilon} w v ds 
\rightarrow \sum_{l=1}^N\int_{\gamma_l} [\partial_{\mathbf{n}}w ] v ds = \sum_{l=1}^N\int_{\gamma_l} g_l v ds,
\ea
\eqs
as $R_\epsilon\rightarrow \cup_{l=1}^N \gamma_l$.
From the discussion above, we have
\bqs
\ba
-\int_\Omega \Delta \tilde u v dx = 0, \quad v \in C_0^\infty(\Omega).
\ea
\eqs
Since $v$ is arbitrary, so it follows that
$$
-\Delta \tilde{u} = 0 \quad \text{in } \Omega,
$$
which together with the boundary condition $\tilde{u}=u-\bar{w}=0$ on $\partial \Omega$ gives $\tilde{u}\equiv 0$ in $\Omega$.
\end{proof}

Theorem \ref{equiv} indicates that the extension $\bar{w}$ of the solution of the transmission problem (\ref{eq:2d}) by (\ref{u_alter}) solves elliptic problem (\ref{eq:Possion}).

\begin{coro}\label{co12}
For $\forall \epsilon>0$, we have the following results,\\
(i) let $u$ be the solution of problem \eqref{eq:Possion}, then $u|_{\Omega\setminus \cup_{l=1}^N\gamma_l}$ solves the transmission problem (\ref{eq:2d});\\
(ii) if $g_l \in H^{\beta_l}(\gamma_l)$, $l=1,\cdots, N$, it follows
$$
u \in H^{\frac{3}{2}-\epsilon}(\Omega) \cap H^{\min\{\alpha+1, \beta+\frac{3}{2}, 2-\epsilon\}}(\Omega\setminus \cup_{l=1}^N\gamma_l).
$$
Further, if
$g_l \in H_0^{\beta_l}(\gamma_l)$, $l=1,\cdots, N$, it follows
\bq
\ba
& u \in H^{\frac{3}{2}-\epsilon}(\Omega) \cap H^{\min\{\alpha+1, \beta+\frac{3}{2}\}}( \Omega\setminus \cup_{l=1}^N\gamma_l), \quad \text{if $\beta+\frac{1}{2}>0$ is not an integer},\\
& u \in H^{\frac{3}{2}-\epsilon}(\Omega) \cap H^{\min\{\alpha+1, \beta+\frac{3}{2}-\epsilon\}}( \Omega\setminus \cup_{l=1}^N\gamma_l), \quad \text{if $\beta+\frac{1}{2}>0$ is an integer},
\ea
\eq
where $\alpha$, $\beta$ are given in Theorem \ref{wreglem}.
\end{coro}
\begin{proof}
The proof follows from Lemma \ref{thm2-2}, Lemma \ref{transweaklem}, Theorem \ref{wreglem} and Theorem \ref{equiv}.
\end{proof}


\begin{rem}
The conclusions in Theorem \ref{wreglem} and Corollary \ref{co12} still hold, if $\gamma_l$, $l=1,\cdots, N$ are sufficiently smooth curved line segments.
\end{rem}

\section{Adaptive finite element method}\label{sec-3}

Let $\mathcal{T} = \{T\}_{T \in \mathcal{T}}$ be a triangulation of $\Omega$ with triangles. 
Denote the set of edges of $\mathcal{T}$ by $\mathcal{E} = \mathcal{E}_I\cup\mathcal{E}_B$, where $\mathcal{E}_I$ and $\mathcal{E}_B$ represent the set of the interior edges and the boundary edges, respectively.
For any triangle $T\in \mathcal{T}$, we denote $h_T$ the diameter of $T$.

The Lagrange finite element space is defined by
\bqs
S(\mathcal{T})=\{v\in \mathcal C^0(\Omega) \cap H_0^1(\Omega):v|_T\in P_k(T), \ \forall \ T \in \mathcal{T}\},
\eqs
where $P_k(T)$ is the space of polynomials with degree less than or equal to $k$ on $T$.

\subsection{Standard finite element method}

We suppose that the mesh $\mathcal T$ consists of quasi-uniform triangles with mesh size $h:=\max h_T$.
Based on the variational formulation (\ref{eqn.weak}) and (\ref{weaktrans}), the standard finite element solution for problem \eqref{eq:Possion} is to find $u_h\in S(\mathcal{T})$ such that
\begin{equation}
\int_\Omega\nabla u_h\cdot\nabla v_hdx= \sum_{l=1}^N \int_{\gamma_l} g_l(s) v_h(s) ds, \quad \forall\ v_h\in S(\mathcal{T}).\label{FEM form}
\end{equation}
Because of the lack of regularity in the solution $u\in H^{\frac{3}{2}-\epsilon}(\Omega)$ for $\forall \epsilon>0$ (see Lemma \ref{thm2-2}), \ the standard error estimate \cite{Ciarlet74} on general quasi-uniform meshes which allow the line fractures pass through the triangles yields only a suboptimal convergence rate,
\begin{equation}\label{H1rate1}
\|u-u_h\|_{H^1(\Omega)}\leq C h^{\frac{1}{2}-\epsilon}.
\end{equation}

If we further assume that the quasi-uniform mesh $\mathcal{T}$ conforms to line fractures $\gamma_{l}$. Namely, $\gamma_{l}$ are the union of some edges in $\mathcal{E}_I$ and
do not cross with any triangles in $\mathcal{T}$.
By Corollary \ref{co12}, the standard error estimate of the finite element approximations on conforming quasi-uniform meshes gives a better convergence rate compared with (\ref{H1rate1}), if all $g_l \in H^{\beta_l}(\gamma_l)$, it follows
\begin{equation}\label{H1rate2}
\|u-u_h\|_{H^1(\Omega)}\leq C h^{\min\{\alpha, \beta+\frac{1}{2}, 1-\epsilon\}},
\end{equation}
and if all $g_l \in H_0^{\beta_l}(\gamma_l)$, it follows
\begin{equation}\label{H1rate3}
\ba
& \|u-u_h\|_{H^1(\Omega)}\leq C h^{\min\{k,\alpha, \beta+\frac{1}{2}\}}, \quad \text{if $\beta+\frac{1}{2}>0$ is not an integer},\\
& \|u-u_h\|_{H^1(\Omega)}\leq C h^{\min\{k,\alpha, \beta+\frac{1}{2}-\epsilon\}}, \quad \text{if $\beta+\frac{1}{2}>0$ is an integer},
\ea
\end{equation}
where $\alpha$, $\beta$ are given in Theorem \ref{wreglem}.

The singularities in the solution can severely slow down the convergence of the standard finite element method associated with the quasi-uniform meshes.
To improve the convergence rate, we introduce an adaptive finite element method to approximate the solution of problem \eqref{eq:Possion}.

\subsection{The adaptive finite element method}
In the following, we first derive a residual-based error estimator and show its reliability and efficiency. Based on the derived error estimator and bisection mesh refinement method, we then propose an adaptive finite element algorithm.

To propose an efficient and reliable residual-based error estimator, one of choices is to regularize the source term such that the regularized source term belongs to $L^2(\Omega)$ or $L^p(\Omega)$ with $1<p<\infty$ \cite{HL21, MMR22}.
Therefore, the residual-based a posteriori error estimator for the usual Poisson equation can be applied.
Let the function $g^r\in L^2(\Omega)$ be a regularized function of the source term $\sum_{l=1}^N g_l\delta_{\gamma_l}$ in (\ref{eq:Possion}), then the classical  residual-based a posteriori error estimator is given by
\begin{align}\label{oldestimator}
\xi = \left(\sum_{T \in \mathcal{T}} \xi_T^2(u_h)\right)^{\frac{1}{2}},
\end{align}
where the local indicator satisfying
\begin{align}\label{poisson estimator}
\xi_T(u_h)^2 = h_T^2\|\Delta u_h+g^r\|_{L^2(T)}^2 + \frac{1}{2}\sum_{e \in \partial T \cap \mathcal{E}_I}h_T\|[\partial_\mathbf{n}u_h]\|_{L^2(e)}^2,
\end{align}
where $[\partial_\mathbf{n}u_h]$ denotes the jump of the normal derivatives of $u_h$ on the interior edges of element $T$.

The regularization technique introduced above is an effective approach to propose adaptive finite element algorithm. However, the corresponding adaptive finite element solution involves not only the discretization error but also the regularization error. When it applied to problem (\ref{eq:Possion}), it may lead to over-refinements on the meshes or low convergence rates for high order approximations. 

For analysis convenience, we extend $g_l$ from $\gamma_l$ to $\mathcal{E}_I$ by defining
\begin{align}
f =
\begin{cases}
g_l,\quad e \in \gamma_l, \quad l=1,\cdots,N,\\
0,\quad e \in \mathcal{E}_I \backslash \cup_{l=1}^N \gamma_l.
\end{cases}
\end{align}
Let $\mathbf{n}$ be the outward unit normal derivative of triangle $T \in \mathcal{T}$.
By Corollary \ref{co12}, we have $[\partial_{\mathbf{n}} u]=g_l=f$ for $e \in \gamma_l$, and $[\partial_\mathbf{n} u]=0=f$ for $e \in \mathcal{E}_I \backslash \cup_{l=1}^N \gamma_l$, so $[\partial_{\mathbf{n}} u]$ is also extended to $\mathcal{E}_I$ in the sense
\begin{align}
[\partial_{\mathbf{n}} u]|_e = {f}|_e,  \quad e \in \mathcal{E}_I.\label{property of fl}
\end{align}

Motivated by the equivalence of the elliptic problem (\ref{eq:Possion}) and the transmission problem (\ref{eq:2d}) in the domain excluding the line fractures,
we propose the following residual-based a posteriori error estimator,
\begin{align}
\eta = \left(\sum_{T \in \mathcal{T}} \eta_T^2(u_h)\right)^{\frac{1}{2}},\label{total error estimator}
\end{align}
where the local indicator on $T \in \mathcal{T}$ is defined by, 
\begin{align}\label{local error estimator}
\eta_T(u_h)^2 = h_T^2\|\Delta u_h\|_{L^2(T)}^2 + \frac{1}{2}\sum_{e \in \partial T \cap \mathcal{E}_I}h_T\|f-[\partial_\mathbf{n}u_h]\|_{L^2(e)}^2.
\end{align}

\begin{rem}
If the equation in (\ref{eq:Possion}) has an additional source term $q(x)\in L^2(\Omega)$, namely,
\begin{equation*}
-\Delta u   = \sum_{l=1}^N g_l\delta_{\gamma_l} + q(x)      \quad \text{in }  \Omega,   
\end{equation*}
then the transmission problem (\ref{eq:2d}a) will be modified as
\begin{subequations}
\begin{align*}
-\Delta w=q(x) &\quad \mbox{ in   }\Omega \setminus \cup_{l=1}^N\gamma_l,
\end{align*}
\end{subequations}
and the local indicator (\ref{local error estimator}) is therefore given by
\begin{align*}
\eta_T(u_h)^2 = h_T^2\|\Delta u_h + q(x)\|_{L^2(T)}^2 + \frac{1}{2}\sum_{e \in \partial T \cap \mathcal{E}_I}h_T\|f-[\partial_\mathbf{n}u_h]\|_{L^2(e)}^2.
\end{align*}
\end{rem}

Before we present the efficiency and reliability of the proposed a posteriori error estimator (\ref{total error estimator}), we first prepare some necessary inequalities and estimates.

\begin{lem}[Trace inequality \cite{BS2008}]\label{Trace inequality}
For any element $T \in \mathcal{T}$, $\forall e  \subset \partial T$, we have
\begin{align*}
\|v\|_{L^2(e)} &\le C h_T^{-1/2}(\|v\|_{L^2(T)} + h_T\|\nabla v\|_{L^2(T)}),\qquad \forall v \in H^{1}(T).
\end{align*}
\end{lem}

\begin{lem}[Inverse inequality \cite{BS2008}]\label{Inverse inequality}
For any element $T\in\mathcal{T}$ and $v \in P_{k}(T)$, $\forall \,e \subset \partial T$, we have
\begin{align*}
\|\nabla^j v\|_{L^2(T)}&\le C h_T^{-j}\|v\|_{L^2(T)},\qquad \forall \, 0\le j\le k.
\end{align*}
\end{lem}
\begin{lem}[Interpolant error estimate \cite{V96}]\label{Interpolation error estimate}
For any $v \in H^{l}(T),\,l\ge1$, it follows
\begin{align*}
\|v-\pi v\|_{H^m(T)} &\le C h^{l-m}\|v\|_{H^l(T)},
\end{align*}where $m = 0,1$ and $\pi v \in S(\mathcal{T})$ represents the nodal interpolant of $v$.
\end{lem}



In the following analysis, we make use of the equivalence of problem (\ref{eq:Possion}) to the transmission problem (\ref{eq:2d}) as discussed in Section \ref{31}, and we pay special attention to handle the flux jumps (\ref{property of fl}) on line fractures $\gamma_l$ in the following reliability analysis.

\begin{thm}[Reliability]
Assume that $u$ and $u_h$ are the solution of (\ref{eq:Possion}) and (\ref{FEM form}), respectively. Then the residual-based a posteriori error estimator $\eta$ satisfies the global bound,
\begin{equation}
\|\nabla(u - u_h)\|_{L^2(\Omega)} \le C\eta(u_h).\label{Reliability inequation}
\end{equation}
\end{thm}
\begin{proof}
Let $e_u = u - u_h$, we have
\begin{align}
\|\nabla e_u\|_{L^2(\Omega)}^2 = \int_{\Omega}\nabla e_u \cdot \nabla e_u\,dx = \int_{\Omega}\nabla e_u \cdot \nabla (e_u - \pi e_u)\,dx,\label{proof 1}
\end{align}
where we have used the Galerkin orthogonality to subtract an interpolant $\pi e_u \in S(\mathcal{T})$ to $e_u$.
Note that by Corollary \ref{co12}, we have
\begin{align}
\Delta e_u=\Delta u - \Delta u_h = - \Delta u_h, \quad \text{in } \Omega\setminus\cup_{l=1}^N\gamma_l. \label{Corollary 2}
\end{align}
Thus splitting \eqref{proof 1} into a sum over the elements and using Green's formula, we have
\begin{align*}
\sum_{T \in \mathcal{T}}\int_{T}\nabla e_u \cdot \nabla (e_u - \pi e_u)\,dx &=\sum_{T \in \mathcal{T}}\int_{T}-\Delta e_u(e_u - \pi e_u)\,dx + \int_{\partial T}\mathbf{n} \cdot \nabla e_u(e_u - \pi e_u)\,ds\\
& = \sum_{T \in \mathcal{T}} \left(\int_T \Delta u_h (e_u - \pi e_u)\,dx+\int_{\partial T\cap \mathcal{E}_I}\mathbf{n} \cdot \nabla e_u(e_u - \pi e_u)\,ds\right),
\end{align*}
where we have used $\pi e_u = e_u = 0$ on $\partial \Omega$.
Note that $e_u$ is continuous by Corollary \ref{co1} and the continuity of the finite element solution, so we have $(e_u^+ - \pi e_u^+)|_e = (e_u^- - \pi e_u^-)|_e$ for any $e = \partial T_+ \cap \partial T_- \in \mathcal{E}_I$. Thus, it follows
\begin{align*}
&\int_{e\cap \partial T^+}\mathbf{n} \cdot \nabla e_u(e_u - \pi e_u)\,ds+\int_{e\cap \partial T^-}\mathbf{n} \cdot \nabla e_u(e_u - \pi e_u)\,ds\\ = &\int_e\mathbf{n}^+\cdot \nabla e_u^+(e_u^+ - \pi e_u^+) +\mathbf{n}^-\cdot \nabla e_u^-(e_u^- - \pi e_u^-)\,ds\\
= &\int_e \left((\mathbf{n}^+\cdot \nabla u^++ \mathbf{n}^-\cdot \nabla u^-)-(\mathbf{n}^+\cdot \nabla u^+_h+ \mathbf{n}^-\cdot \nabla u^-_h)\right)(e_u - \pi e_u)\,ds\\
=& \int_e[\partial_{\mathbf{n}}u](e_u - \pi e_u)\,ds-\int_e[\partial_{\mathbf{n}}u_h](e_u - \pi e_u)\,ds.
\end{align*}
This, together with \eqref{property of fl}, implies that
\begin{align*}
\sum_{T \in \mathcal{T}}\int_{\partial T\cap \mathcal{E}_I}\mathbf{n} \cdot \nabla e_u(e_u - \pi e_u)\,ds &= \sum_{e \in \mathcal{E}_I}\left(\int_e[\partial_{\mathbf{n}} u](e_u - \pi e_u)\,ds-\int_e[\partial_{\mathbf{n}} u_h](e_u - \pi e_u)\,ds\right)\\
&= \sum_{e \in \mathcal{E}_I}\int_e(f-[\partial_{\mathbf{n}} u_h])(e_u - \pi e_u)\,ds.
\end{align*}
Returning to the sum over the elements with simply distributing half of $f-[\partial_{\mathbf{n}} u_h]$ on $T_+$ and the remaining half on $T_-$, we have
\bq\label{proof 4}
\begin{aligned}
\|\nabla e_u\|_{L^2(\Omega)}^2 =  \sum_{T \in \mathcal{T}} \left(\int_T \Delta u_h (e_u - \pi e_u)\,dx  + \frac{1}{2}\sum_{e \in\partial T\cap \mathcal{E}_I}\int_{e}(f-[\partial_{\mathbf{n}} u_h])(e_u - \pi e_u)\,ds\right).
\end{aligned}
\eq
Next, we estimate the terms on the right hand side of (\ref{proof 4}) one by one.

Using Cauchy-Schwarz inequality and Lemma \ref{Interpolation error estimate}, we have
\begin{align}
\int_T \Delta u_h (e_u - \pi e_u)\,dx &\le \|\Delta u_h\|_{L^2(T)}\|e_u - \pi e_u\|_{L^2(T)} \le Ch_T\|\Delta u_h\|_{L^2(T)}\|\nabla e_u\|_{L^2(T)}.\label{proof 6}
\end{align}
Then, using Cauchy-Schwarz inequality, Lemma \ref{Trace inequality}, and Lemma \ref{Interpolation error estimate}, we have
\begin{align}
\int_{e}&(f -[\partial_{\mathbf{n}} u_h])(e_u -\pi e_u)\,ds\le \|f -[\partial_{\mathbf{n}} u_h]\|_{L^2(e)}\|e_u - \pi e_u\|_{L^2(e)}\nonumber\\
& \le C\left(h_T^{-1}\|e_u - \pi e_u\|_{L^2(T)}^2 + h_T\|\nabla (e_u -\pi e_u)\|_{L^2(T)}^2\right)^{1/2}\|f -[\partial_{\mathbf{n}} u_h]\|_{L^2(e)}\nonumber\\
& \le Ch_T^{1/2}\|f -[\partial_{\mathbf{n}} u_h]\|_{L^2(e)}\|\nabla e_u\|_{L^2(T)}. \label{proof 5}
\end{align}
The estimate (\ref{Reliability inequation}) now follows from \eqref{proof 4}-\eqref{proof 5}.
\end{proof}

Let $\overline{f} \in P_{k}(e)$ be the $L^2$-projection of $f$. We define the oscillation on $e \in \mathcal{E}_I$ by
\begin{align*}
osc(e)^2 = h_e\|f - \overline{f}\|_{L^2(e)}^2,
\end{align*}
where $h_e$ is the length of $e$. 
Let $e = \partial T_+ \cap \partial T_-$ with $T_+$ and $T_-$ being two adjacent triangles, and we set $\omega_e = T_+ \cup T_-$, then for any $T \in \omega_e$ there exist positive constants $C_1$ and $C_2$ such that
$$
C_1h_T \le h_e \le C_2h_T.
$$

For a triangle $T \in \mathcal{T}$ with vertices $x_1,x_2,x_3$, 
we denote $(\lambda_{x_1}, \lambda_{x_2}, \lambda_{x_3})$ the barycentric coordinates on $T$.
We define a bubble function $b_T$ in $T$ by
\begin{align}\label{bubbbT}
b_T = 27\lambda_{x_1}\lambda_{x_2}\lambda_{x_3}.
\end{align}
For an edge $e=x_ix_j \in \partial T \subset \mathcal{E}$, we define an edge bubble function $b_e$ in $T$ by
\begin{align}\label{bubbbe}
b_e = 4\lambda_{x_i}\lambda_{x_j}.
\end{align}

For the bubble functions (\ref{bubbbT}) and (\ref{bubbbe}), we have the following results. 

\begin{lem}[\cite{V1994}]\label{bT properties}
For the element bubble function $b_T$ in (\ref{bubbbT}), it has the following properties,
\begin{align}
0 &\le b_T(x) \le 1,\quad \forall x \in T,\qquad
{\rm and}\qquad b_T(x) = 0,\quad \forall x \in \partial T,\label{bT property 1}
\end{align}
Moreover, for any $v \in P_k$, it follows
\begin{align}
\|v\|_{L^2(T)}&\le C\|b_T^{1/2}v\|_{L^2(T)}.\label{bT property 3}
\end{align}
\end{lem} 

\begin{lem}[\cite{V1994}]\label{be properties}
For $e = \partial T_+ \cap \partial T_-$, the edge bubble function $b_e$ defined by (\ref{bubbbe}) has the following properties,
\begin{align}
0 &\le b_e(x) \le 1,\quad \forall x \in \omega_e,\qquad {\rm and}\qquad 
b_e(x) = 0,\quad \forall x \in \partial \omega_e \setminus e, \label{be property 1}
\end{align}
where $\partial \omega_e = \partial T_+ \cup \partial T_-$.
Moreover, for any $v \in P_k$, it follows
\begin{align}
\|v\|_{L^2(e)}&\le C\|b_e^{1/2}v\|_{L^2(e)},\label{be property 3}\\
\|\nabla (b_ev)\|_{L^2(\omega_e)} &\le C h_e^{-1/2}\|v\|_{L^2(e)},\label{be property 4}\\
\|b_ev\|_{L^2(\omega_e)} &\le C h_e^{1/2}\|v\|_{L^2(e)}.\label{be property 5}
\end{align}
\end{lem} 
\begin{thm}[Efficiency]\label{Efficiency}
For the local indicator $\eta_T$ defined in (\ref{local error estimator}), it follows
\begin{align}
\eta_T(u_h) \le C\left(\|\nabla e_u\|_{L^2(\omega_T)} +  osc(\partial T)\right), \quad \forall T \in \mathcal{T},
\end{align}
where $w_T = \cup_{e\in \partial T} w_e$, and
\begin{align*}
osc(\partial T)^2 = \sum_{e \in \partial T }osc(e)^2.
\end{align*}
\end{thm}
\begin{proof}
Using Green's formula, \eqref{Corollary 2} and \eqref{bT property 1}, we have
\begin{align}
\int_T \nabla e_u \, \nabla(\Delta u_h b_T)\,dx & = -\int_T \Delta e_u \,\Delta u_h b_T\,dx + \int_{\partial T}\nabla e_u \cdot \mathbf{n}\,\Delta u_h b_T\,ds  = \int_T \Delta u_h \,\Delta u_h b_T\,dx,\label{Efficiency proof 0}
\end{align}
Since $\Delta u_h$ is a piecewise polynomial over $\mathcal{T}$, according to \eqref{bT property 3} we have
\begin{align*}
\|\Delta u_h\|_{L^2(T)}^2&\le C \|\Delta u_h b_T^{1/2}\|_{L^2(T)}^2.
\end{align*}
Using the Cauchy-Schwarz inequality, Lemma \ref{Inverse inequality}, and \eqref{bT property 1}, it follows that
\begin{align*}
\|\Delta u_h\|_{L^2(T)}^2&\le C\int_T \nabla e_u \, \nabla(\Delta u_h b_T)\,dx  \le C\|\nabla e_u\|_{L^2(T)}\|\nabla(\Delta u_h b_T)\|_{L^2(T)}\\
& \le Ch_T^{-1}\|\nabla e_u\|_{L^2(T)}\|\Delta u_h b_T\|_{L^2(T)} \le Ch_T^{-1}\|\nabla e_u\|_{L^2(T)}\|\Delta u_h\|_{L^2(T)},
\end{align*}
which gives
\begin{align}
h_T\|\Delta u_h\|_{L^2(T)}\le C\|\nabla e_u\|_{L^2(T)}.\label{Efficiency proof 3}
\end{align}

We now extend $\overline{f} -[\partial_{\mathbf{n}} u_h]$ from edge $e$ to $w_e$ by taking constants along the normal on $e$. The resulting extension $E(\overline{f}-[\partial_{\mathbf{n}} u_h])$ is a piecewise polynomial in $\omega_e$, then according to \eqref{be property 4}-\eqref{be property 5}, we have
\begin{align}
\|\nabla (E(\overline{f} - [\partial_{\mathbf{n}} u_h])b_e)\|_{L^2(\omega_e)}\le C h_e^{-\frac{1}{2}}\|\overline{f} - [\partial_{\mathbf{n}} u_h]\|_{L^2(e)}, \label{extest 1} \\
\|E(\overline{f} - [\partial_{\mathbf{n}} u_h])b_e\|_{L^2(\omega_e)}\le C h_e^{\frac{1}{2}}\|\overline{f} - [\partial_{\mathbf{n}} u_h]\|_{L^2(e)}.\label{extest 2}
\end{align}
Using arguments similar to those leading to \eqref{Efficiency proof 0}, it follows
\begin{align*}
\int_{\omega_e}& \nabla e_u \, \nabla(E(\overline{f} -[\partial_{\mathbf{n}} u_h]) b_e)\,dx = \sum_{T \in \omega_e}\int_{T} \nabla e_u \, \nabla(E(\overline{f} -[\partial_{\mathbf{n}} u_h]) b_e)\,dx \\
& = \sum_{T \in \omega_e} \left(\int_{T} - \Delta e_u \,E(\overline{f} -[\partial_{\mathbf{n}} u_h]) b_e\,dx + \int_{\partial T}\nabla e_u \cdot \mathbf{n}\,E(\overline{f} -[\partial_{\mathbf{n}} u_h]) b_e\,ds \right)\\
& = \sum_{T \in \omega_e}  \left(\int_{T} \Delta u_h \,E(\overline{f} -[\partial_{\mathbf{n}} u_h]) b_e\,dx+ \int_{\partial T}\nabla e_u \cdot \mathbf{n}\,E(\overline{f} -[\partial_{\mathbf{n}} u_h])b_e\,ds\right).
\end{align*}
Note that $\overline{f} -[\partial_{\mathbf{n}} u_h]$ and $b_e$ are continuous on $e \in \mathcal{E}_I$, and $b_e =0$ on $\left(\cup_{T \in \omega_e} \partial T \setminus e\right)$, so we have
\begin{align*}
\sum_{T \in \omega_e}&\int_{\partial T}\nabla e_u \cdot \mathbf{n}\,E(\overline{f} -[\partial_{\mathbf{n}} u_h])b_e\,ds\\
&= \int_e \left((\mathbf{n}^+\cdot \nabla u^++ \mathbf{n}^-\cdot \nabla u^-)-(\mathbf{n}^+\cdot \nabla u^+_h+ \mathbf{n}^-\cdot \nabla u^-_h)\right)(\overline{f} -[\partial_{\mathbf{n}}u_h]) b_e\,ds\\
& = \int_e[\partial_{\mathbf{n}} u](\overline{f} -[\partial_{\mathbf{n}} u_h]) b_e\,ds-\int_e[\partial_{\mathbf{n}}  u_h](\overline{f} -[\partial_{\mathbf{n}}  u_h]) b_e\,ds\\
&=\int_e(f-[\partial_{\mathbf{n}}  u_h])(\overline{f} -[\partial_{\mathbf{n}}  u_h]) b_e\,ds,
\end{align*}
where we used \eqref{property of fl} in the last equality.
Therefore, we get
\begin{align*}
\int_{\omega_e} \nabla e_u \, \nabla(E(\overline{f} -[\partial_{\mathbf{n}} u_h]) b_e)\,dx 
= &\int_{\omega_e} \Delta u_h \,E(\overline{f} -[\partial_{\mathbf{n}} u_h]) b_e\,dx\\ &+\int_e(f-[\partial_{\mathbf{n}}  u_h])(\overline{f} -[\partial_{\mathbf{n}}  u_h]) b_e\,ds\\
=& \int_{\omega_e} \Delta u_h \,E(\overline{f} -[\partial_{\mathbf{n}} u_h]) b_e\,dx + \int_e(\overline{f} -[\partial_{\mathbf{n}} u_h])^2 b_e\,ds \\
&+ \int_e(f-\overline{f})(\overline{f} -[\partial_{\mathbf{n}} u_h]) b_e\,ds
\end{align*}
It follows from \eqref{be property 3}, we obtain
\begin{align*}
\|\overline{f}-[\partial_{\mathbf{n}} u_h]\|_{L^2(e)}^2 \le C\|(\overline{f}-[\partial_{\mathbf{n}} u_h]) b_e^{1/2}\|_{L^2(e)}^2.
\end{align*}
Using Cauchy-Schwarz inequality  and \eqref{extest 1}-\eqref{extest 2}, \eqref{be property 1}, we have
\begin{align*}
\|\overline{f}-[\partial_{\mathbf{n}} u_h]\|_{L^2(e)}^2 \le &C\left(\int_{\omega_e} \nabla e_u \, \nabla(E(\overline{f} -[\partial_{\mathbf{n}} u_h]) b_e)\,dx  -\int_{\omega_e} \Delta u_h \,E(\overline{f} -[\partial_{\mathbf{n}} u_h]) b_e\,dx \right.\\
& \left.\quad- \int_e(f-\overline{f})(\overline{f} -[\partial_{\mathbf{n}} u_h]) b_e\,ds \right)\\
\le & C\left(\|\nabla e_u\|_{L^2(\omega_e)}\|\nabla(E([\overline{f} -\partial_{\mathbf{n}} u_h]) b_e)\|_{L^2(\omega_e)}\right.\\
&\quad + \|\Delta u_h\|_{L^2(\omega_e)}\|E([\overline{f} -\partial_{\mathbf{n}} u_h]) b_e\|_{L^2(\omega_e)} \\
& \left. \quad + \|(f-\overline{f})b_e\|_{L^2(e)}\|\overline{f} -[\partial_{\mathbf{n}} u_h]\|_{L^2(e)}\right)  \\
\le &C \left( h_e^{-1/2} \|\nabla e_u\|_{L^2(\omega_e)}\|\overline{f} -[\partial_{\mathbf{n}} u_h]\|_{L^2(e)}\right.\\
&\quad + h_e^{1/2}\|\Delta u_h\|_{L^2(\omega_e)}\|\overline{f} -[\partial_{\mathbf{n}} u_h]\|_{L^2(e)} \\
&\left.\quad+ \|f-\overline{f}\|_{L^2(e)}\|\overline{f} -[\partial_{\mathbf{n}} u_h]\|_{L^2(e)} \right)\\
\le &C h_e^{-1/2} \|\overline{f} -[\partial_{\mathbf{n}} u_h]\|_{L^2(e)}\left(\|\nabla e_u\|_{L^2(\omega_e)}+ h_e\|\Delta u_h\|_{L^2(\omega_e)} + osc(e)\right),
\end{align*}
which gives
\begin{align}
h_e^{\frac{1}{2}} \|\overline{f}-[\partial_{\mathbf{n}} u_h]\|_{L^2(e)} \le  C\left(\|\nabla e_u\|_{L^2(\omega_e)}+ h_e\|\Delta u_h\|_{L^2(\omega_e)} + osc(e)\right).\label{Efficiency proof 5}
\end{align}
Together with the triangle inequality, \eqref{Efficiency proof 3} and \eqref{Efficiency proof 5}, we obtain the estimation
\begin{align}
h_e^{\frac{1}{2}} \|f-[\partial_{\mathbf{n}} u_h]\|_{L^2(e)} \le  C\left(\|\nabla e_u\|_{L^2(\omega_e)} + osc(e)\right).\label{Efficiency proof 6}
\end{align}
The required estimation now follows form \eqref{Efficiency proof 3} and \eqref{Efficiency proof 6}.
\end{proof}

The corresponding algorithm is summarized as follows.
\begin{breakablealgorithm}
\caption{Adaptive finite element algorithm for the elliptic equation.} 
\label{alg1} 
\begin{algorithmic}[1] 
\STATE Input: an initial mesh $\mathcal{T}^0$; a constant $0<\theta <1$; the maximum number of mesh refinements $n$.
\STATE Output: the numerical solution $u_h^n$; a new refined mesh $\mathcal{T}^n$.
\STATE for $i=0$ to $n$ do\\
       \qquad Solve the discrete equation (\ref{FEM form}) for the finite element solution $u_h^i$ on $\mathcal{T}^i$;\\
       \qquad Computing the local error estimation $\eta_{T}^i(u_h^i)$ and the total error estimation $\eta^i(u_h^i)$ by \eqref{local error estimator} and \eqref{total error estimator};\\
       \qquad if $i< n$ then\\
       \qquad \qquad Mark a subset $\widetilde{\mathcal{T}}^i \subset \mathcal{T}^i$ of elements to refined such that,
       $$(\sum_{T\in \widetilde{\mathcal{T}}^i}{\eta_T^i(u_h^i)}^2)^{1/2}\ge \theta \eta^i(u_h^i);$$
       \qquad \qquad Refine each element $T \in \widetilde{\mathcal{T}}^i$ by longest edge bisection to obtain a new mesh $\mathcal{T}^{i+1}$.\\
       \qquad end if\\
       end for
\end{algorithmic}
\end{breakablealgorithm}

\section{Numerical examples}\label{sec-4}

\subsection{The standard finite element method}
In this subsection, we present numerical examples to verify the convergence rate of the standard finite element method solving equation (\ref{eq:Possion}). The quasi-uniform meshes are considered in this subsection, that is, each triangle is divided into four equal triangles in each mesh refinement. Since the solution $u$ is unknown, we use the following numerical convergence rate
\begin{eqnarray}\label{rate}
\mathcal R=\log_2\frac{|u_h^j-u_h^{j-1}|_{H^1(\Omega)}}{|u_h^{j+1}-u_h^j|_{H^1(\Omega)}},
\end{eqnarray}
where $u_h^j$ is the finite element solution on the  mesh $\mathcal T^j$ obtained after $j$  refinements of the initial triangulation $\mathcal T^0$.

\begin{example}\label{P1h} In this example, we test the convergence rates of the finite element solutions on quasi-uniform meshes.
We consider problem (\ref{eq:Possion}) in a square domain $\Omega=(0,1)^2$ with one line fracture $\gamma_1=Q_1Q_2$ for $Q_1=(0.25,0.5)$ and $Q_2=(0.75,0.5)$.
We take the function $g_1=((x-0.25)(0.75-x))^{r_0}+r_1$ on $\gamma_1$. 
For different parameters $r_0$, $r_1$ in Case 1-6 listed in Table \ref{cases}, we show the smoothness of the corresponding function $g_1$, and the regularity for the solution $u$ of problem \eqref{eq:Possion} followed by Corollary \ref{co12}.

\begin{table}
\centering
\caption{Case 1-6 in Example \ref{P1h}.}\vspace{-3mm}
\begin{tabular}{|c||c|c||c||c|}
\hline
Case number &$r_0$  &  $r_1$  &  $g_1 \in$ & $u \in$ \\
\hline
Case 1 &$-\frac{1}{4}+10^{-3}$  &  $1$  &  $H^{\frac{1}{4}}(\gamma_1)$ & $ H^{\frac{3}{2}-\epsilon}(\Omega)\cap H^{\frac{7}{4}}(\Omega \backslash \gamma_1)$ \\
\hline
Case 2 & $\frac{1}{4}+10^{-3}$  &  1  & $H^{\frac{3}{4}}(\gamma_1)$ & $H^{\frac{3}{2}-\epsilon}(\Omega)\cap H^{2-\epsilon}(\Omega \backslash \gamma_1)$ \\
\hline
Case 3 & $0$  & $1$  & $C^\infty(\gamma_1)$ & $H^{\frac{3}{2}-\epsilon}(\Omega)\cap H^{2-\epsilon}(\Omega \backslash \gamma_1)$ \\
\hline
Case 4 & $\frac{1}{4}+10^{-3}$  & $0$  & $H_0^{\frac{3}{4}}(\gamma_1)$& $ H^{\frac{3}{2}-\epsilon}(\Omega)\cap H^{\frac{9}{4}}(\Omega \backslash \gamma_1)$\\
\hline
Case 5 & $\frac{1}{2}+10^{-3}$  &  $0$  &  $H_0^1(\gamma_1)$ & $H^{\frac{3}{2}-\epsilon}(\Omega)\cap H^{\frac{5}{2}}(\Omega \backslash \gamma_1)$ \\
\hline
Case 6 &$1+10^{-3}$  &  $0$  & $H_0^{\frac{3}{2}}(\gamma_1)$ &  $H^{\frac{3}{2}-\epsilon}(\Omega)\cap H^{3-\epsilon}(\Omega \backslash \gamma_1)$ \\
\hline
\end{tabular}\label{cases}
\end{table}

\begin{figure}
\centering
\subfigure[The initial Union-Jack mesh]{\includegraphics[width=0.3\textwidth]{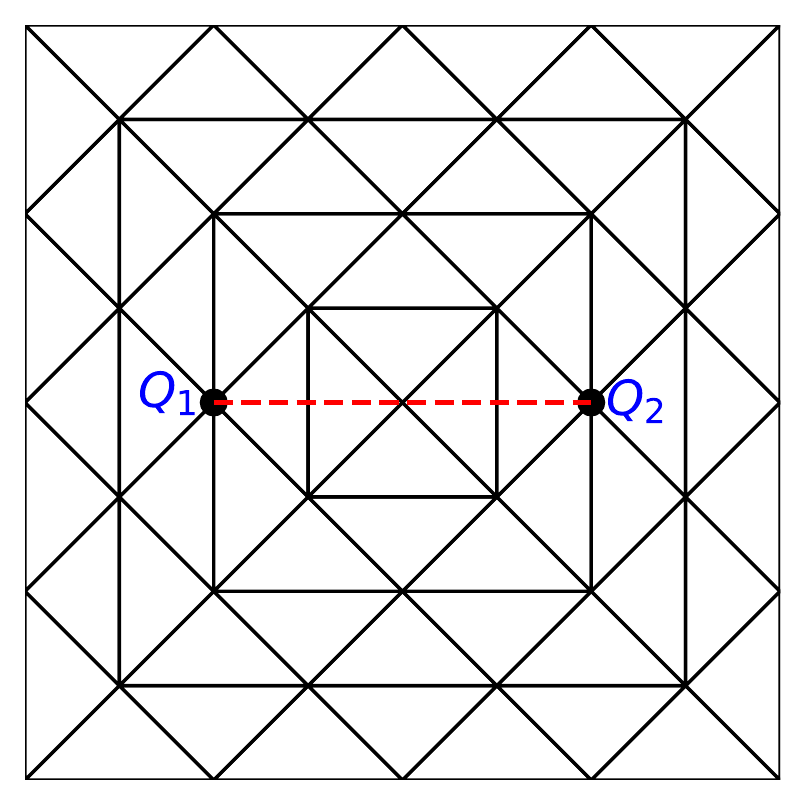}}\hspace{10mm}
\subfigure[The mesh conforming to $\gamma_1$]{\includegraphics[width=0.3\textwidth]{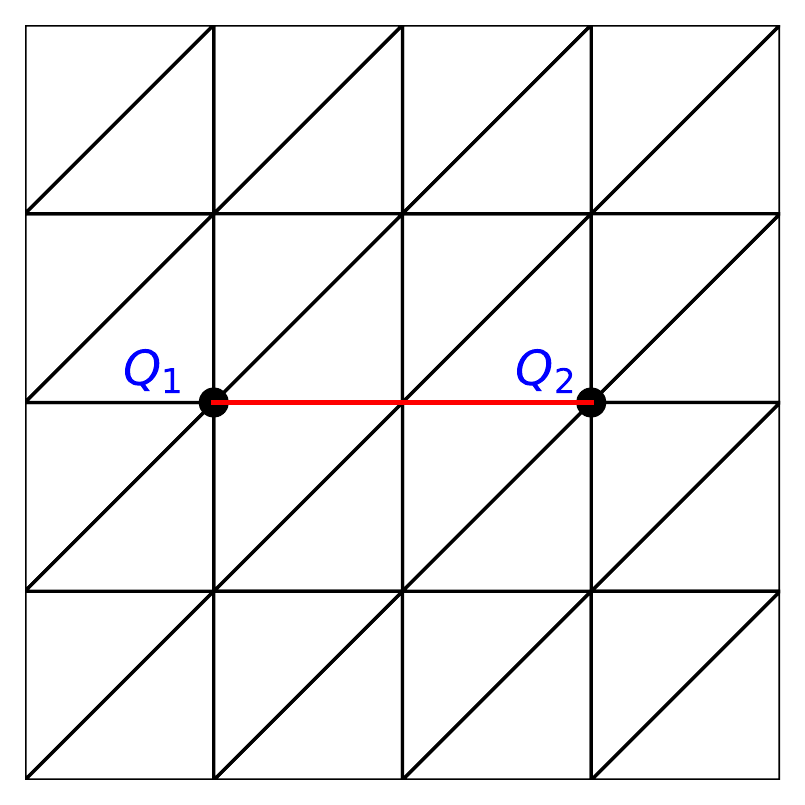}}
\caption{Example \ref{P1h}: the initial meshes.}\label{Mesh_Init}
\end{figure}

\begin{figure}
\centering
\subfigure[Case 1]{\includegraphics[width=0.32\textwidth]{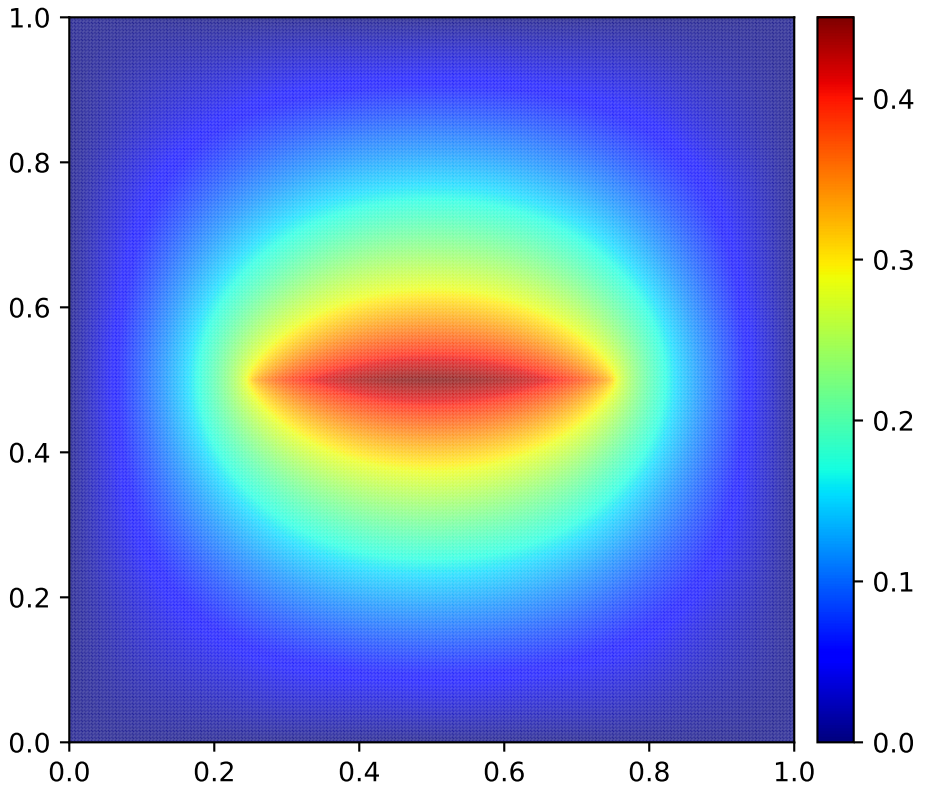}}
\subfigure[Case 2]{\includegraphics[width=0.32\textwidth]{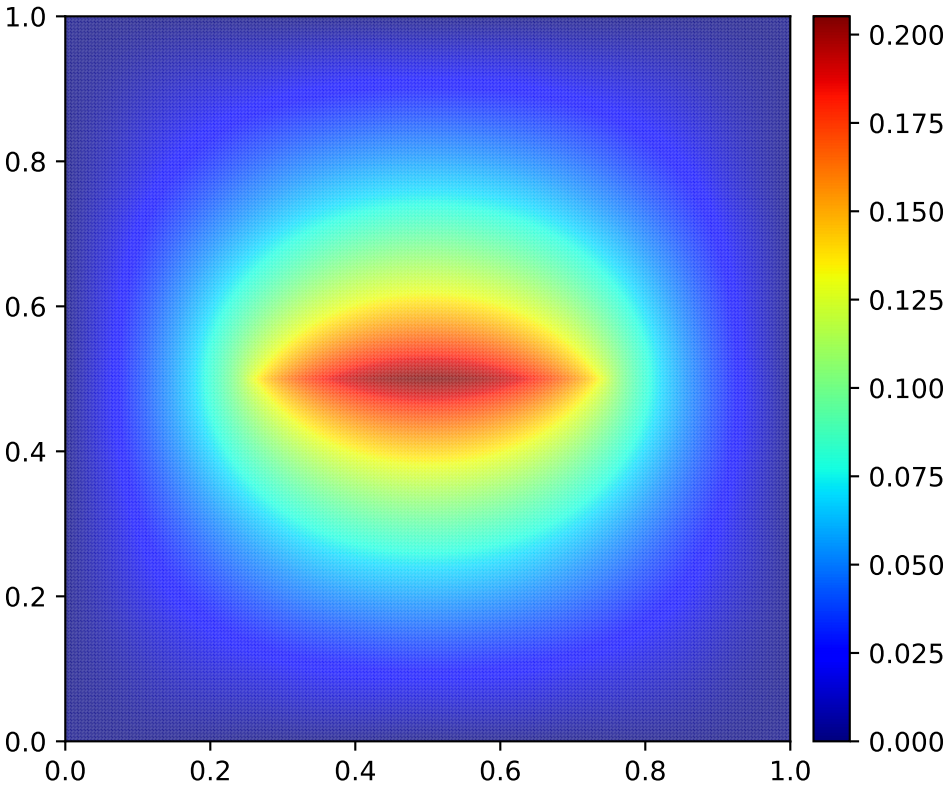}}
\subfigure[Case 3]{\includegraphics[width=0.315\textwidth]{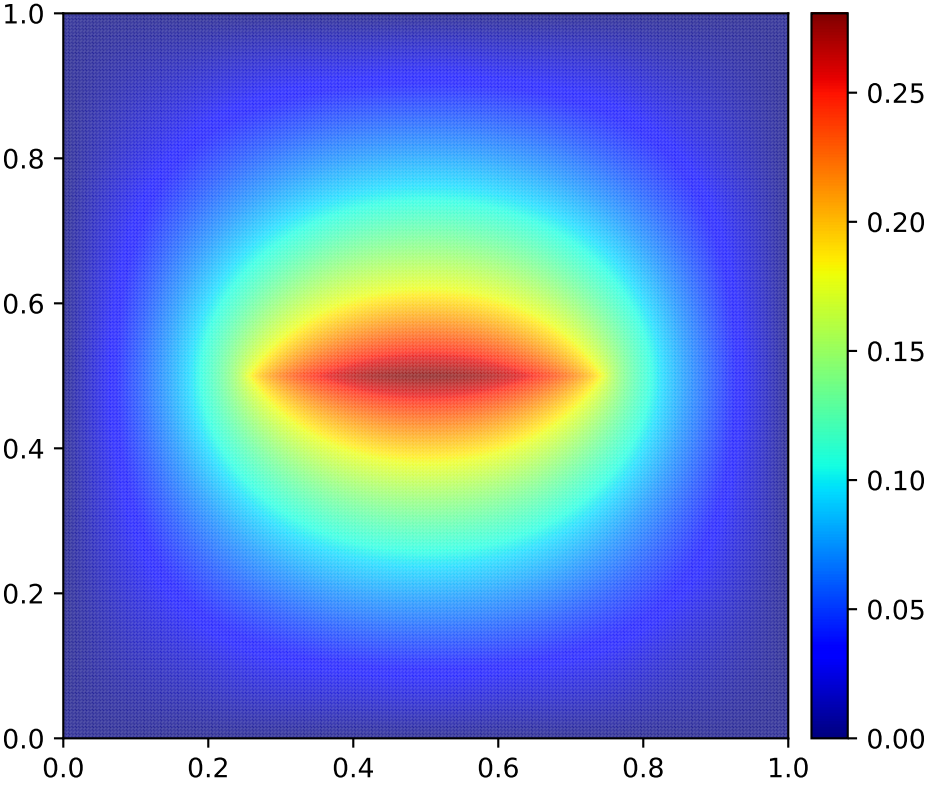}}
\subfigure[Case 4]{\includegraphics[width=0.32\textwidth]{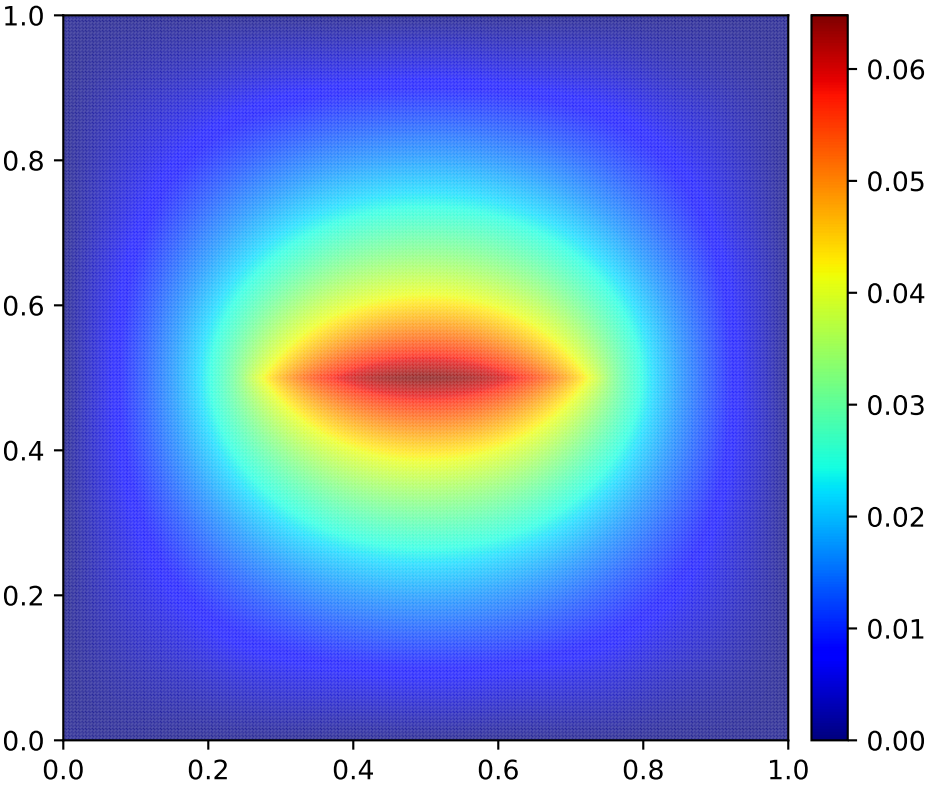}}
\subfigure[Case 5]{\includegraphics[width=0.32\textwidth]{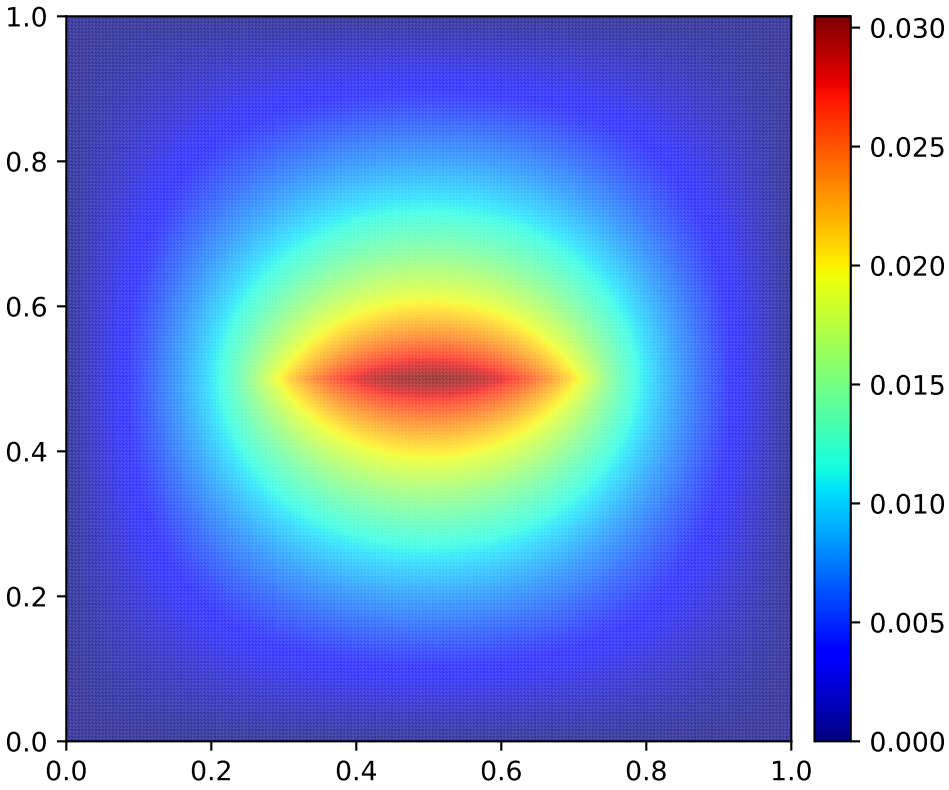}}
\subfigure[Case 6]{\includegraphics[width=0.32\textwidth]{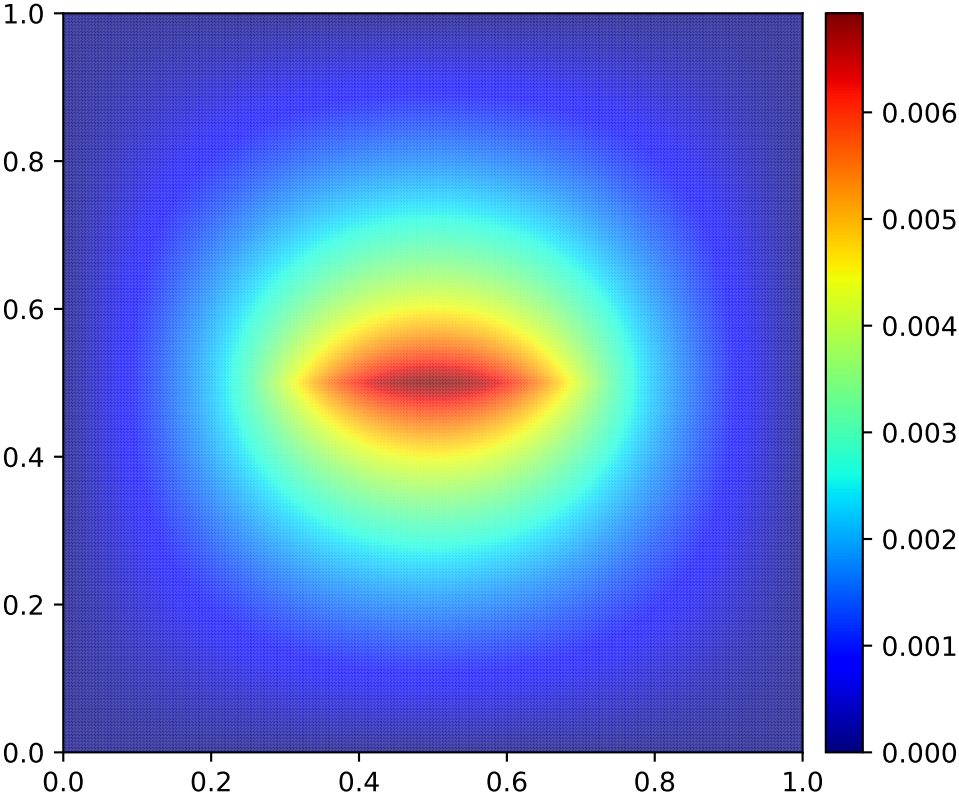}}
\caption{Example \ref{P1h} Test 1: finite element solutions based on $P_1$ polynomials after 8 mesh refinements.}\label{Ex1Con1}
\end{figure}

\noindent\textbf{Test 1.} We take the initial mesh as the Union-Jack mesh and the line fracture $\gamma_1$ pass through the triangles in the mesh as shown in Figure \ref{Mesh_Init}(a). The convergence rates of the finite element solutions based on $P_1,\,P_2$ polynomials are shown in Table \ref{TabConRate}, and we find that suboptimal convergence rates $\mathcal R \approx 0.5$ are obtained for Case 1$-$6, which is due to $u\in H^{\frac{3}{2}-\epsilon}(\Omega),\,\forall \epsilon>0$ regardless of the smoothness of $g_1$ as indicated by Lemma \ref{thm2-2}. 
The contours of the finite element solution for Case 1$-$6 are shown in Figure \ref{Ex1Con1}.

\begin{table}
\centering
\caption{$H^1$ convergence history of finite element solutions in Example \ref{P1h} Test 1 on Union-Jack meshes.}
\begin{tabular}{|c||c|c|c|c||c|c|c|c|c|c|c|c|}
\hline
&  \multicolumn{3}{c}{$\qquad \qquad P_1$} & &\multicolumn{3}{c}{$\qquad  \qquad P_2$} & \\
\hline
j &6 & 7 & 8 &9 & 4&5 & 6 & 7\\
\hline
Case1 &  0.477 &  0.485 &  0.490  &  0.493 &  0.484 &  0.489 &  0.493 &  0.495\\
\hline
Case2 &  0.475 &  0.486 &  0.492 &  0.496 &0.493 &  0.497 &  0.498 &  0.499\\
\hline
Case3 &  0.485 &  0.491 &  0.495 &  0.497&  0.495 &  0.498 &  0.499 &  0.499\\
\hline
Case4 &  0.476 &  0.487 &  0.493 &  0.496 &  0.499 &  0.500   &  0.500   &  0.500\\
\hline
Case5 &  0.476 &  0.487 &  0.493 &  0.497 &  0.503 &  0.501 &  0.500   &  0.500\\
\hline
Case6 &  0.474 &  0.487 &  0.493 &  0.497 &  0.505 &  0.501 &  0.500   &  0.500\\
\hline
\end{tabular}\label{TabConRate}
\end{table}

\noindent\textbf{Test 2.} 
We take the initial mesh as Figure \ref{Mesh_Init}(b), whose elements conforming to the line fracture $\gamma_1$. The convergence rates of the finite element solutions based on $P_1,\,P_2$ polynomials are shown in Table \ref{TabConRateP2}. From the results, we can find that the convergence rates $0.5 < \mathcal R < 2$ depends on the smoothness of the function $g_1$ and  the degree of the polynomials. The results in Table \ref{TabConRateP2} satisfy the theoretical 
expectations shown in Corollary \ref{co12}. 

\begin{table}
\centering
\vspace{-2mm}
\caption{$H^1$ convergence history of finite element solutions in Example \ref{P1h} Test 2 on conforming quasi-uniform meshes.}
\begin{tabular}{|c||c|c|c|c||c|c|c|c|c|c|c|c|}
\hline
&  \multicolumn{3}{c}{$\qquad \qquad P_1$} & &\multicolumn{3}{c}{$ \qquad \qquad P_2$} & \\
\hline
j &6 & 7 & 8 &9 & 5&6 & 7 & 8\\
\hline
Case1 &  0.786 &  0.786 &  0.785 &  0.783&  0.792    &  0.786 &  0.781 &  0.777\\
\hline
Case2 &  0.927 &  0.937 &  0.945 &  0.951 &  1.045   &  1.039 &  1.033 &  1.028\\
\hline
Case3 &  0.905 &  0.916 &  0.925 &  0.932&  1.000    &  1.000    &  1.000    &  1.000\\
\hline
Case4 &  0.969 &  0.979 &  0.986 &  0.990 &  1.253    &  1.252 &  1.251 &  1.251\\
\hline
Case5 &  0.988 &  0.994 &  0.997 &  0.999  &  1.500   &  1.501 &  1.501 &  1.501\\
\hline
Case6  &  0.996 &  0.999 &  1.000    &  1.000 &  1.865 &  1.886 &  1.902 &  1.914\\
\hline
\end{tabular}\label{TabConRateP2}
\end{table}

From the two tests above, we confirm that the finite element solution on the meshes conforming to the line fracture shows better convergence rates than that on meshes with the line fracture passing through the triangles. So we will always consider the initial meshes that conform to line fractures for the remaining examples.

\end{example}

\subsection{Adaptive finite element method}

The parameter $\theta$ in Algorithm \ref{alg1} is taken as $\theta=0.25$ in following examples. On adaptive meshes, the convergence rate of the a posteriori error estimator $\xi$ in (\ref{oldestimator}) or $\eta$ in (\ref{total error estimator}) for $P_k$ polynomials is called quasi-optimal if
$$
\xi \approx N^{-0.5k}, \quad \text{or} \quad \eta \approx N^{-0.5k}.
$$
Here and in what follows, we abuse the notation $N$ to represent the total number of degrees of freedom. 

\begin{example}\label{example3}
We apply the AFEM to the Example \ref{P1h} to test the performance of the proposed a posterior error estimator (\ref{total error estimator}) and the corresponding Algorithm \ref{alg1}.
We take the mesh in Figure \ref{Mesh_Init}(b) as the initial mesh.
The convergence rates of the error estimator $\eta$ based on $P_1$ and $P_2$ polynomials are shown
Figure \ref{fig:errors exam3}. From the results, we find that the convergence rates of $\eta$ are quasi-optimal. The contours of the AFEM approximations for different cases are shown in Figure \ref{fig:solution1 exam3}, from which we can find that these solutions are almost identical to these in Example \ref{P1h} Test 1.

For Case 1$-$6, the function $g_1$ is sufficiently smooth on $\gamma_1$ except near the endpoints $Q_1$ and $Q_2$ of the line fracture $\gamma_1$, so the solution is more singular near these two endpoints compared with any other regions in the domain.
Figure \ref{fig:mesh1 exam3} and Figure \ref{fig:mesh2 exam3} show the adaptive meshes of $P_1,\,P_2$ approximations, respectively. We can see clearly that the error estimator guide the mesh refinements densely around the endpoints $Q_1$ and $Q_2$. 
We also find that the more regular the solution is, the less dense the mesh concentrates at the endpoints $Q_1$ and $Q_2$.
Here, Case 3 is an example in \cite{LWYZ21} solved by the graded finite element method, which showed optimal convergence rates with mesh refinements concentrating at the singular points $Q_1$ and $Q_2$ as well. 

\begin{figure}
\centering
\subfigure[Case 1]
{\includegraphics[width=0.32\textwidth]{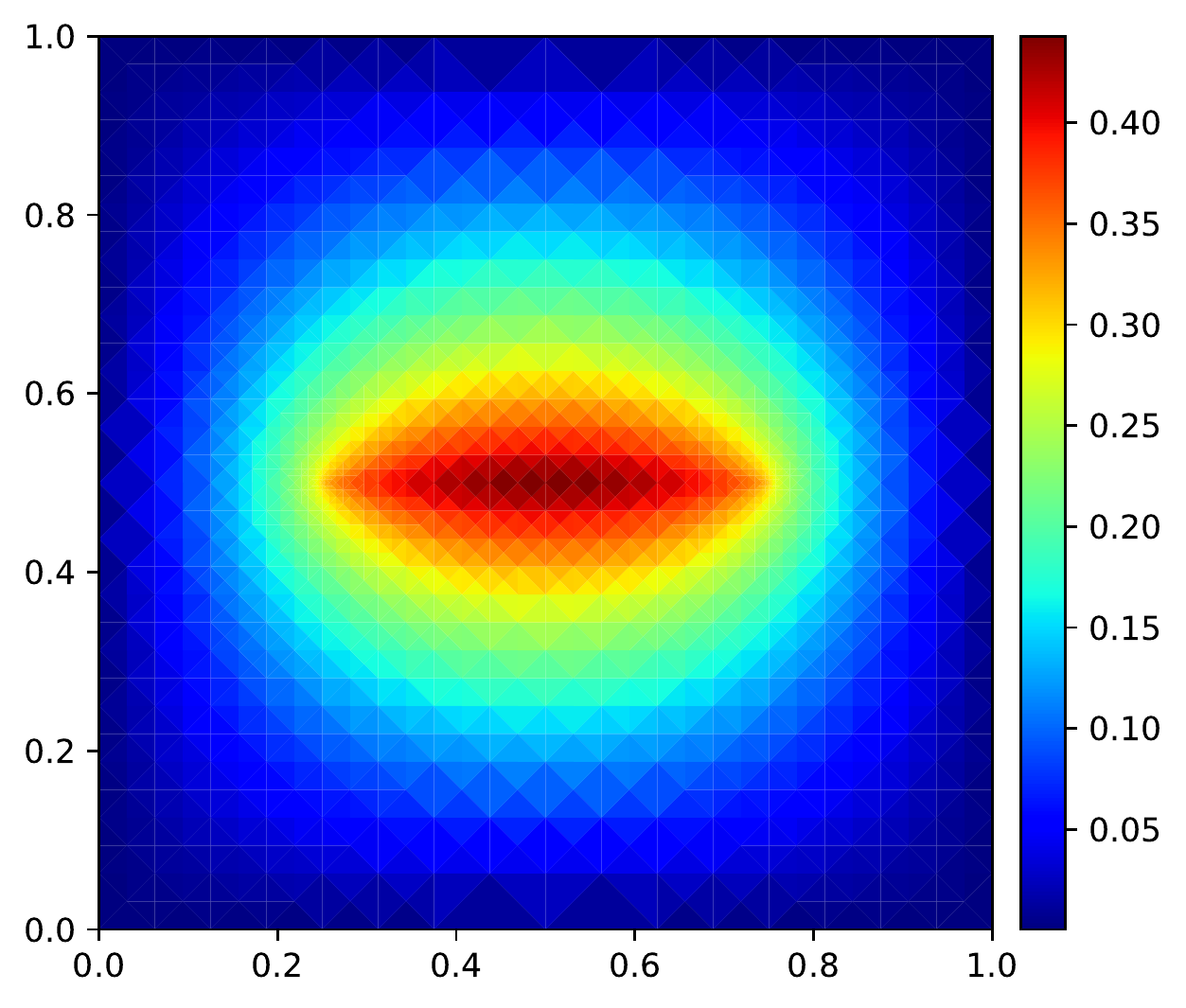}}
\subfigure[Case 2]
{\includegraphics[width=0.32\textwidth]{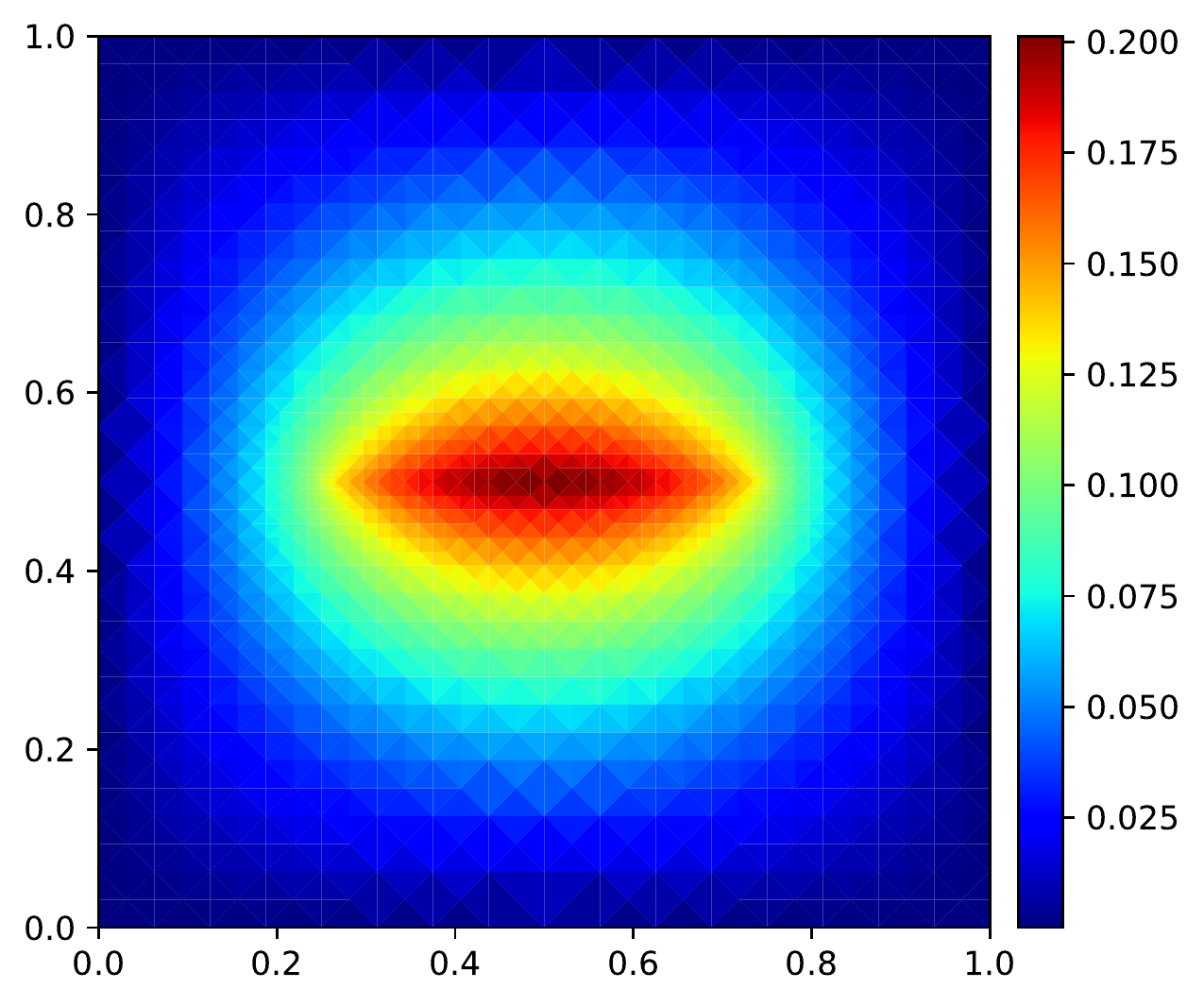}}
\subfigure[Case 3]
{\includegraphics[width=0.32\textwidth]{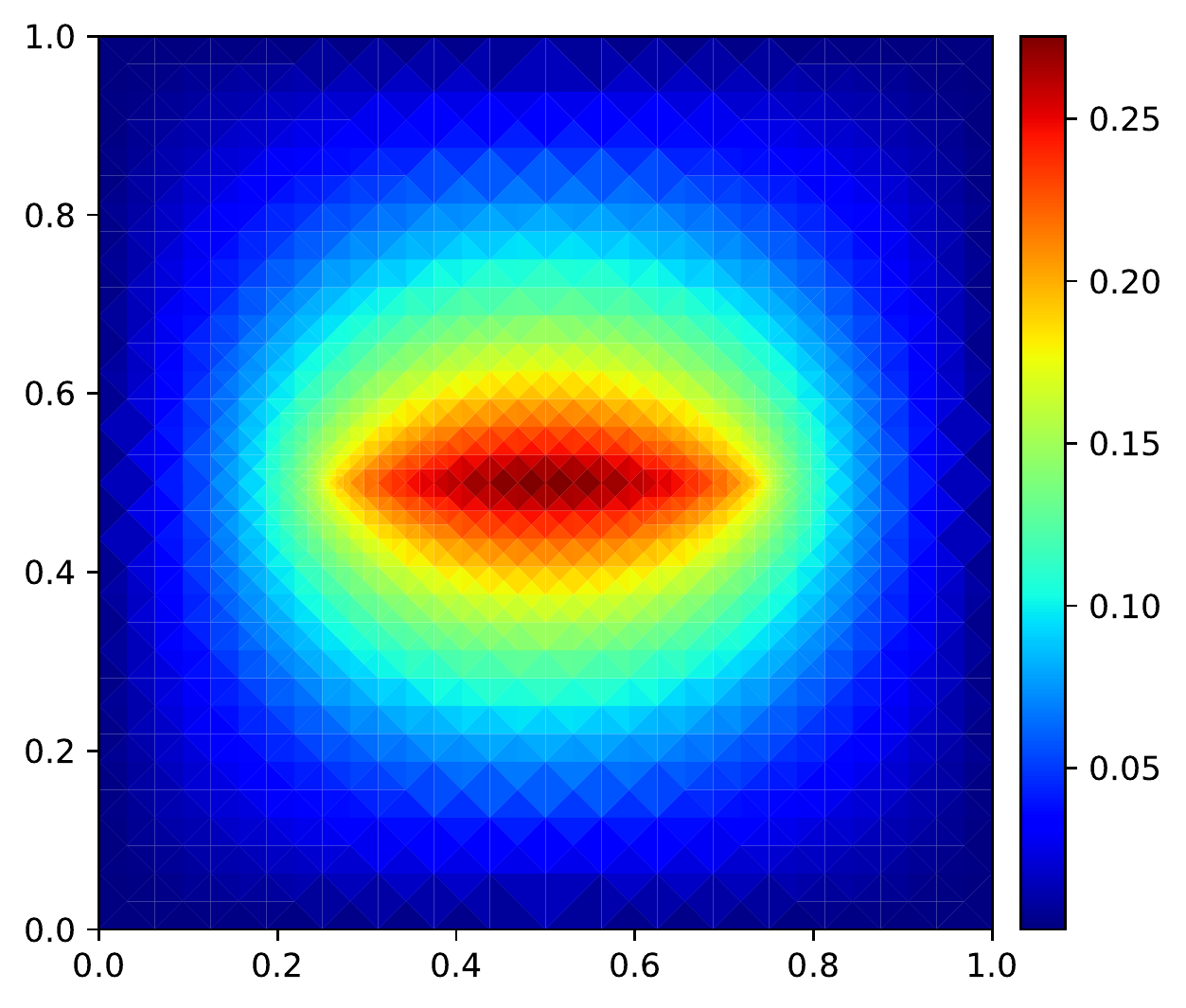}}
\subfigure[Case 4]
{\includegraphics[width=0.32\textwidth]{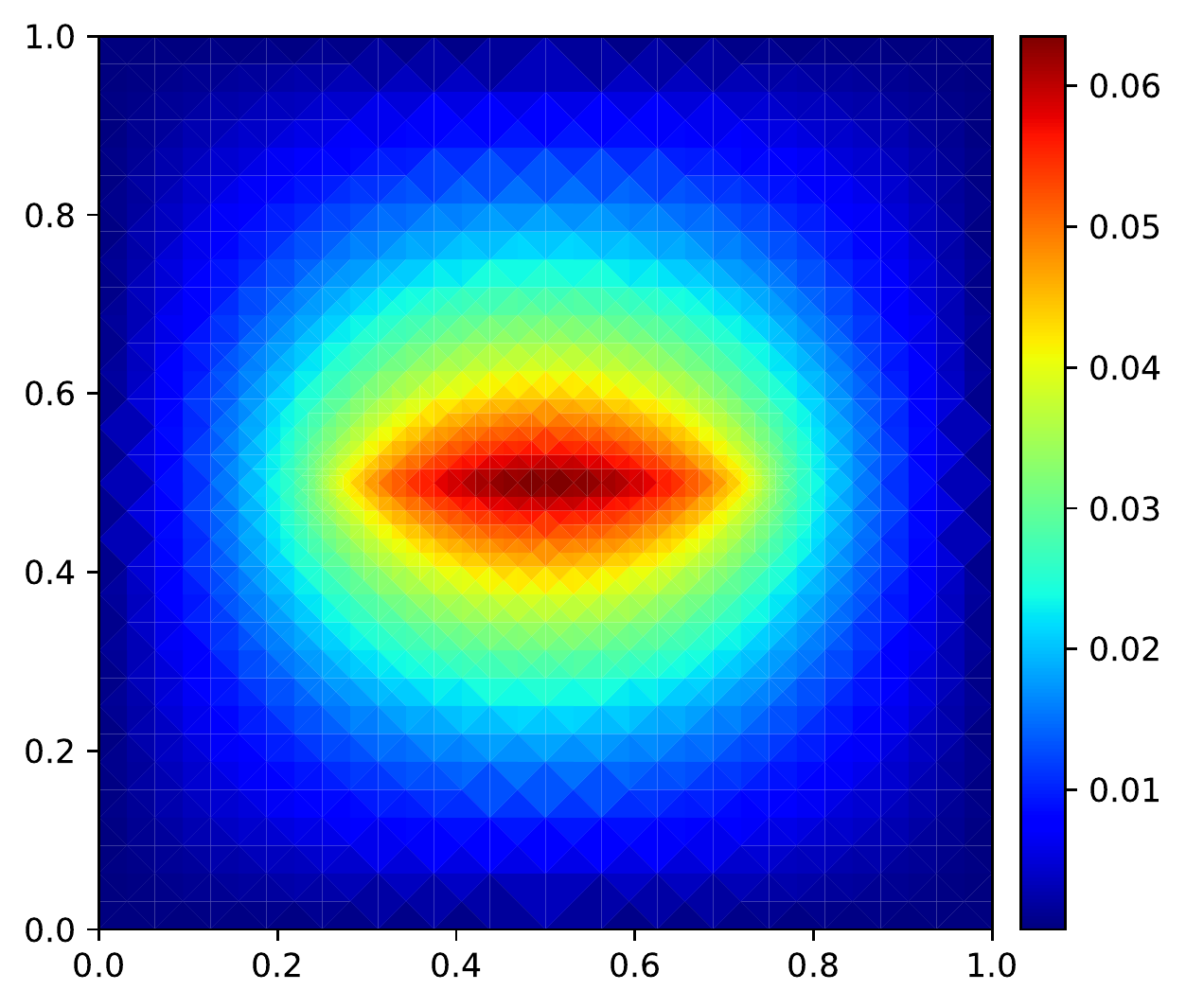}}
\subfigure[Case 5]
{\includegraphics[width=0.32\textwidth]{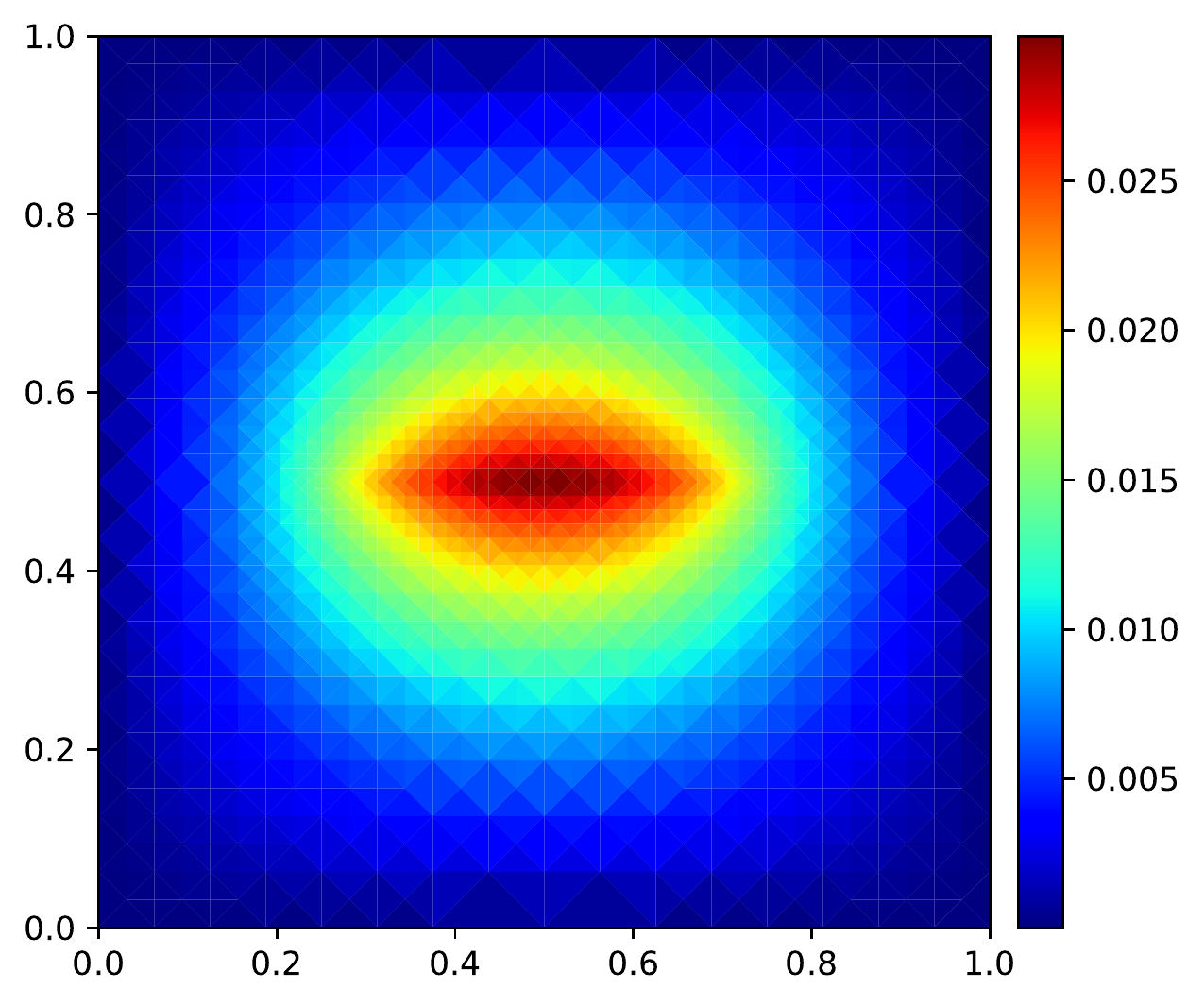}}
\subfigure[Case 6]
{\includegraphics[width=0.32\textwidth]{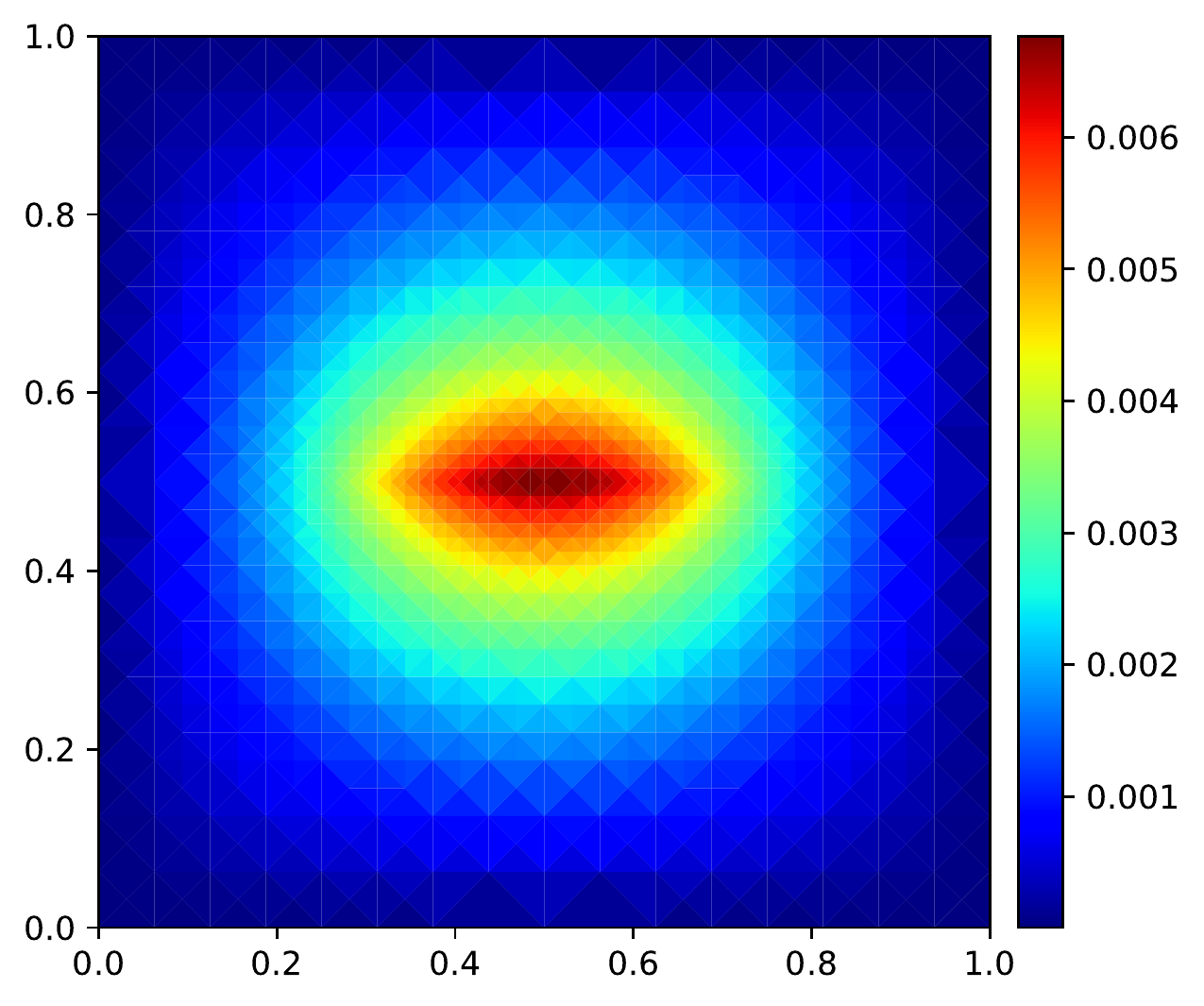}}
\vspace{-2mm}
\caption{Example \ref{example3}: AFEM solutions based on $P_1$ polynomials.}\label{fig:solution1 exam3}
\end{figure}

\begin{figure}
\centering
\subfigure[$P_1$]
{\includegraphics[width=0.35\textwidth]{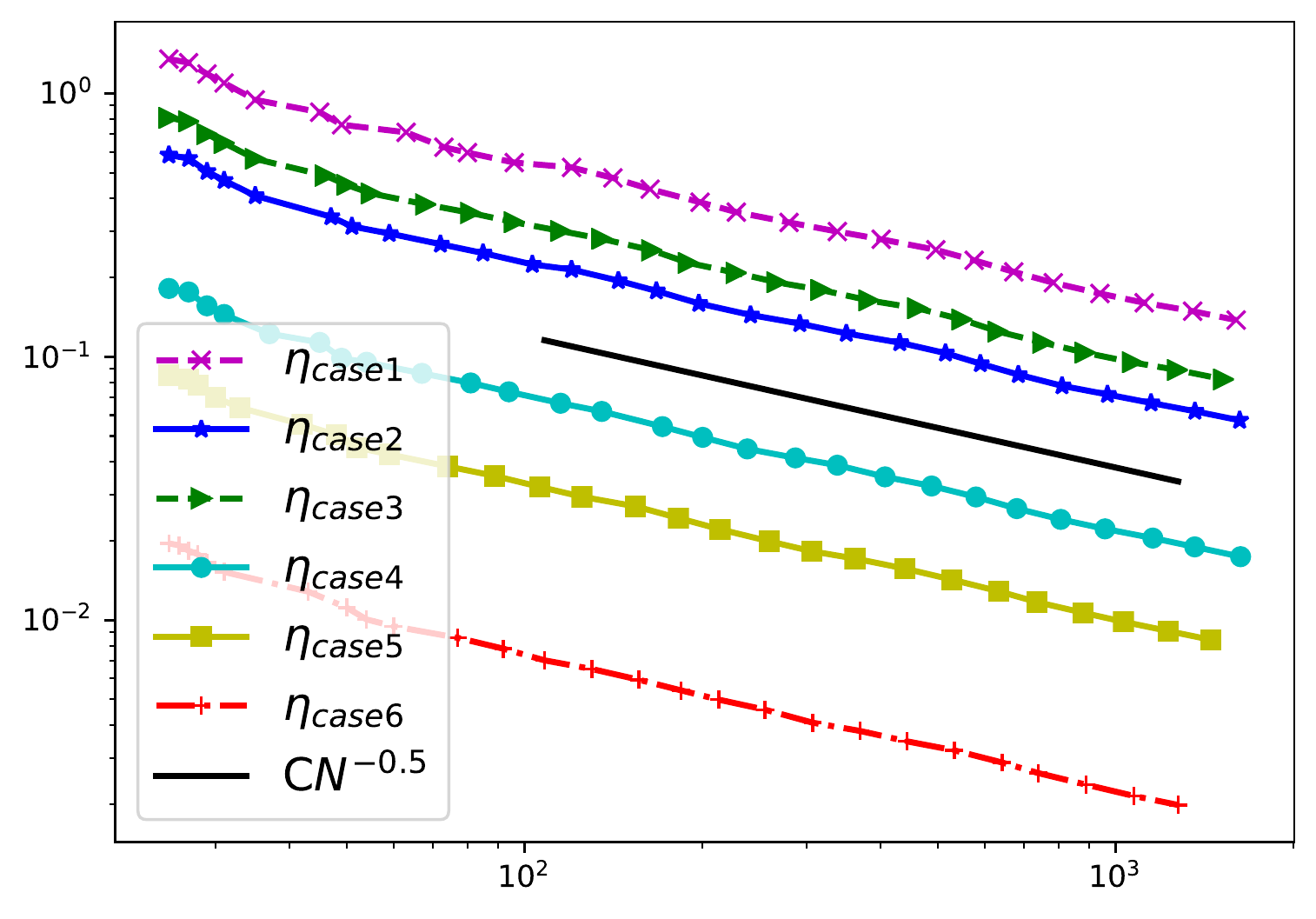}}\hspace{8mm}
\subfigure[$P_2$]
{\includegraphics[width=0.35\textwidth]{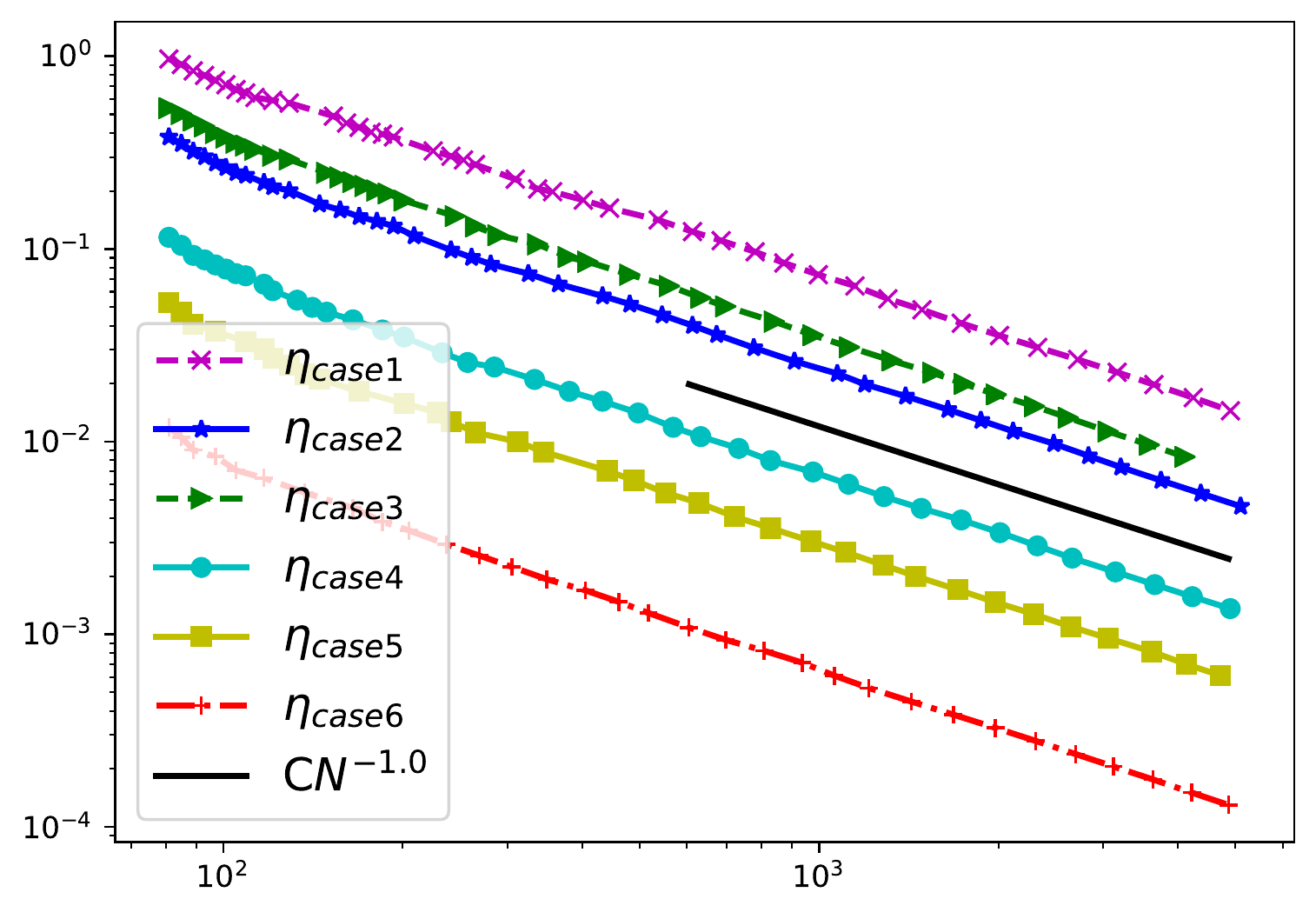}}
\vspace{-2mm}
\caption{Example \ref{example3}: error estimators.}\label{fig:errors exam3}
\end{figure}

\begin{figure}
\centering
\subfigure[Case 1]
{\includegraphics[width=0.29\textwidth]{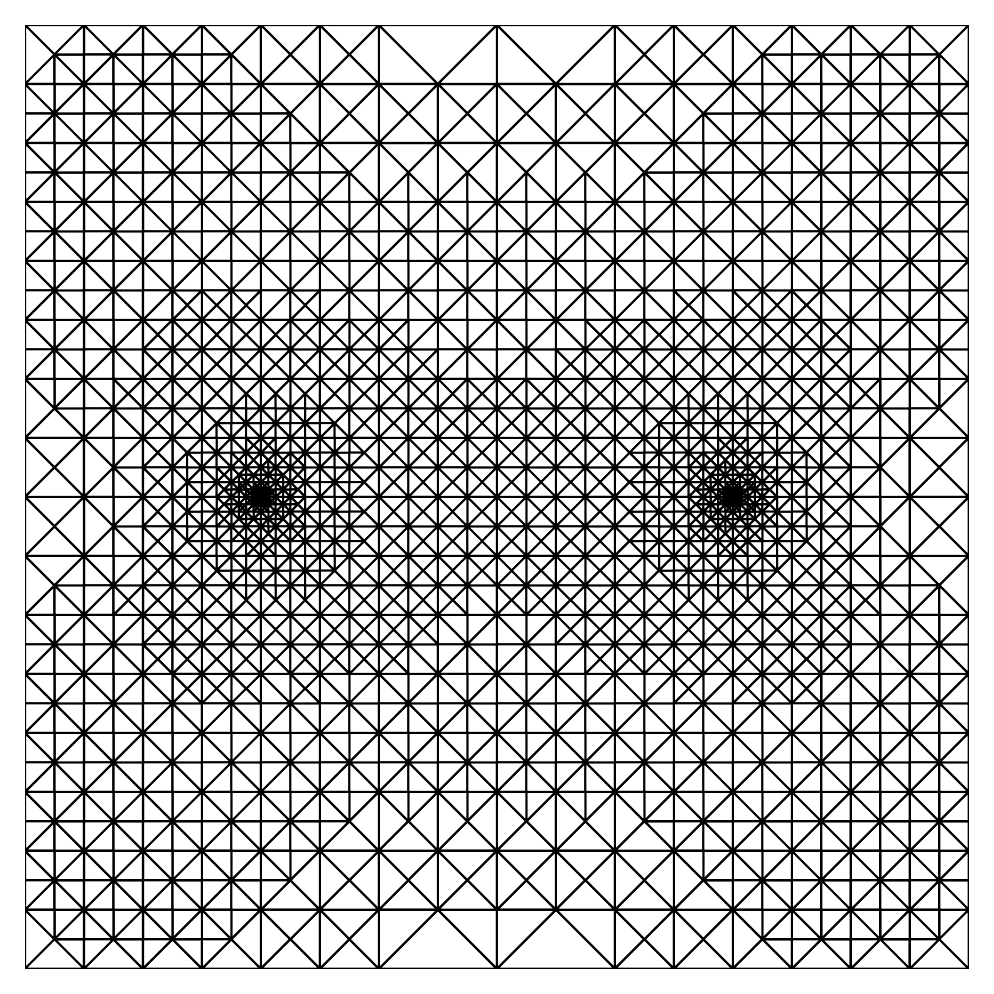}}\hspace{2mm}
\subfigure[Case 2]
{\includegraphics[width=0.29\textwidth]{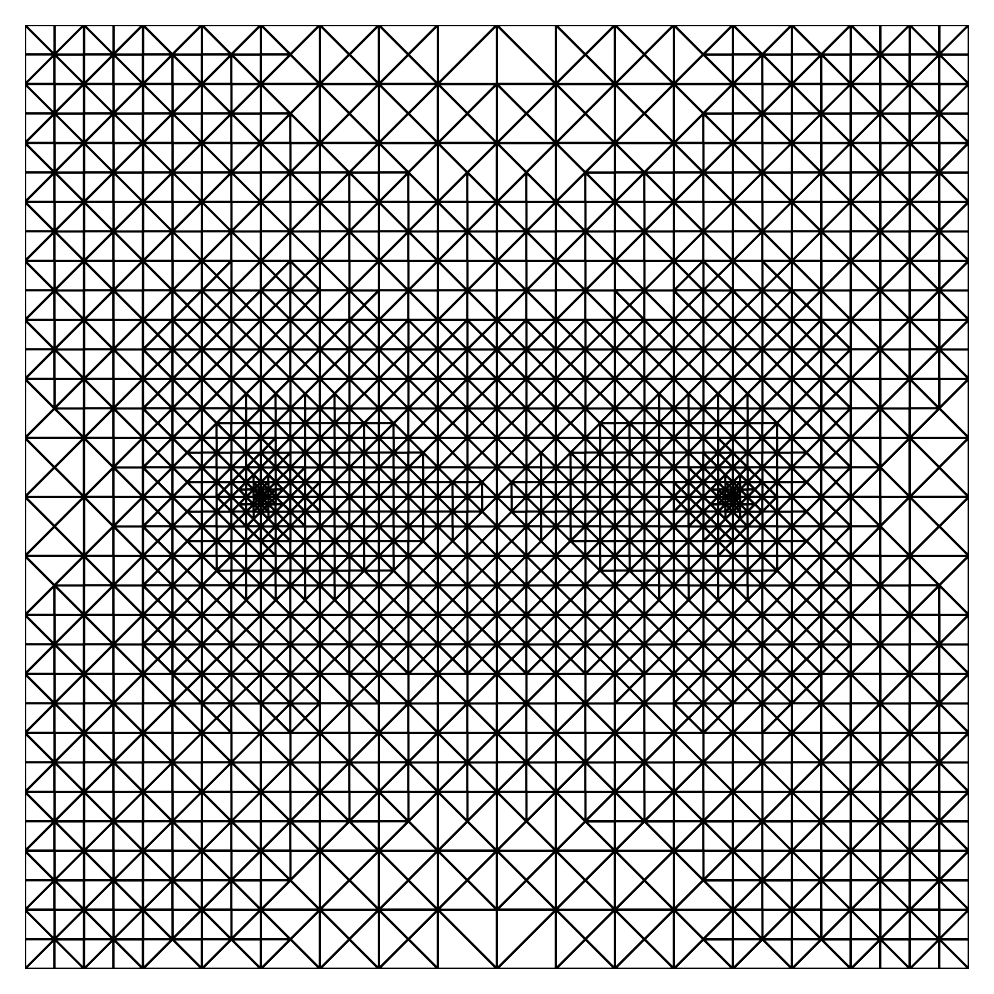}}\hspace{2mm}
\subfigure[Case 3]
{\includegraphics[width=0.29\textwidth]{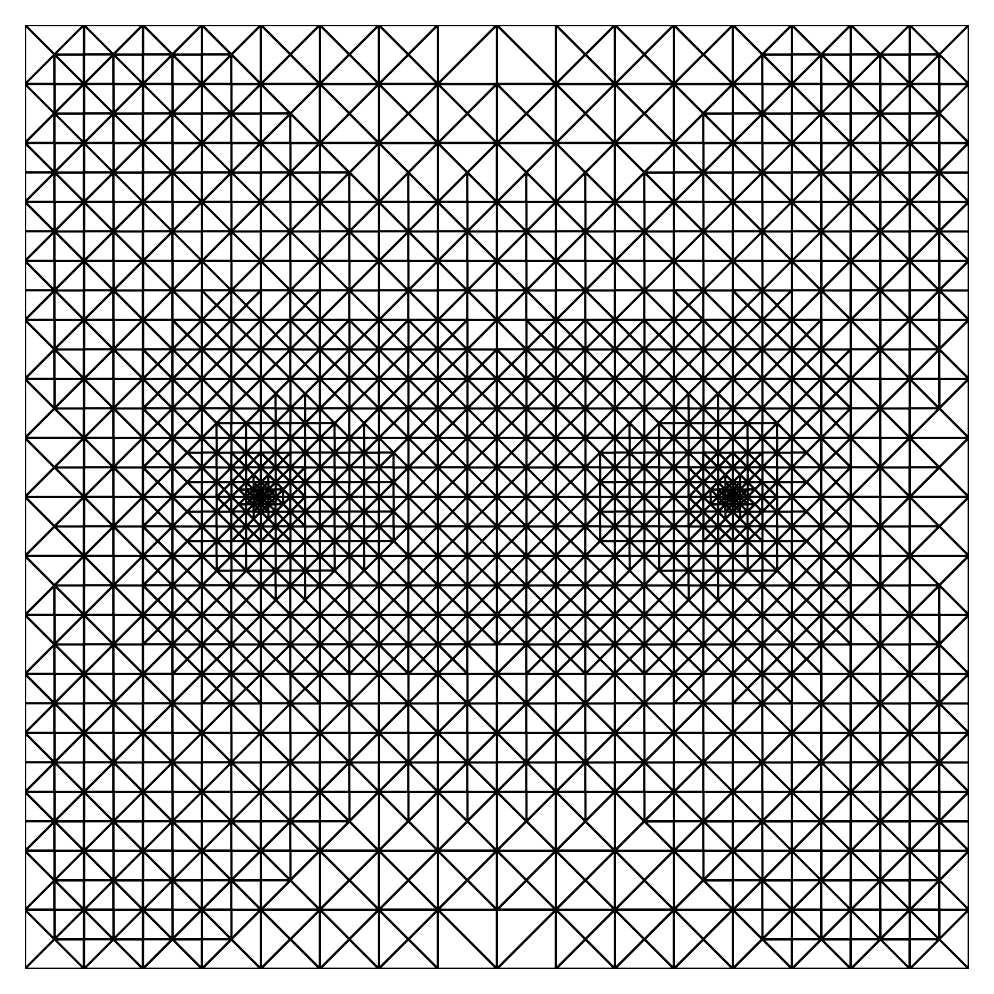}}\hspace{2mm}
\subfigure[Case 4]
{\includegraphics[width=0.29\textwidth]{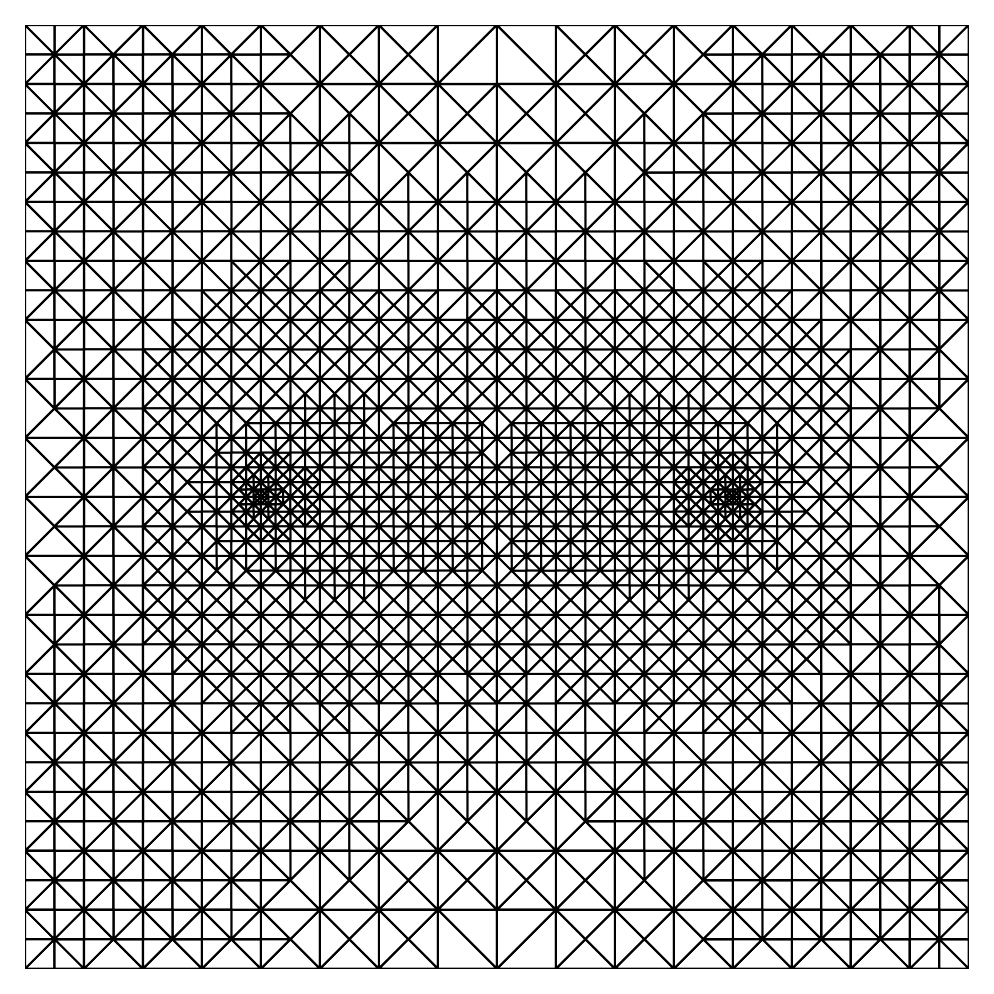}}\hspace{2mm}
\subfigure[Case 5]
{\includegraphics[width=0.29\textwidth]{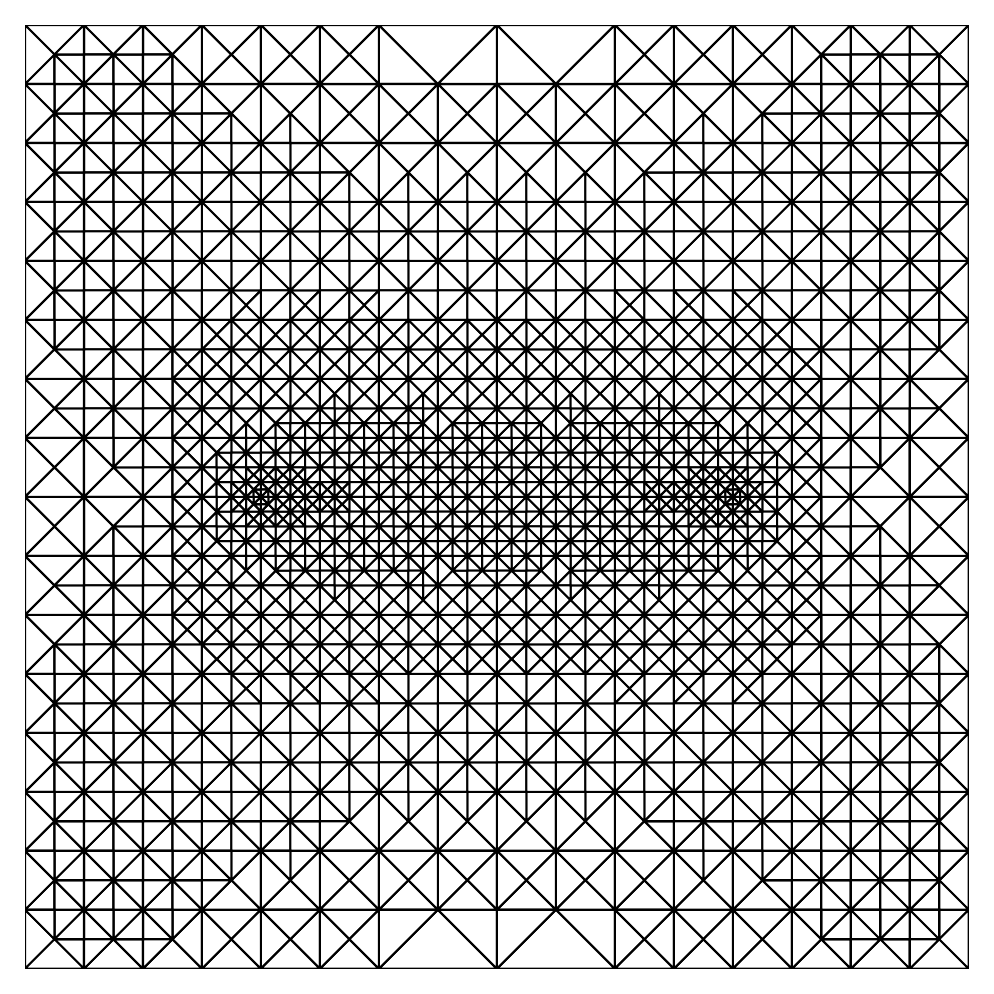}}\hspace{2mm}
\subfigure[Case 6]
{\includegraphics[width=0.29\textwidth]{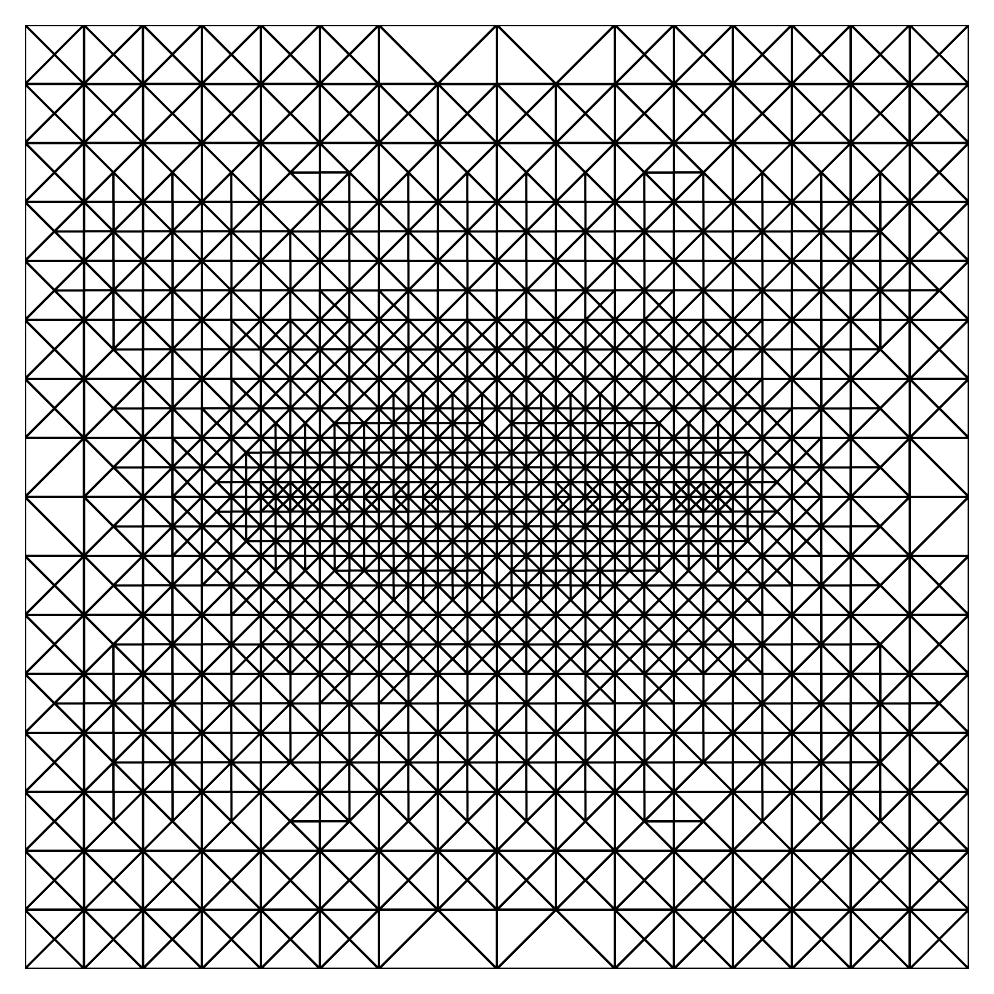}}
\vspace{-2mm}
\caption{Example \ref{example3}: adaptive meshes based on $P_1$ polynomials.}\label{fig:mesh1 exam3}
\end{figure}

\begin{figure}
\centering
\subfigure[Case 1]
{\includegraphics[width=0.3\textwidth]{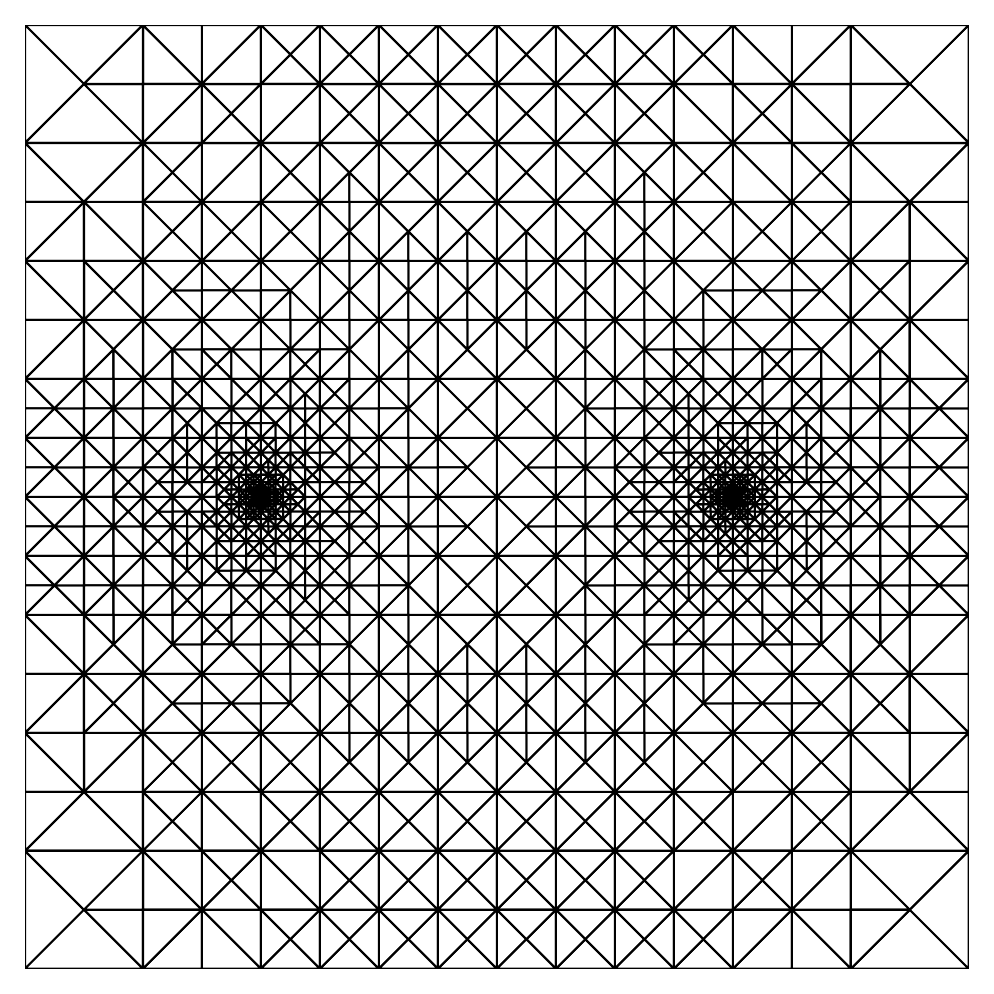}}\hspace{2mm}
\subfigure[Case 2]
{\includegraphics[width=0.3\textwidth]{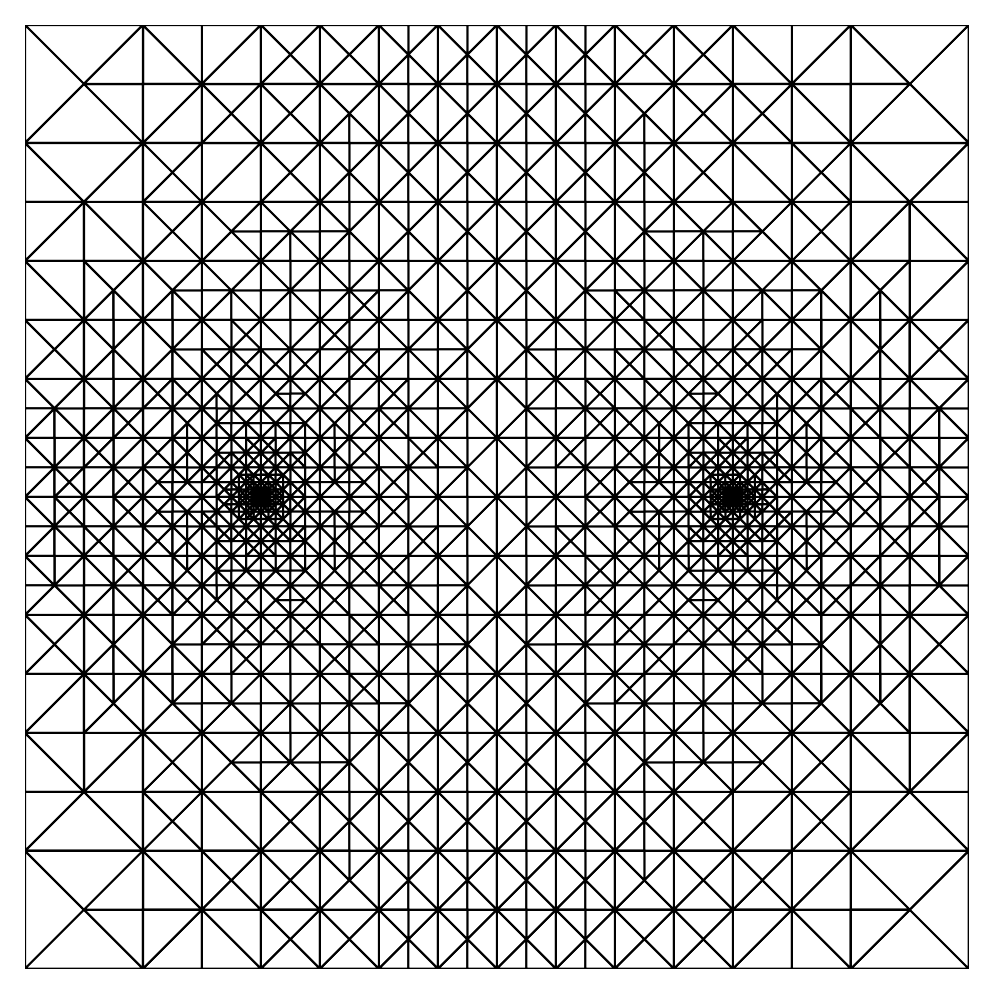}}\hspace{2mm}
\subfigure[Case 3]
{\includegraphics[width=0.3\textwidth]{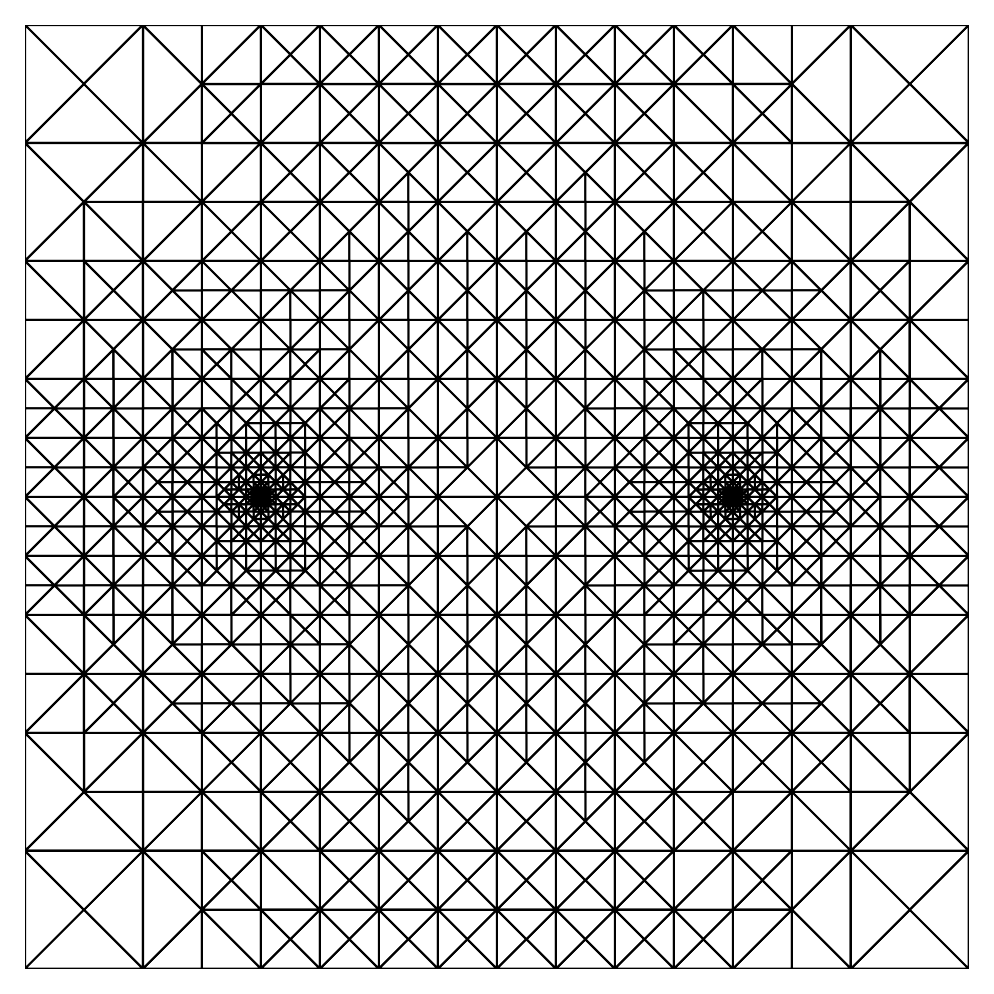}}\hspace{2mm}
\subfigure[Case 4]
{\includegraphics[width=0.3\textwidth]{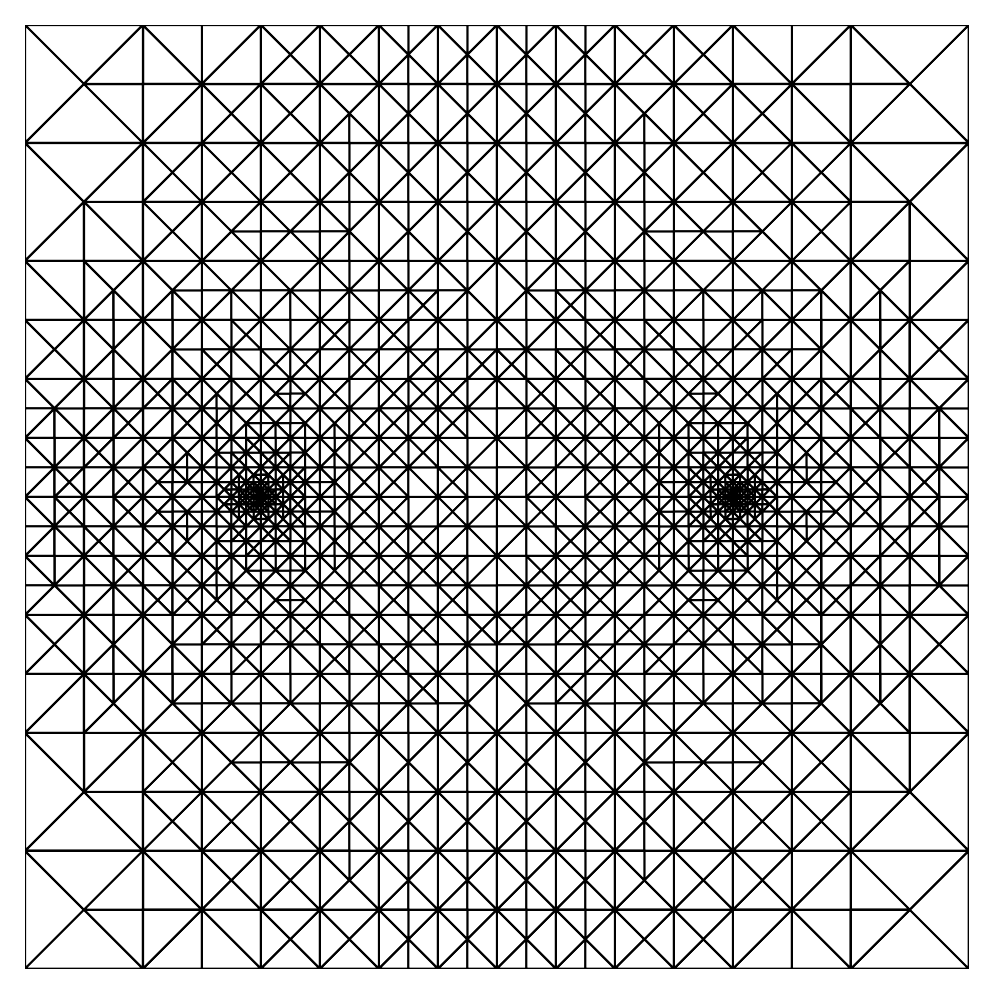}}\hspace{2mm}
\subfigure[Case 5]
{\includegraphics[width=0.3\textwidth]{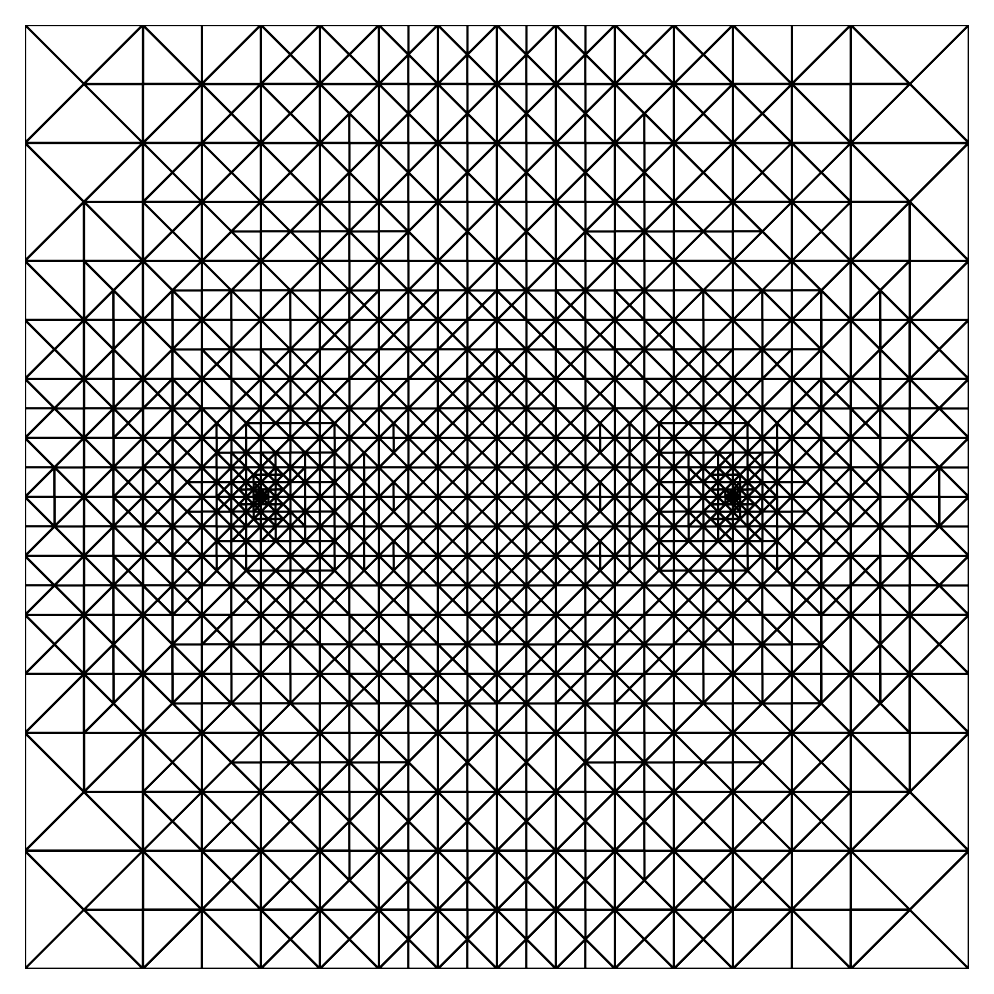}}\hspace{2mm}
\subfigure[Case 6]
{\includegraphics[width=0.3\textwidth]{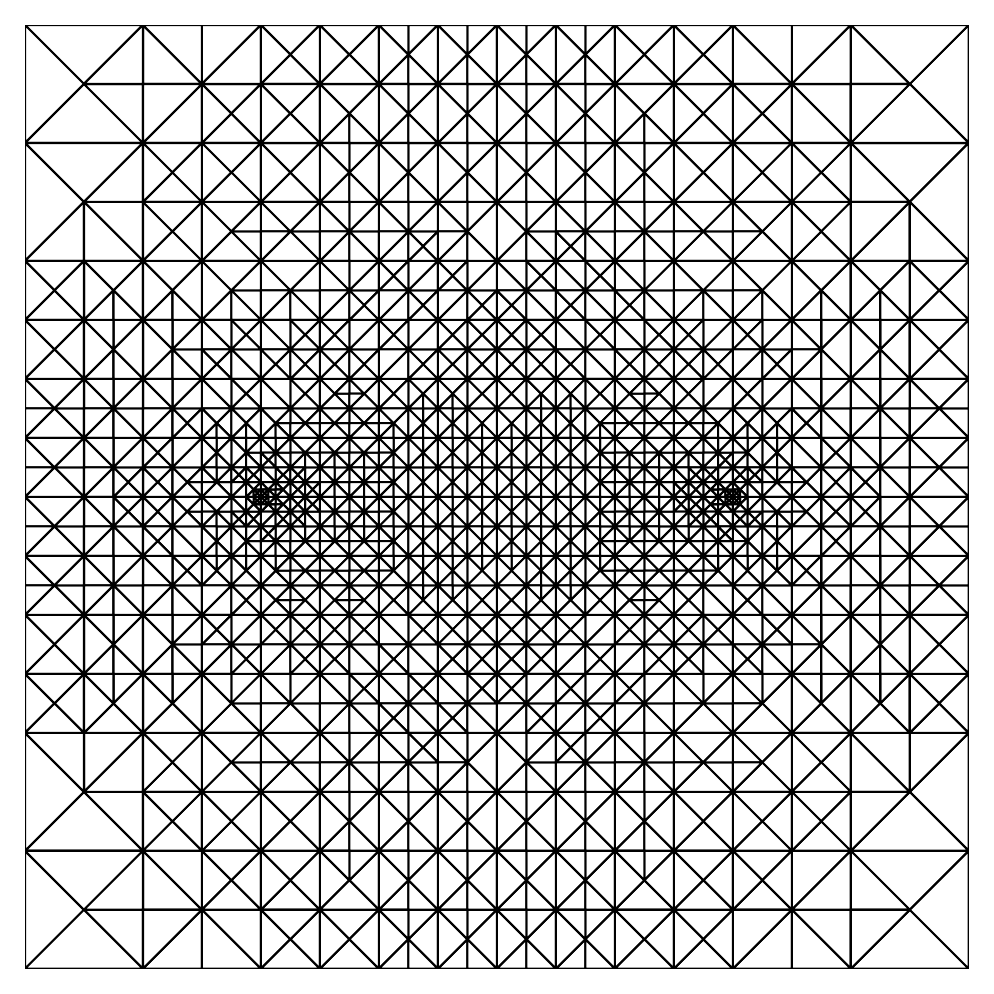}}
\vspace{-2mm}
\caption{Example \ref{example3}: adaptive meshes based on $P_2$ polynomials.}\label{fig:mesh2 exam3}
\end{figure}
\end{example}

\begin{example}\label{example2} 
We take this example from \cite{HL21}. More specifically, we consider problem \eqref{eq:Possion} on an L-shaped domain $\Omega = (-1,1)^2 \backslash [0,1)^2$ and take the line fractures $\cup_{l=1}^6\gamma_l = \partial \Omega_1 $ with $\Omega_1 = (-0.8,-0.2)^2 \backslash [-0.5,-0.2)^2$ as shown in Figure \ref{fig:meshes exam2}(a). The function $g_l = 5$ on $\gamma_l,\,l = 1,\cdots,6$. We apply the AFEMs based on the residual-based a posteriori error estimators $\xi$ in (\ref{oldestimator}) and $\eta$ in (\ref{total error estimator}) to solve this problem, respectively. Both AFEMs take the mesh in Figure \ref{fig:meshes exam2}(a) as their initial mesh.

\noindent\textbf{Test 1.} We first consider the AFEM based on the residual-based a posteriori error estimator $\xi$ in (\ref{oldestimator}). 
For simplicity of presentation, we denote $\gamma = \cup_{l=1}^6\gamma_l$, and $g|_{\gamma_l} = g_l$, $l=1,\cdots, 6$. Instead of directly discretizing (\ref{eq:Possion}), one discretize its regularized problem, which is to replace the line Dirac source term $\sum_{l=1}^N g_l\delta_{\gamma_l}$ by its regularized data \cite{HL21},
\begin{align*}
g^r(x) = \int_{\gamma} g(y)\delta^r(y-x)\,\mathrm{d}y \in L^2(\Omega).
\end{align*}
Here, the line Dirac approximation $\delta^r$ of the 2-dimensional Dirac distribution $\delta$ is defined by
\begin{align*}
\delta^r(x) = \frac{1}{r^2}\psi\left(\frac{x}{r}\right),
\end{align*}
satisfying
\begin{align*}
\lim_{r \to 0} \delta^r(x) = \lim_{r \to 0}\frac{1}{r^2}\psi\left(\frac{x}{r}\right) = \delta(x),
\end{align*}
where $r$ is the regularization parameter depending on the local mesh size, and $\psi(x)$ is the Dirac approximation \cite{HNS16,T2002,HL21}. Here, we take $r=0.05$, and 
\begin{align*}
\psi(x) = \frac{1}{4}\Pi_{i=1}^2\chi_{[-1,1]}(x_i),
\end{align*}
in which $\chi_{[-1,1]}(x_i)$ is the characteristic function. 
The contour of the finite element solution based on $P_1$ polynomials is shown in Figure \ref{fig:meshes exam2}(b).

\noindent\textbf{Test 2.} We then consider the AFEM based on the residual-based a posteriori error estimator $\eta$ in (\ref{total error estimator}), namely, the Algorithm \ref{alg1}, for problem \eqref{eq:Possion}. The contour of the finite element solution based on $P_1$ polynomials is shown in Figure \ref{fig:meshes exam2}(c), which is comparable to the contour in Test 1 as shown in Figure \ref{fig:meshes exam2}(b).

Since $g_l \in C^\infty$ are sufficiently smooth on line fractures $\gamma_l$, so the solution is more singular at the endpoints of line fractures $\gamma_l$ and the reentrant corner of the domain.
The adaptive meshes from Test 1 and Test 2 based on $P_1$ polynomials are shown in Figure \ref{fig: exam2 P1}(a) and Figure \ref{fig: exam2 P1}(b), respectively. From the results, we find that both meshes are densely refined at the endpoints of the line fractures $\gamma_l$ and the reentrant corner of the domain, but the mesh from Test 1 is also densely refined on the whole line fractures $\gamma_l$, $l=1,\cdots,6$. 
Similar adaptive meshes can also be found for Test 1 and Test 2 based on $P_2$ polynomials as shown in Figure \ref{fig: exam2 P2}(a)-(b).
These results imply that the error estimator $\eta$ in (\ref{total error estimator}) guides the mesh refinements effectively by only densely refining the triangles around the endpoints of the line fractures, where the solution is more singular.

The convergence rates of the error estimator $\xi$ and $\eta$ based on $P_1$ polynomials are shown in Figure \ref{fig: exam2 P1}(c). We can find that the error estimators from both Test 1 and Test 2 are quasi-optimal with $\xi \approx N^{-0.5}$ and $\eta \approx N^{-0.5}$.
The convergence rates based on $P_2$ polynomials are shown in Figure \ref{fig: exam2 P2}(c). From the results, we can find that the error estimator $\eta \approx N^{-1}$ for Test 2 is quasi-optimal, but the error estimator $\xi (\approx N^{-0.5})$ for Test 1 does not achieve the quasi-optimal rate even with more dense refined meshes.

\begin{figure}
\centering
\subfigure[Domain with initial mesh]
{\includegraphics[width=0.28\textwidth]{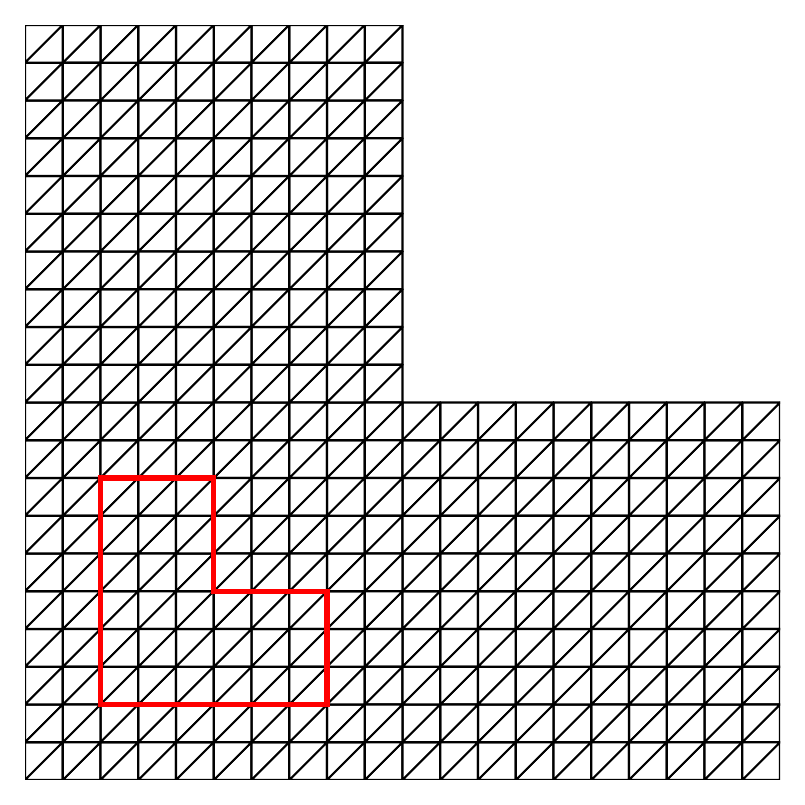}}
\subfigure[Solution of Test 1]
{\includegraphics[width=0.35\textwidth]{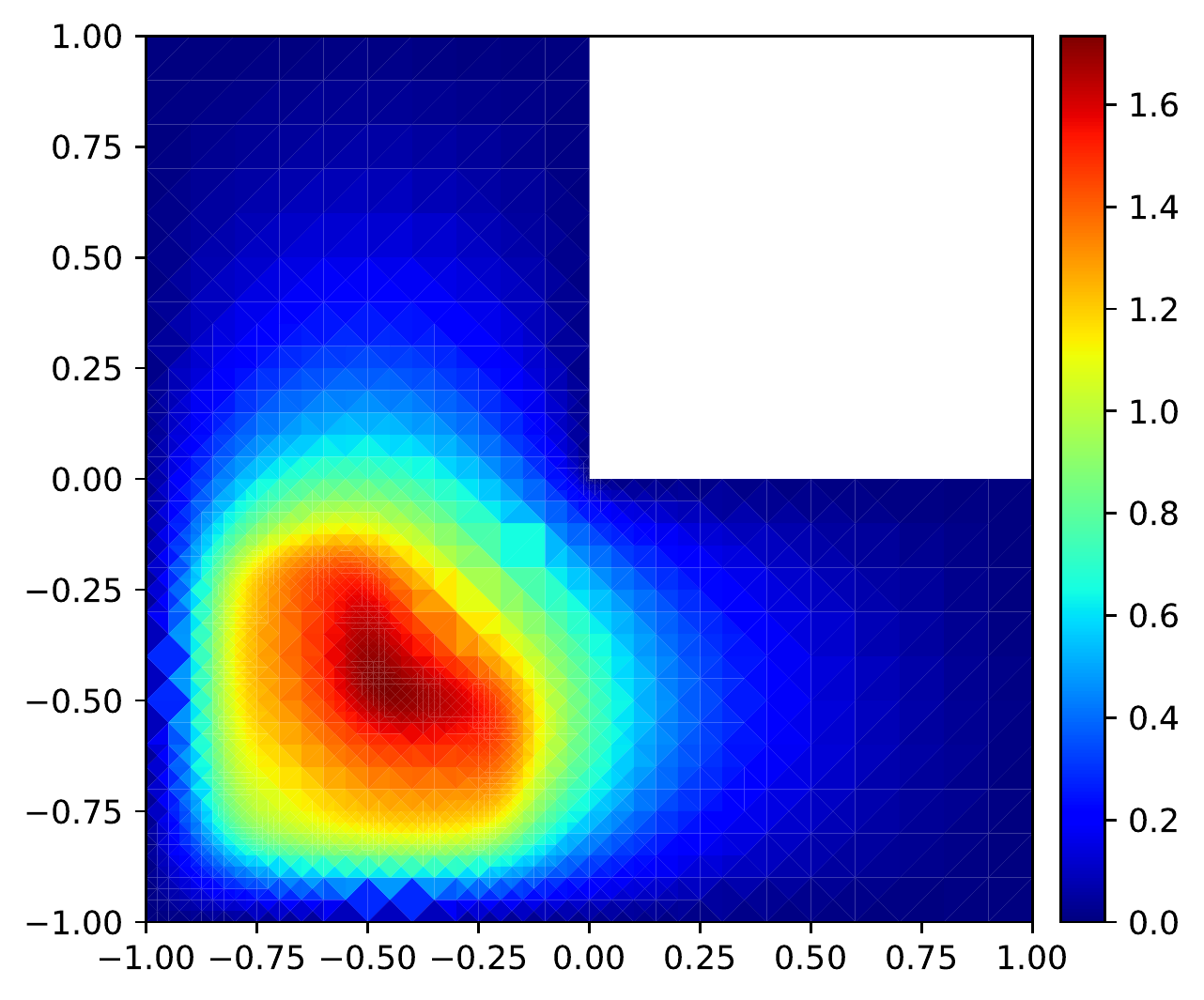}}
\subfigure[Solution of Test 2]
{\includegraphics[width=0.35\textwidth]{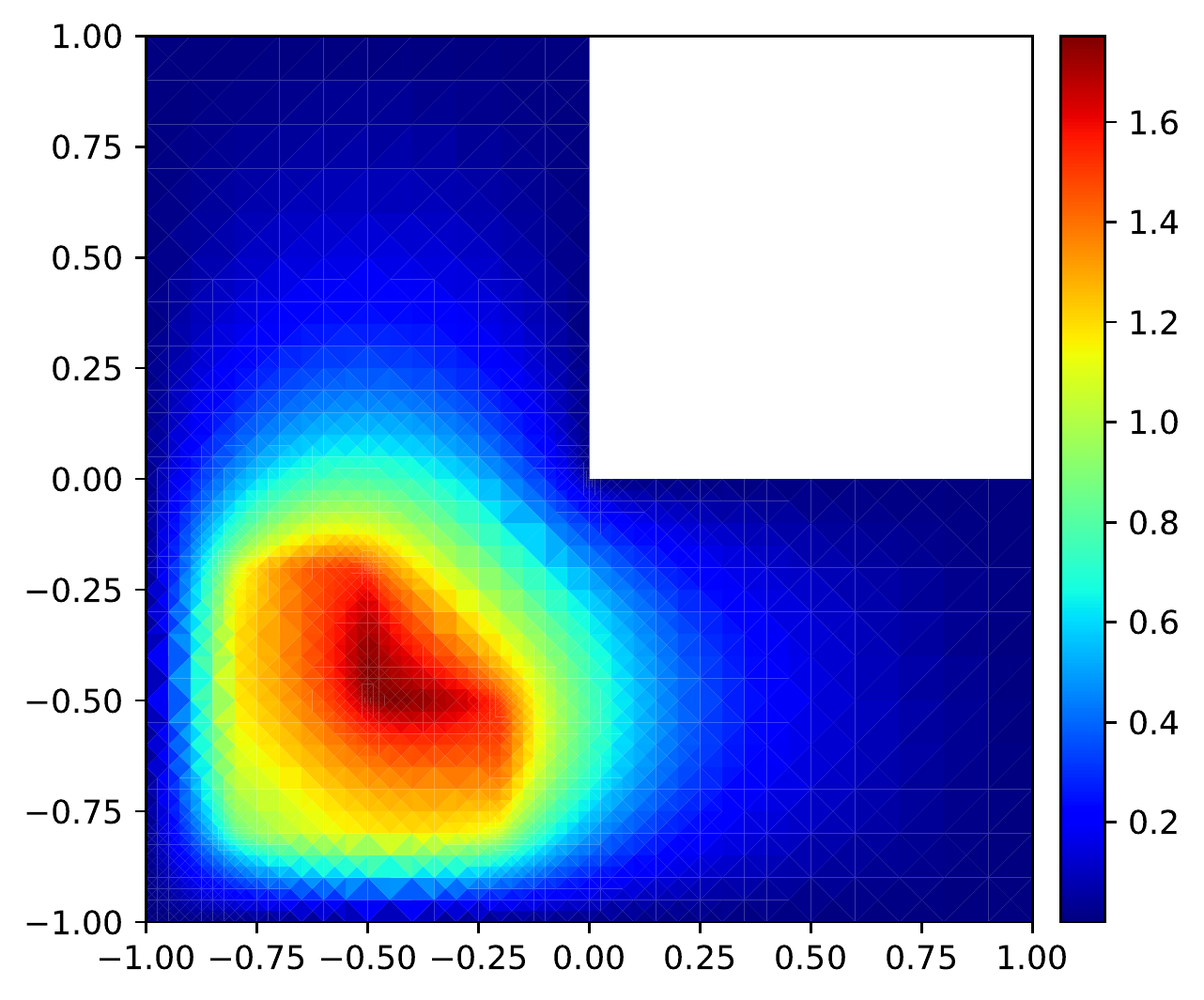}}
\vspace{-2mm}
\caption{Example \ref{example2}: Initial mesh and finite element solutions.}\label{fig:meshes exam2}
\end{figure}

\begin{figure}
\centering
\subfigure[Test 1: adaptive mesh based on regularization]
{\includegraphics[width=0.3\textwidth]{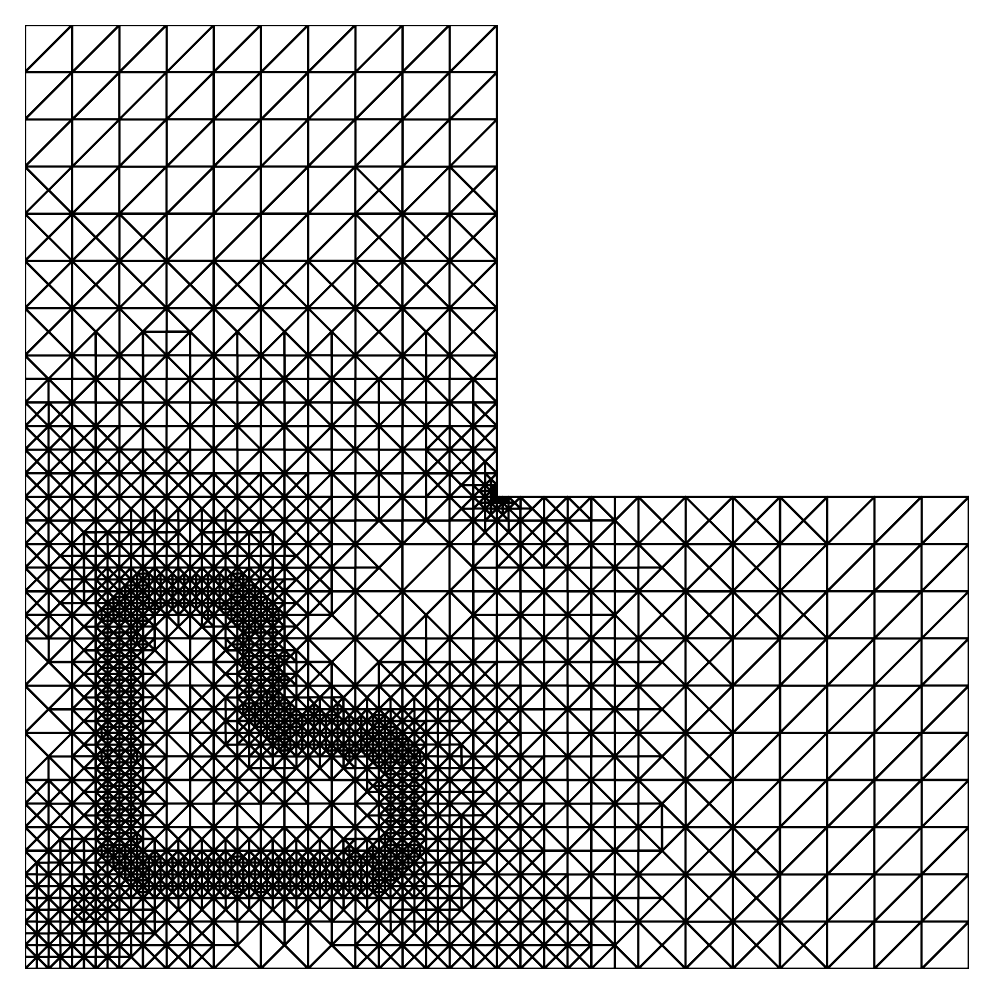}}
\subfigure[Test 2: adaptive mesh based on Algorithm \ref{alg1}]
{\includegraphics[width=0.3\textwidth]{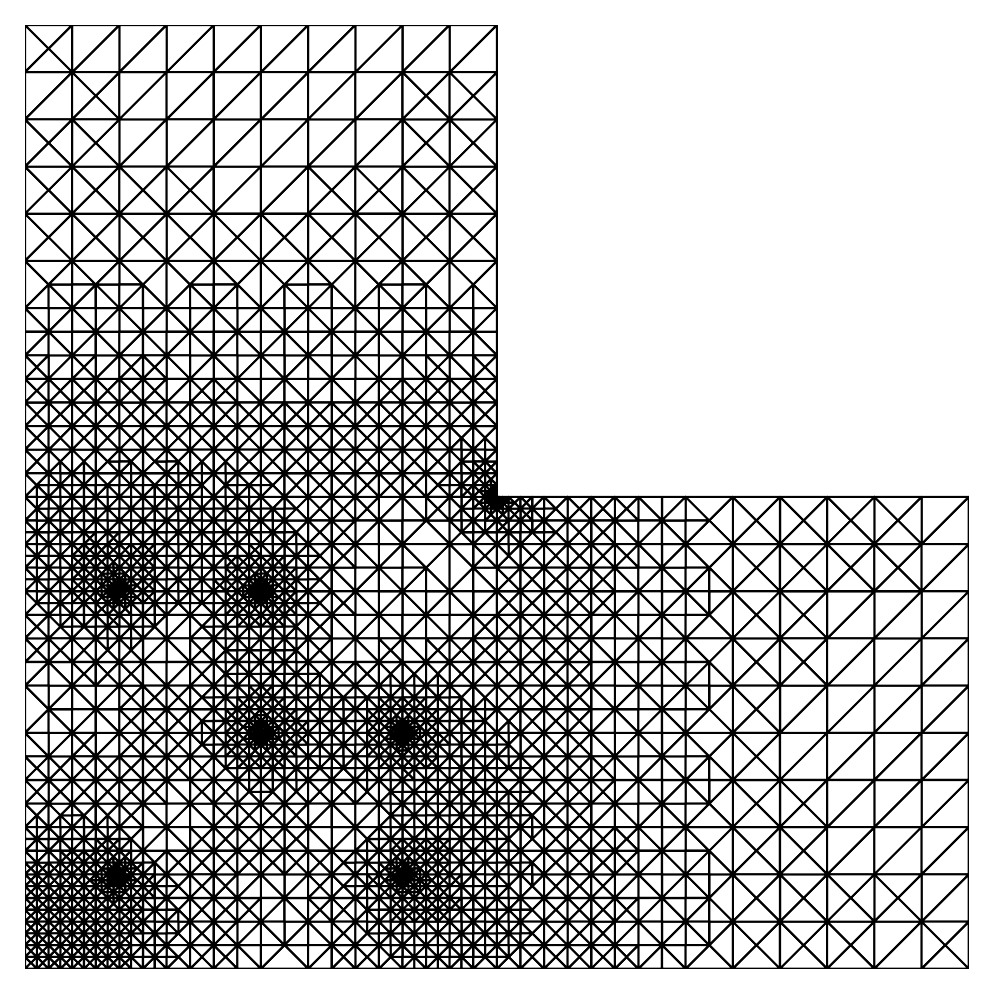}}
\subfigure[Error estimators]
{\includegraphics[width=0.37\textwidth]{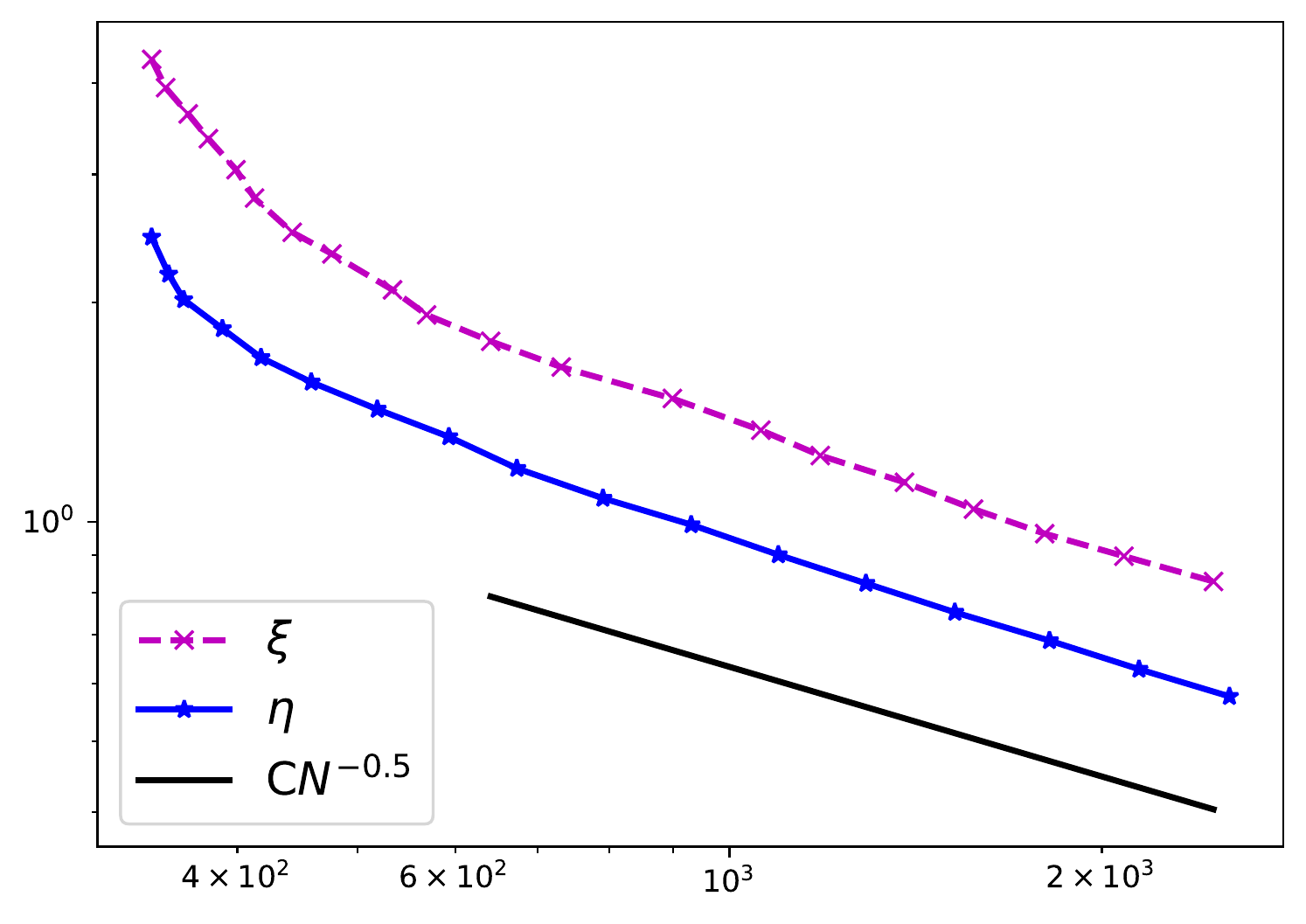}}
\vspace{-2mm}
\caption{Example \ref{example2}: Adaptive meshes and error estimators based on $P_1$ polynomials.}\label{fig: exam2 P1}
\end{figure}
\begin{figure}
\centering
\subfigure[Test 1: adaptive mesh based on regularization]
{\includegraphics[width=0.3\textwidth]{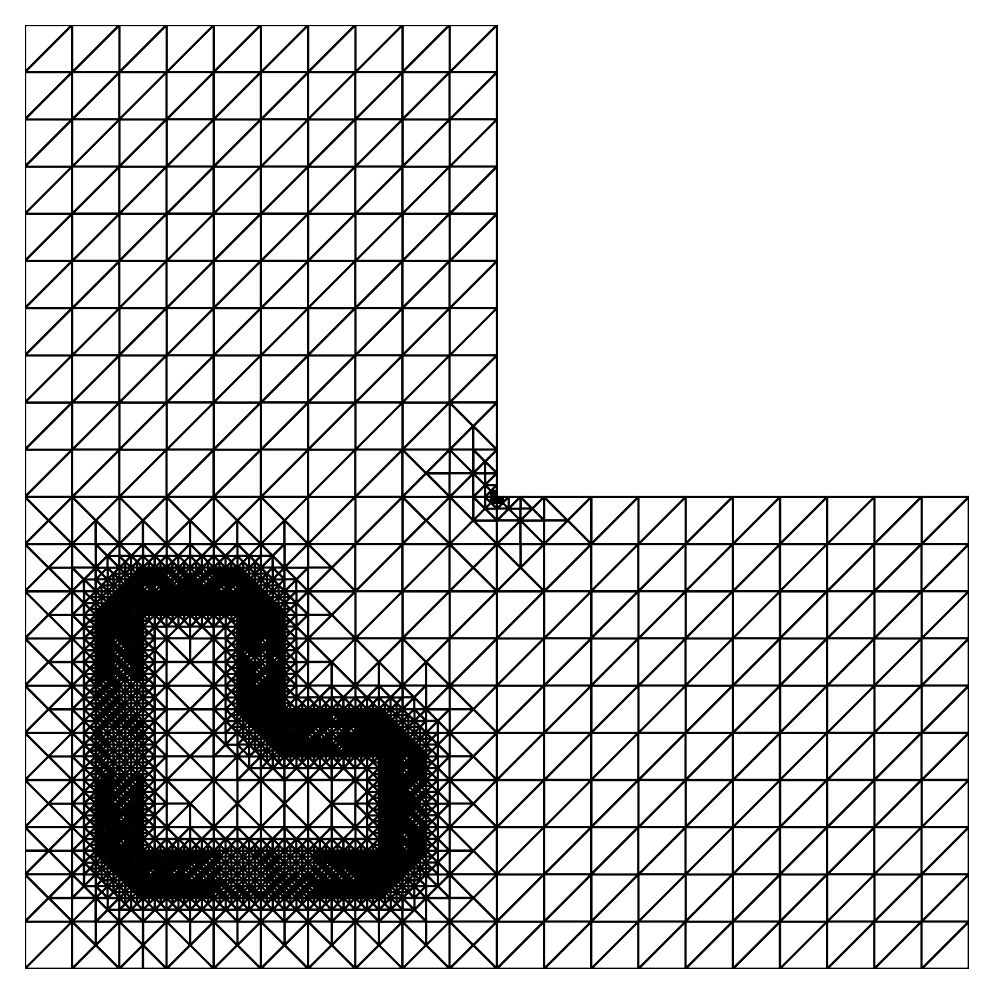}}
\subfigure[Test 2: adaptive mesh based on Algorithm \ref{alg1}]
{\includegraphics[width=0.3\textwidth]{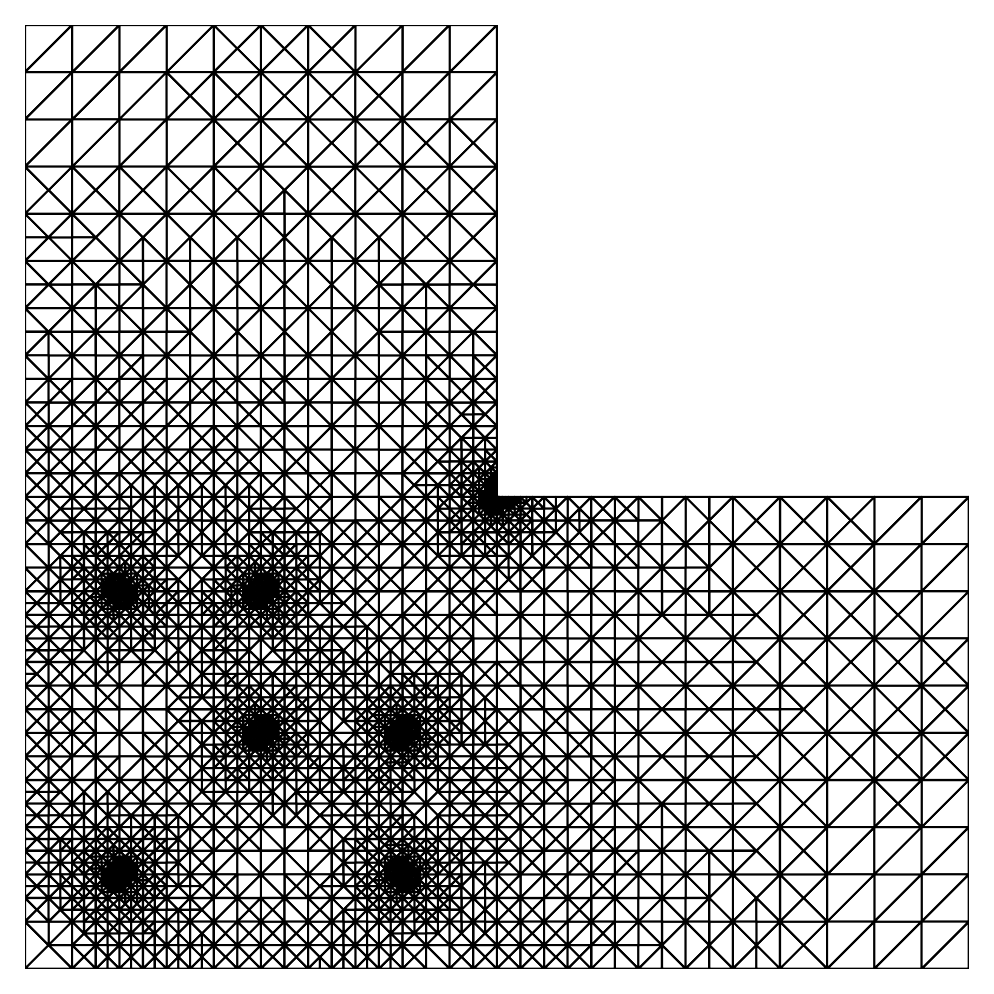}}
\subfigure[Error estimators]
{\includegraphics[width=0.37\textwidth]{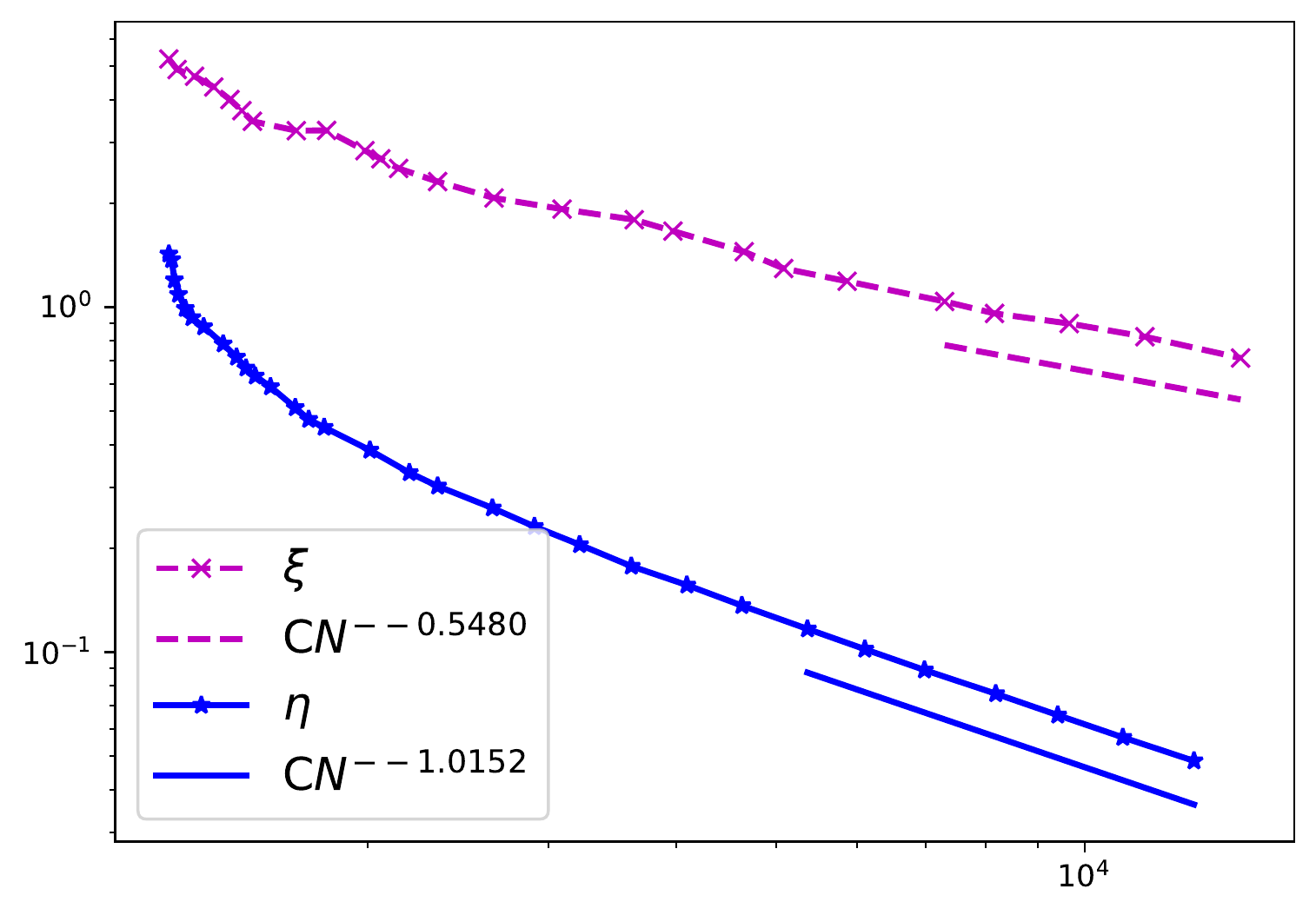}}
\vspace{-2mm}
\caption{Example \ref{example2}: Adaptive meshes and error estimators based on $P_2$ polynomials.}\label{fig: exam2 P2}
\end{figure}
\end{example}

\begin{figure}
\centering
\subfigure[Geometry 1]
{\includegraphics[width=0.3\textwidth]{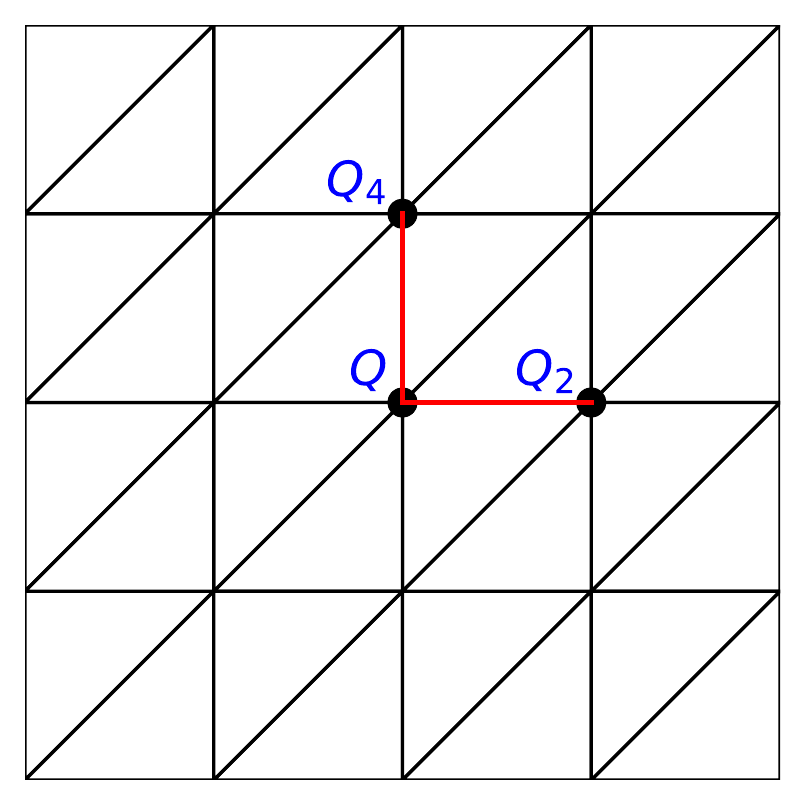}}\hspace{2mm}
\subfigure[Geometry 2]
{\includegraphics[width=0.3\textwidth]{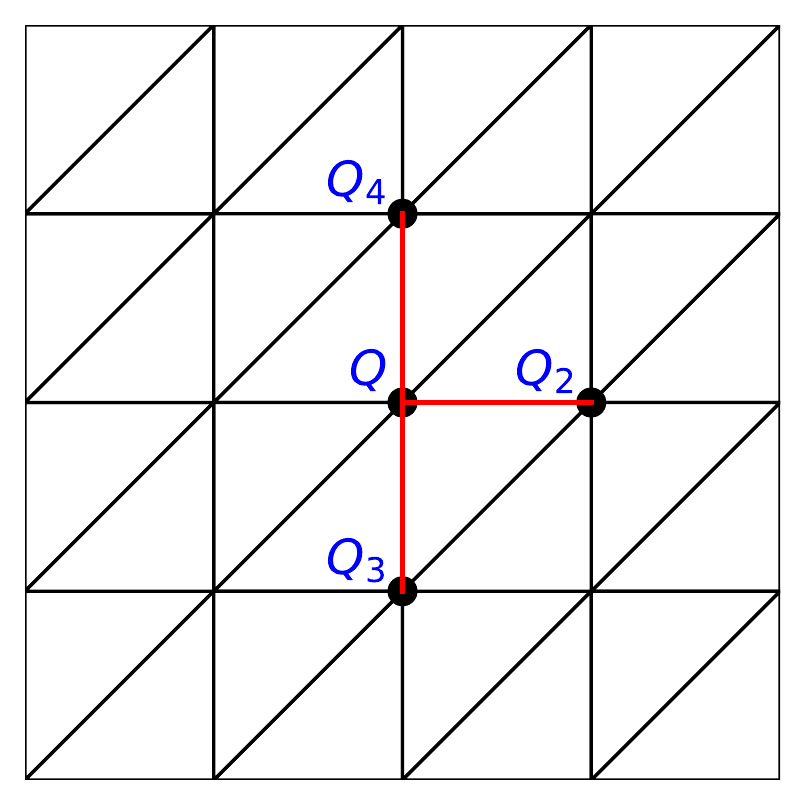}}\hspace{2mm}
\subfigure[Geometry 3]
{\includegraphics[width=0.3\textwidth]{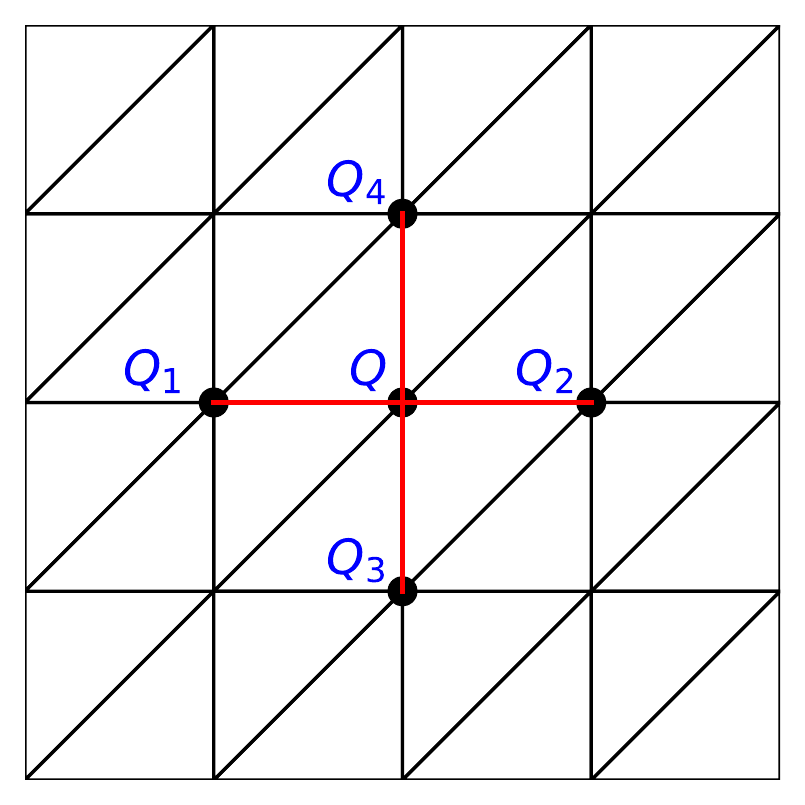}}\vspace{-2mm}
\caption{Example \ref{example4}: initial meshes.}\label{fig:initialmeshes exam4}
\end{figure}

\begin{example}\label{example4}
In this example, we first introduce four intersecting line fractures $\gamma_l=QQ_l$, $l=1,\cdots, 4$, where $Q(0.5,0.5)$, $Q_1(0.25,0.5)$, $Q_2(0.75,0.5)$, $Q_3(0.5,0.25)$ and $Q_4(0.5, 0.75)$. 
Here, we consider three types of geometries of $\Omega$. Geometry 1 consists of two line fractures $\gamma_2$ and $\gamma_4$; Geometry 2 consists of three line fractures $\gamma_2$, $\gamma_3$ and $\gamma_4$; Geometry 3 consists of all line fractures $\gamma_l$, $l=1,\cdots, 4$. The initial meshes of Geometry 1$-$3 are shown in Figure \ref{fig:initialmeshes exam4}.
The functions $g_l$ on each line fracture $\gamma_l$ are taken as the following,
\bq
\ba
g_1 = -g_2 = -g_3 =g_4 =-1.
\ea
\eq
\end{example}
\begin{figure}
\centering
\subfigure[Geometry 1]
{\includegraphics[width=0.32\textwidth]{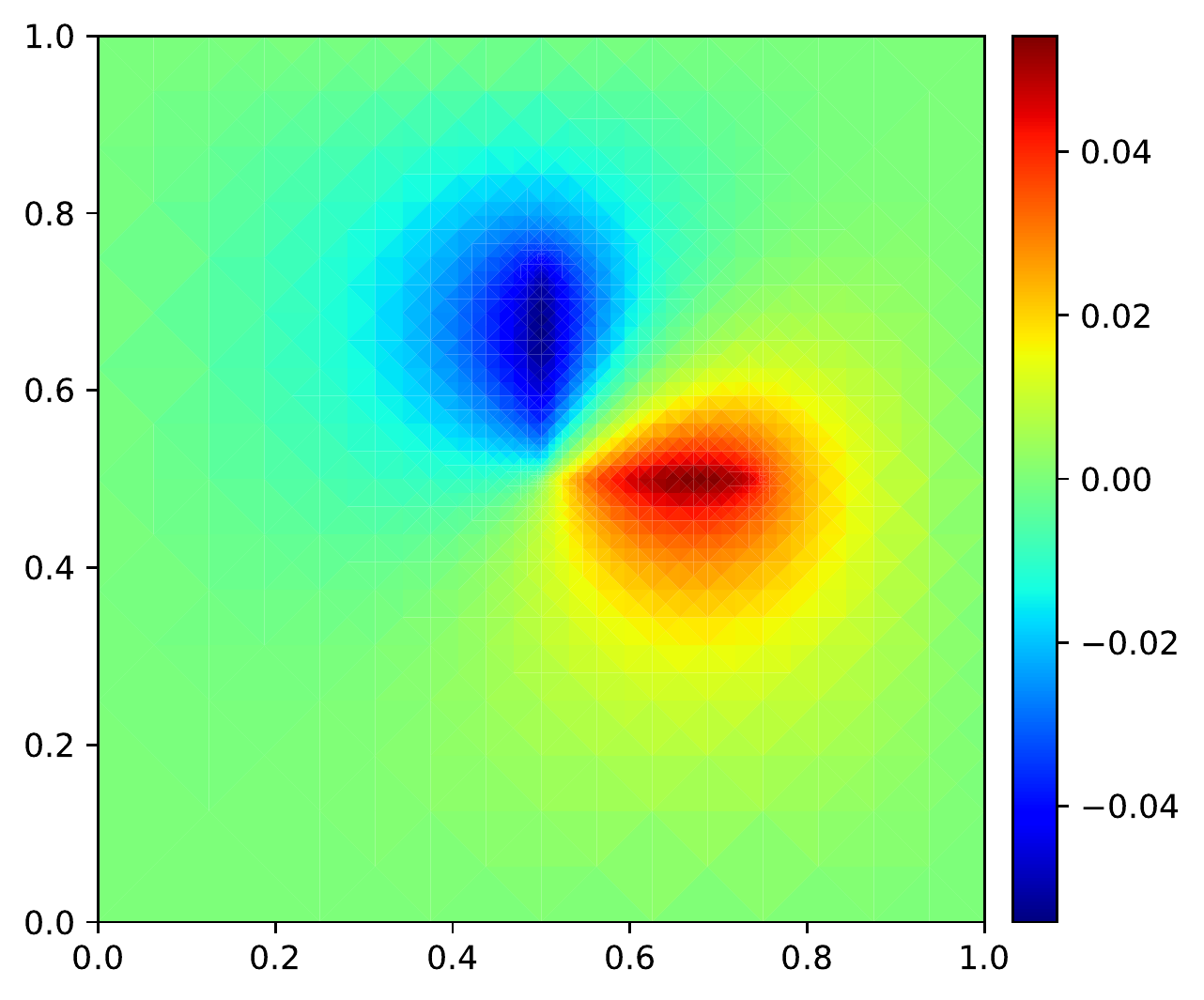}}
\subfigure[Geometry 2]
{\includegraphics[width=0.32\textwidth]{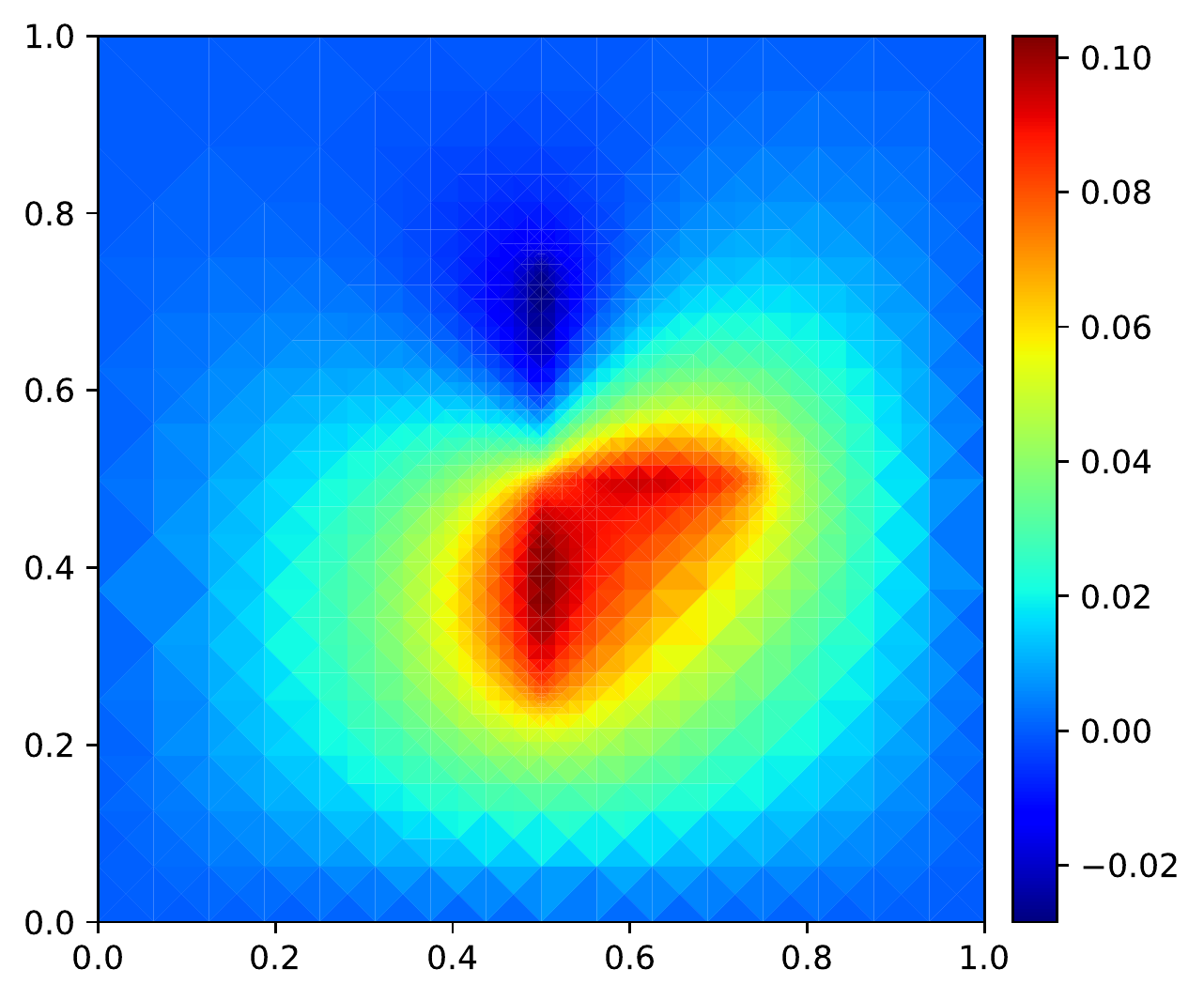}}
\subfigure[Geometry 3]
{\includegraphics[width=0.32\textwidth]{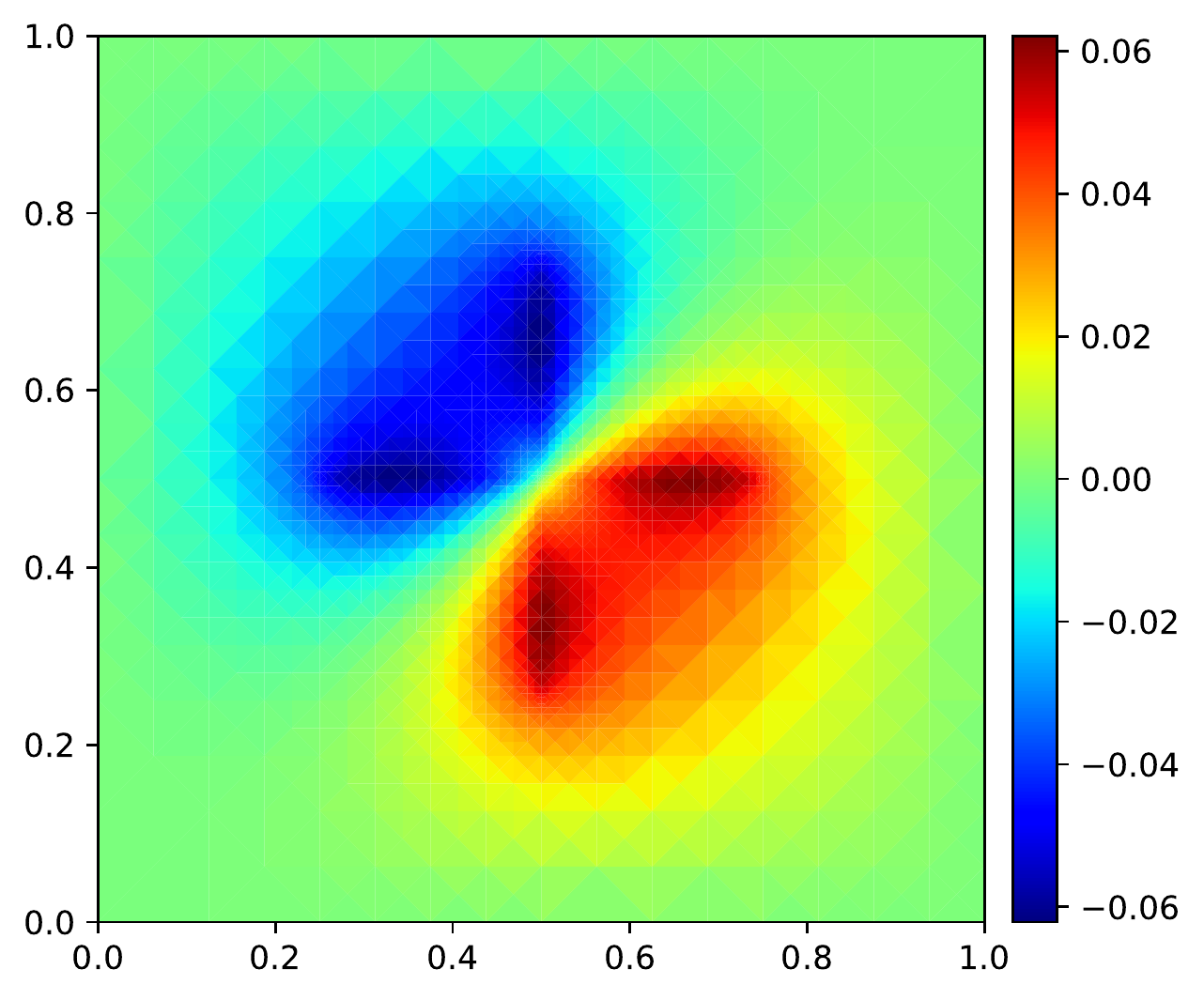}}
\vspace{-2mm}
\caption{Example \ref{example4}: AFEM solutions based on $P_1$ polynomials.}\label{fig:solution1 exam4}
\end{figure}

\begin{figure}
\centering
\subfigure[$P_1$]
{\includegraphics[width=0.37\textwidth]{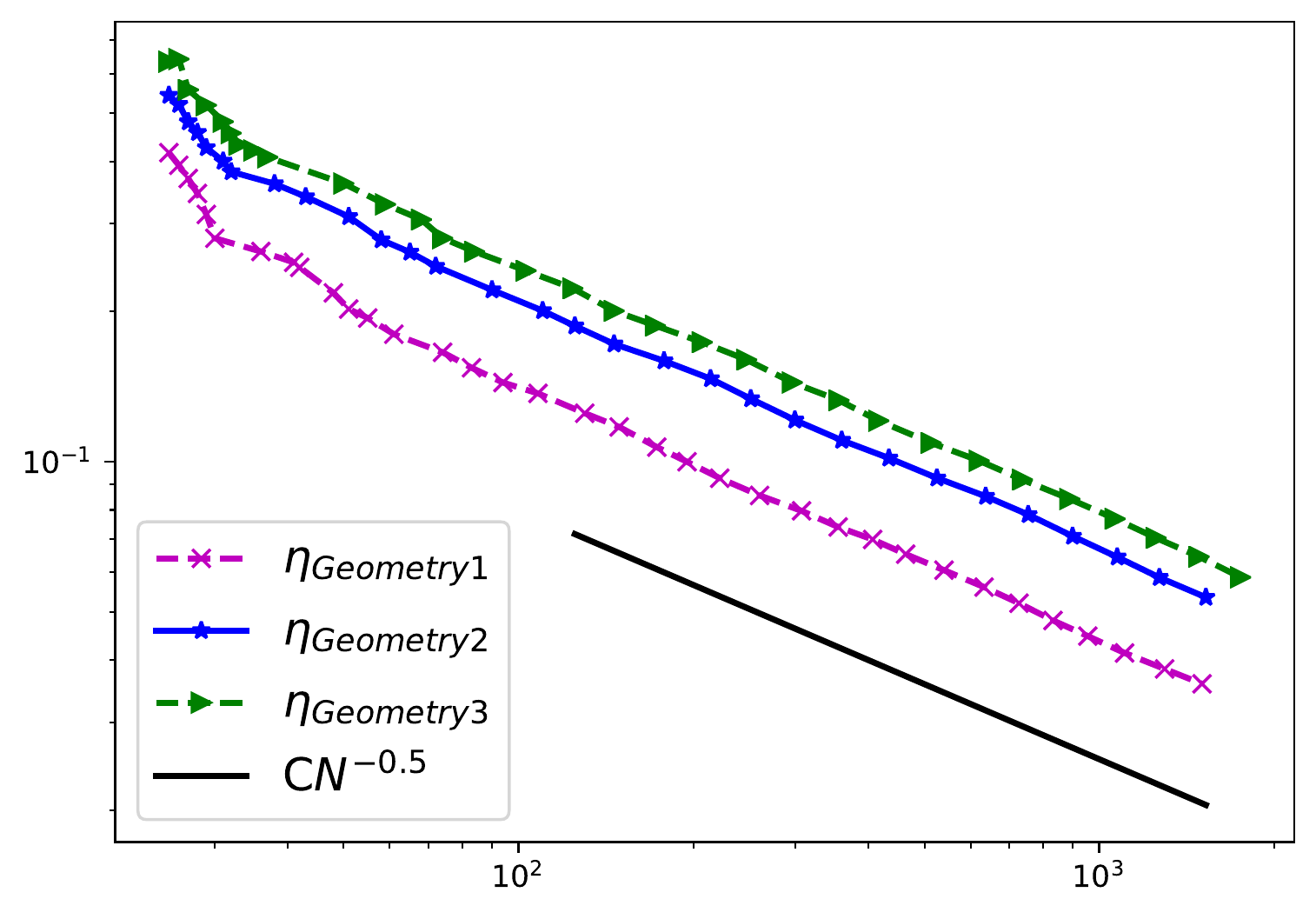}}\hspace{8mm}
\subfigure[$P_2$]
{\includegraphics[width=0.37\textwidth]{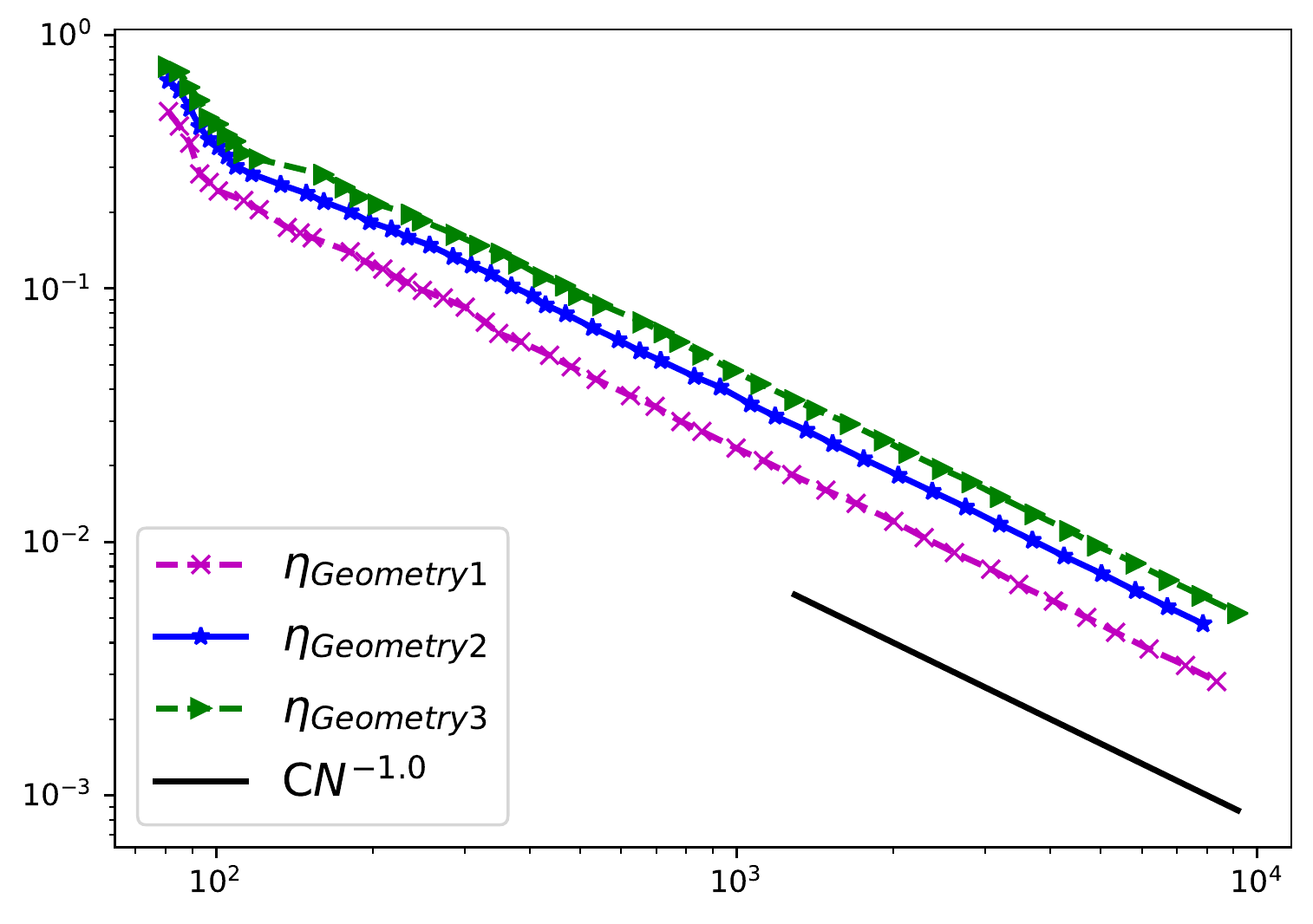}}
\vspace{-3mm}
\caption{Example \ref{example4}: error estimators.}\label{fig:errors exam4}
\end{figure}

\begin{figure}
\centering
\subfigure[Geometry 1]
{\includegraphics[width=0.28\textwidth]{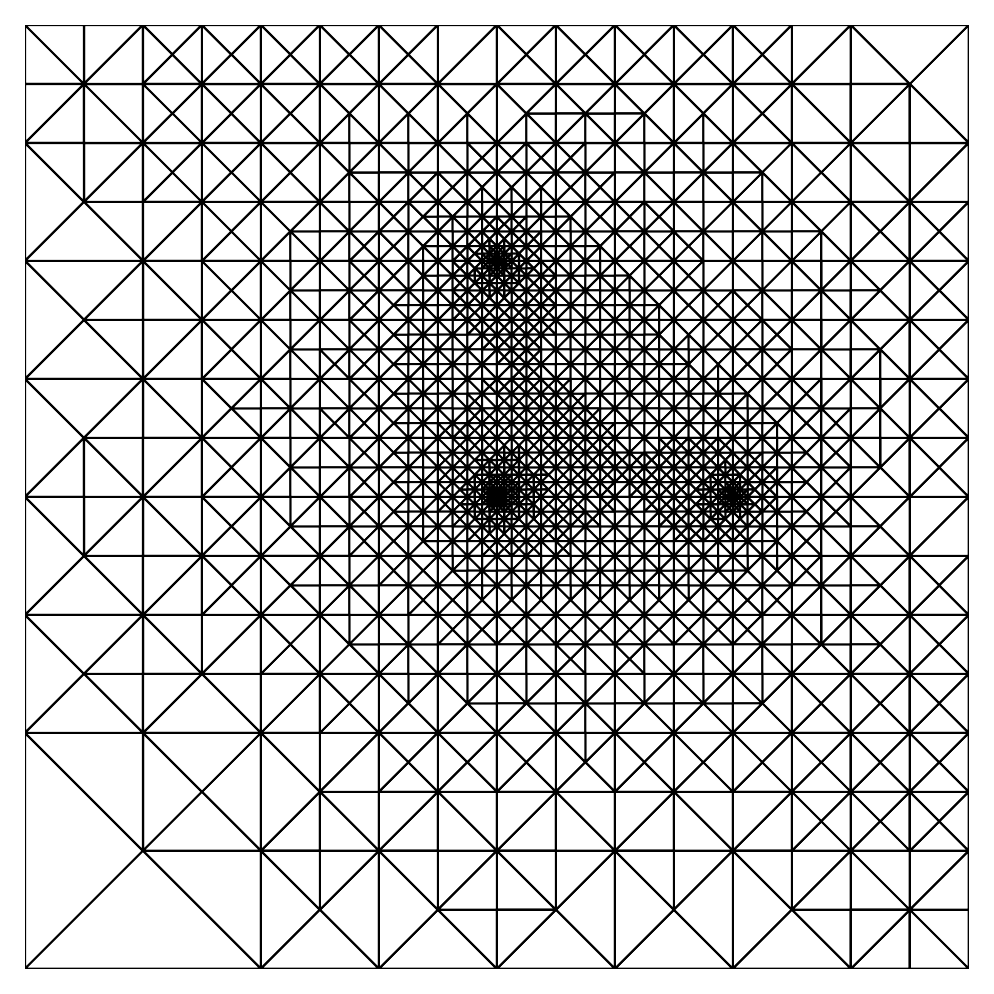}}\hspace{2mm}
\subfigure[Geometry 2]
{\includegraphics[width=0.28\textwidth]{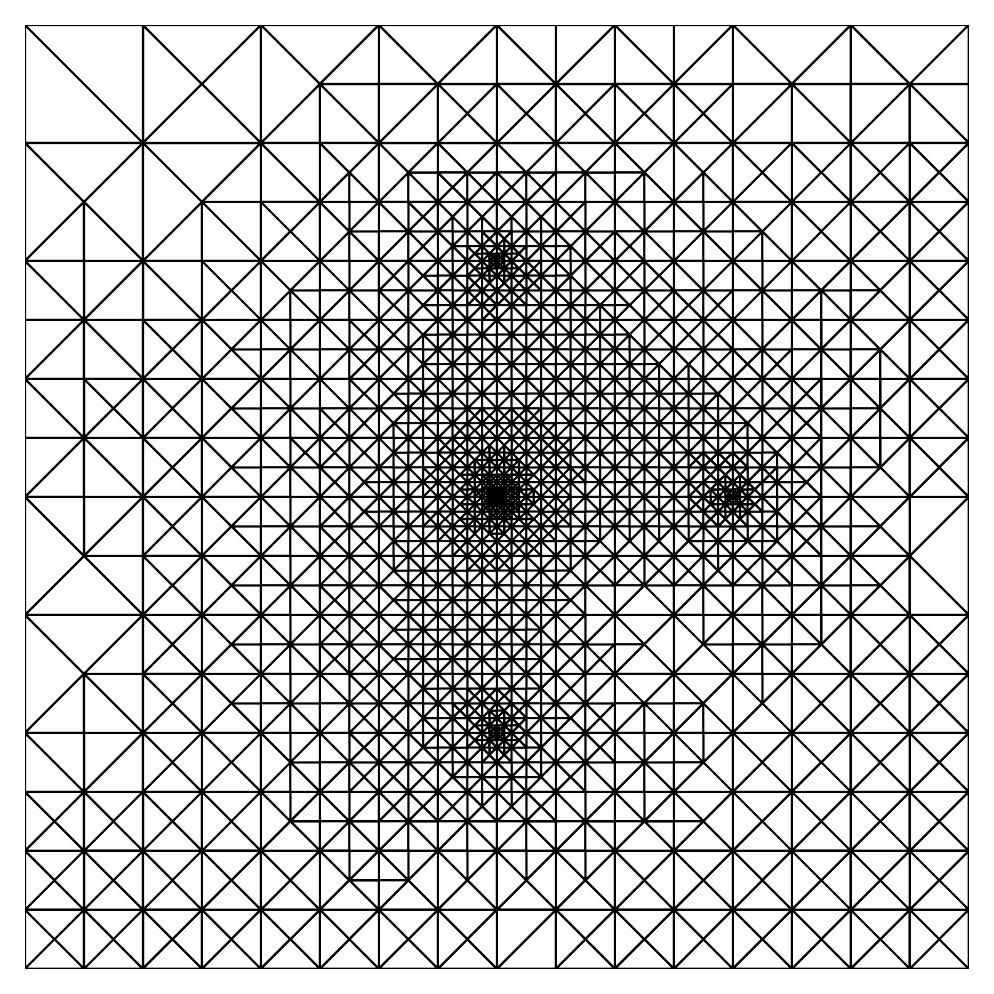}}\hspace{2mm}
\subfigure[Geometry 3]
{\includegraphics[width=0.28\textwidth]{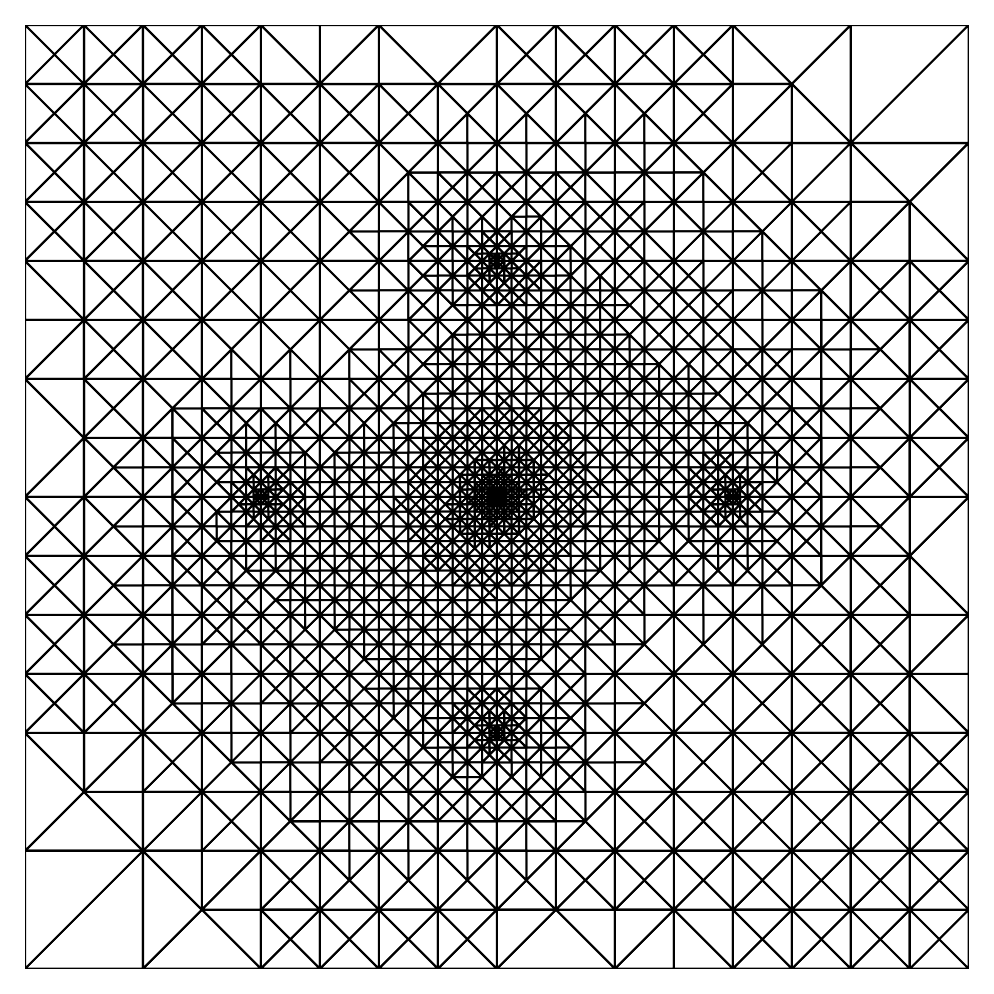}}
\vspace{-3mm}
\caption{Example \ref{example4}: adaptive meshes based on $P_1$ polynomials.}\label{fig:mesh1 exam4}
\end{figure}

\begin{figure}
\centering
\subfigure[Geometry 1]
{\includegraphics[width=0.28\textwidth]{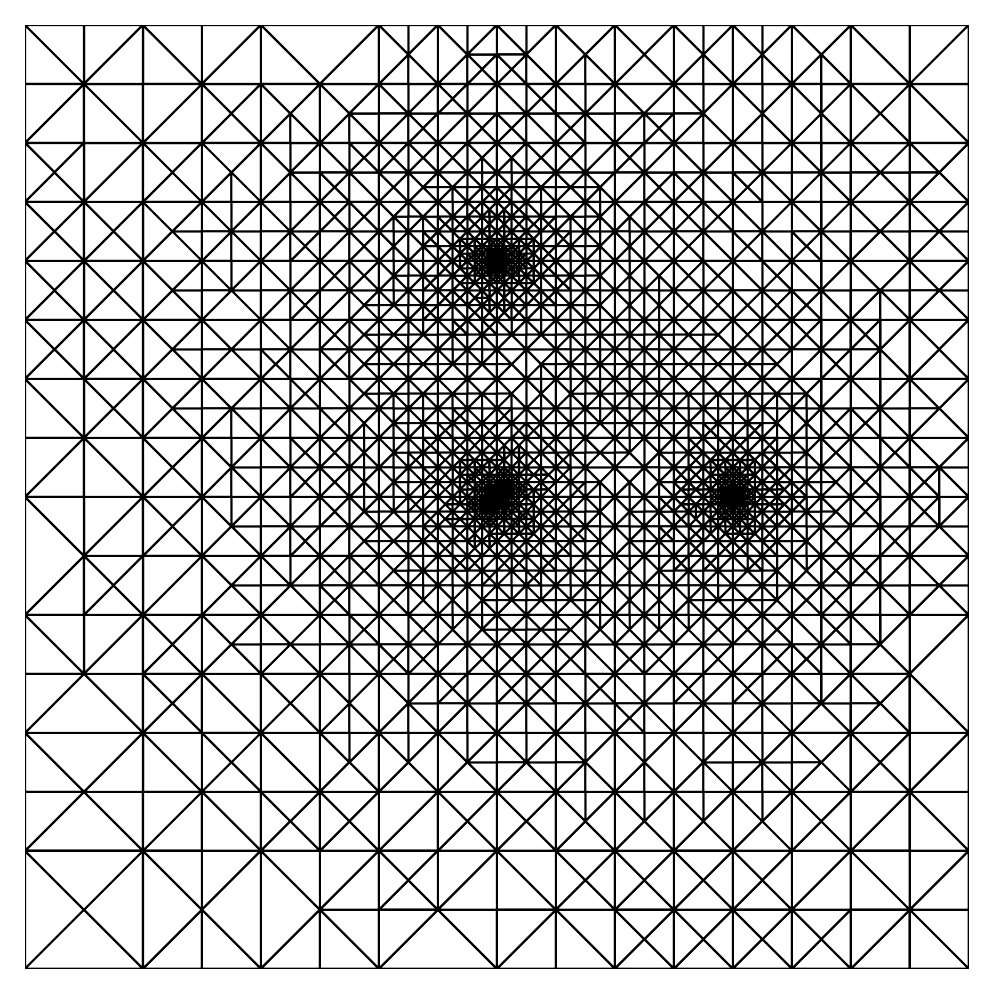}}\hspace{2mm}
\subfigure[Geometry 2]
{\includegraphics[width=0.28\textwidth]{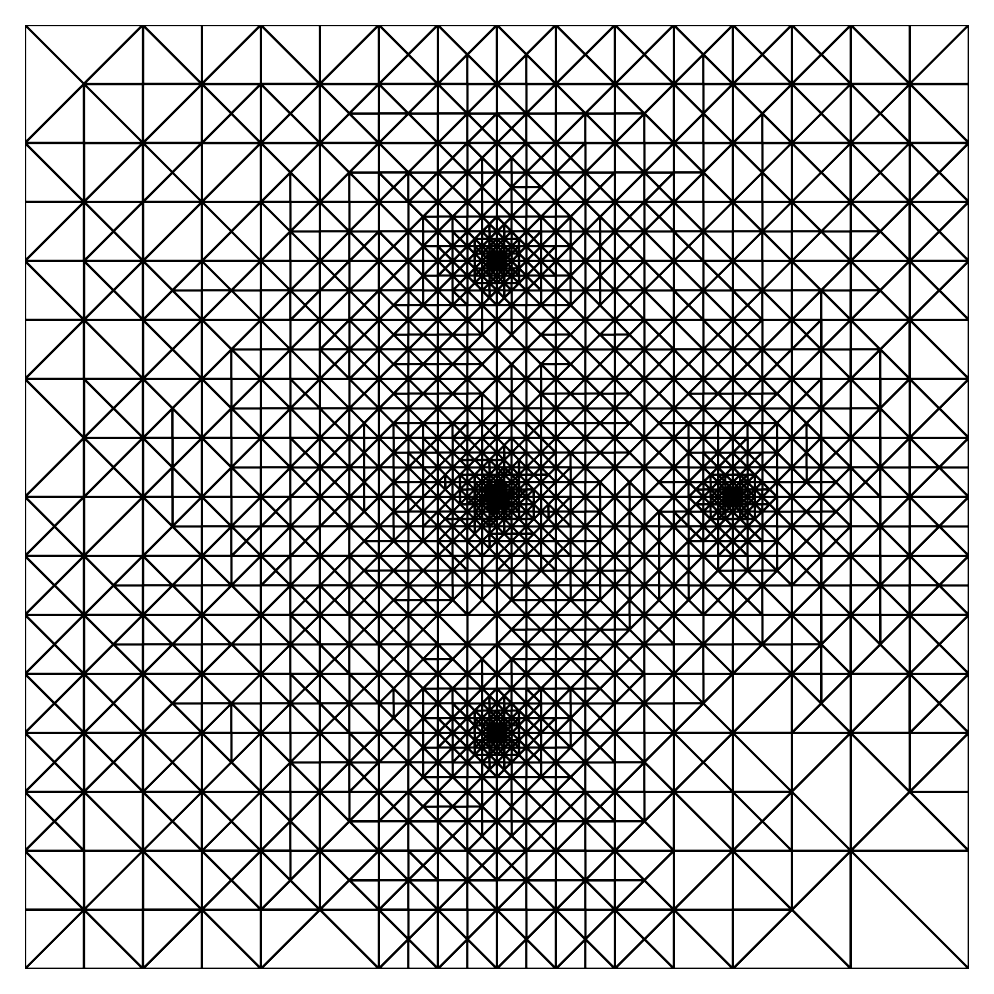}}\hspace{2mm}
\subfigure[Geometry 3]
{\includegraphics[width=0.28\textwidth]{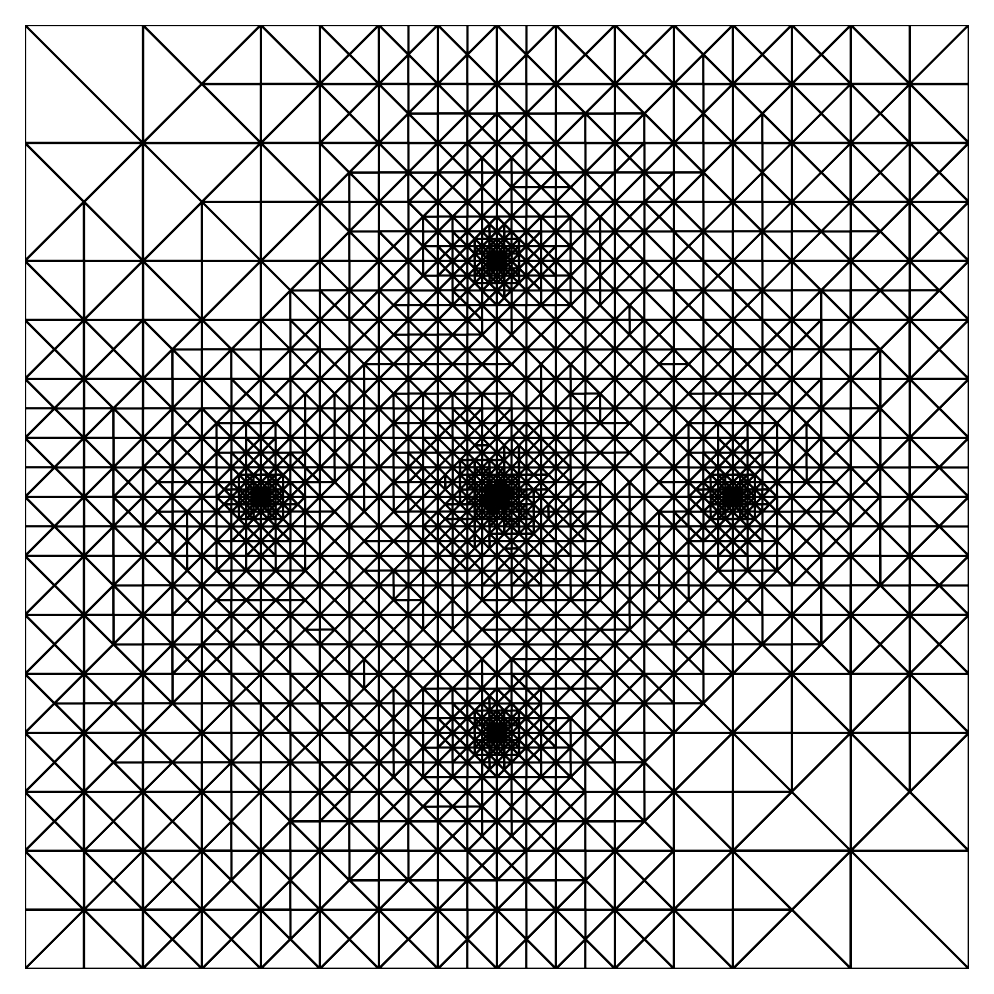}}
\vspace{-3mm}
\caption{Example \ref{example4}: adaptive meshes based on $P_2$ polynomials.}\label{fig:mesh2 exam4}
\end{figure}
The history of the error estimators are reported in Figure \ref{fig:errors exam4}, which shows that the convergence rates of the error estimators are quasi-optimal for all the three cases. 
Figure \ref{fig:mesh1 exam4}-\ref{fig:mesh2 exam4} and Figure \ref{fig:solution1 exam4} show the corresponding adaptive mesh refinements and the numerical solutions, respectively. We can see clearly that the error estimator successfully guide the mesh refinement around the singular points $Q_i$, where the solution shows singularity.

\section*{Acknowledgments}
H. Cao was supported by Hunan Provincial
Innovation Foundation for Postgraduate (CX20200619).
H. Li was supported in part by the National Science Foundation Grant DMS-1819041 and by the Wayne State University Faculty Competition for Postdoctoral Fellows Award.
N. Yi was supported by NSFC Project (12071400,1191410) and China's National Key R\&D Programs (2020YFA0713500).

\end{document}